\documentclass[11pt]{article}

\usepackage[
  shownumpages,               
  bgcolor={245,245,250},      
  braincolor={19,102,240},    
  citingstyle=authoryear,     
  bibliostyle=plainnat,       
  bibfile=refs                
]{brainlab}

\usepackage{microtype}
\usepackage{graphicx}
\usepackage{subfigure}
\usepackage{booktabs} 
\usepackage{hyperref}
\usepackage{wrapfig}
\usepackage{multirow}
\usepackage{cuted}
\usepackage{multicol}

\setbrainmeta{
  title={Markovian Compression: Looking to the Past Helps Accelerate the Future},
  authors={
    Andrey Veprikov\textsuperscript{1, 2}, Vladimir Solodkin\textsuperscript{1, 2}, Mikhail Rudakov\textsuperscript{2, 3}, Petr Babkin\textsuperscript{2}, Aleksandr Beznosikov\textsuperscript{1, 2}
  },
  affiliations={
    \textsuperscript{1}Federated Learning Problems Laboratory \\
    \textsuperscript{2}Basic Research of Artificial Intelligence Laboratory (BRAIn Lab) \\
    \textsuperscript{3}Innopolis University
  },
  abstract={
    This paper deals with distributed optimization problems that use compressed communication to achieve efficient performance and mitigate communication bottleneck. We propose a family of compression schemes in which operators transform vectors fed to their input according to a Markov chain, i.e. the stochasticity of the compressors depends on previous iterations. The compressors are implemented in the vanilla Quantized Stochastic Gradient Descent (\texttt{QSGD}) algorithm, and, to further improve the efficiency and convergence rate, in the momentum accelerated \texttt{QSGD}. We provide convergence results for our algorithms with Markovian compressors, the analysis covers non-convex, Polyak-Lojasiewicz, and strongly convex cases. To demonstrate the applicability of our approach to distributed data-parallel optimization problems, we conduct experiments on the CIFAR-10 and GLUE datasets with the Resnet-18 and DeBERTaV3 models. Practical results show the superiority of methods that use our compressor design over existing schemes.
  },
}

\usepackage{url}            
\usepackage{nicefrac}       
\usepackage{xcolor}         
\usepackage{amsmath,amsfonts,amssymb,amsthm,array}
\usepackage[font=small]{caption}
\usepackage{subcaption}
\usepackage{bm}
\usepackage{color}
\usepackage{enumitem}
\usepackage{bbm}
\usepackage{relsize}
\usepackage{makecell} 
\allowdisplaybreaks

\newtheorem{lemma}[theorem]{Lemma}
\newtheorem{corollary}[theorem]{Corollary}
\theoremstyle{remark}

\newtheorem{example}[theorem]{Example}

\graphicspath{ {./Figures/} }

\usepackage{amsfonts}
\usepackage{amsmath}
\usepackage{amsthm}
\usepackage{amssymb}
\usepackage{dsfont}

\usepackage{xcolor}
\usepackage{colortbl}
\usepackage{color}
\usepackage{graphicx}
\usepackage{multirow}
\usepackage{pifont}

\usepackage{mathtools}

\usepackage{verbatim}
\usepackage{xspace} %

\newcommand{\dotprod}[2]{\left\langle #1,#2 \right\rangle}
\newcommand{\norms}[1]{\left\| #1 \right\|}
\newcommand{\expect}[1]{\mathbb{E}\left[ #1 \right]}
\newcommand{\prob}[1]{\mathbb{P}\left\{ #1 \right\}}

\definecolor{PineGreen}{HTML}{008B72}

\providecommand{\norm}[1]{\left\lVert#1\right\rVert}

  \providecommand{\R}{\mathbb{R}} %

  \DeclareMathOperator{\E}{{\mathbb E}}
  \providecommand{\prob}[1]{{\rm Pr}\left[#1\right] }

  \providecommand{\e}{\varepsilon}

  \usepackage{bm}

\newcommand{\circledOne}{\text{\ding{172}}}
\newcommand{\circledTwo}{\text{\ding{173}}}
\newcommand{\circledThree}{\text{\ding{174}}}
\newcommand{\circledFour}{\text{\ding{175}}}

\usepackage{url}

\def\<#1,#2>{\langle #1,#2\rangle}

\providecommand{\mycomment}[3]{\todo[inline, caption={},size=footnotesize,color=#1!20]{\textbf{#2: }#3}}%
\newcommand\commenter[2]%
{%
  \expandafter\newcommand\csname #1\endcsname[1]{\mycomment{#2}{#1}{##1}}
}

\commenter{Vladimir}{red} 
\commenter{andrew}{blue} 
\begin{document}
\begin{mainpart}
\addtocontents{toc}{\protect\setcounter{tocdepth}{0}}

\section{Introduction}

The optimization problem is currently a key issue in many practical applications, such as optimization in neural network training, resource allocation in computational systems, and parameter tuning in algorithmic trading strategies.

Consequently, a variety of algorithms for optimization on a single device, including \texttt{SGD} \cite{robbins1951stochastic}, \texttt{Adam} \cite{kingma2014adam} and \texttt{Lion} \cite{yazdani2016lion}, have emerged and have been subjected to theoretical analysis.
However, in the contemporary landscape of deep learning, there is an increasing trend towards adopting complex and expansive models that pose significant training challenges.
These challenges are evident in advanced deep learning frameworks for image analysis \cite{dosovitskiy2021imageworth16x16words}, sophisticated natural language processing structures similar to transformers \cite{vaswani2017attention}, and complex reinforcement learning methodologies designed for autonomous system operations \cite{kiran2021deep}.
As a result, the training of such models has become impractical for implementation on a single device due to their requirement for massive datasets for training, which are unfeasible to store at one place.
Consequently, optimization algorithms have been specifically developed for distributed learning \cite{verbraeken2020survey, chen2021distributed}.
Such type of methods utilize a large number of devices, with each device processing distinct data subsets and participating in an effective data exchange mechanism. This collaborative effort facilitates the training of computationally intensive models.
Thus, the problem of classical optimization evolves into a distributed optimization form:
\begin{equation}
    \label{eq:problem}
        \min_{x \in \mathbb{R}^d} \Big\{ f(x) := \frac{1}{n} \sum_{i=1}^n f_i(x) \Big\} ,
\end{equation}
where $f_i$ is a function, located on a device $i$.
Note that this formulation encompasses two concepts. The first is distributed learning, in which data are spread across multiple devices to facilitate the storage of large amounts of data. The second is federated learning \cite{konevcny2016federated, li2020review, kairouz2021advances}, in which data distribution is motivated by the architecture of the system itself, allowing decentralized model training while maintaining data privacy and integrity between different devices.

The downside of this approach is the complexity associated with the transmission of large-scale data, resulting in "communication bottleneck" \cite{gupta2021localnewton}. This bottleneck can significantly impede the efficiency of the system, particularly in scenarios involving extensive data exchange across a distributed network. Furthermore, the challenge is amplified in environments with limited bandwidth, requiring solutions to mitigate the impact of data transmission delays and ensure seamless data flow.

The primary solution at present is the compression of the transmitted information \cite{bekkerman2011scaling, chilimbi2014project, alistarh2017qsgd}. That is, no entire package of data is sent, but rather a selected subset. The idea involves strategically selecting and compressing the most informative data segments for transmission. Consequently, this results in the significant advantage of reducing the volume of data that must be communicated across the network, thereby alleviating the communication bottleneck.

Motivated by the aforementioned, a number of methods employing compression have been developed and extensively analyzed \cite{mishchenko2019distributed, gorbunov2021marina, richtarik2021ef21}. However, the majority of studies have focused on unbiased compression operators due to their simplicity and amenability to theoretical analysis. These compression techniques, including random sparsification and value rounding \cite{doi:10.1137/100802001, alistarh2017qsgd, horvath2022natural, beznosikov2023biased}, do not in any way take into account the information transmitted in previous iterations. We hence highlight a potential research gap regarding the usage of previously transferred data in compression operators and optimization algorithms.

This omission raises the following research questions that we address in our paper:
\begin{itemize}
    \item \textit{Is it possible to design compression operators that take into account information about what and how we forwarded in previous iterations?}

    \item \textit{What methods can we integrate this kind of compression operators into? How does it affect the convergence rate of the methods, both in theory and in practice?}

    \item \textit{Can the methods be made even more efficient, e.g., by using additional momentum acceleration techniques?}
\end{itemize}

In this paper, we focus on compression-based methods that take into account information that has already been transferred, employing Markovian compression operators.

\subsection{Our contributions}
\begin{itemize}
    \item \textbf{New type of compression operator.} We introduce a novel type of compressor that utilizes stochasticity transmitted over previous iterations. We refer to this type of compressors as Markovian because the states of these compressors can be viewed as a Markov chain. Moreover, we provide two naive examples of such compression operators (Definitions \ref{def:banlastk} and \ref{def:kawasaki}) that have been shown to be efficient and simple to implement in real-life problems.

\item \textbf{In-depth analysis of strongly convex and non-convex cases.}
Motivated by various applications primarily from machine learning, we study QSDG and Accelerated QSGD with Markovian compression and provide in-depth discussion of the convergence rate (see Section \ref{sec:discussion}). Moreover, we provide theoretical analysis in both the strongly convex (see Theorem \ref{theorem:GD_acc}) and the non-convex\,/\,PL-condition (see Theorem \ref{theorem:GD_odd}) case of the target function $f$.

\item \textbf{Experimental Validation.}
We conduct experiments with Markovian compressors in a data-parallel setup for a series of optimization problems and datasets. In more details, we study the \texttt{QSGD}, \texttt{Accelerated QSGD} and \texttt{DIANA} for logistic regression on MNIST and LIBSVM datasets, and  \texttt{SGD} for image classification with Resnet-18 on CIFAR-10. Moreover, we consider \texttt{Adam} for fine-tuning DeBERTaV3-base for GLUE benchmark. In all setups, we observe an acceleration of convergence for methods employing our compressors compared to the baselines.
\end{itemize}

\subsection{Related work}

\textbf{Compressed communications.}
The use of compressed communications is a fairly well-known idea in distributed learning \cite{Seide20141bitSG}. As soon as the main property of compressed messages is that they are much easier to transfer, it can be reached in different ways, such as by quantizing the entries of the input vector \cite{alistarh2017qsgd, mayekar2019ratq, gandikota2020vqsgd, horvath2022natural}, or by sparsifying it \cite{richtarik2016parallel, alistarh2018convergence}, or even by combining these ideas \cite{albasyoni2020optimal, beznosikov2023biased}. However, all of the compression operators could be roughly \cite{condat2023efbv} separated into two large groups: \textit{unbiased} and \textit{biased}.

The first group is much easier to analyze and is therefore more broadly represented in the literature. The basic method with unbiased compression was presented in \cite{alistarh2017qsgd}. Later this algorithms were modified using variance reduction technique with compression of gradient differences \cite{mishchenko2019distributed, horváth2019stochastic, gorbunov2021marina} in order to improve the theoretical convergence guarantees. One can also note the works \cite{gorbunov2019unified} and \cite{khaled2020tighter}, where the authors developed a general theory for \texttt{SGD}-type methods with unbiased compression.

On the other hand, our understanding of distributed optimization with biased compressors is more complicated. In particular, biased compression implies the use of error compensation techniques \cite{stich2018sparsified}. Distributed \texttt{SGD} with biased compression and linear rate of convergence in a multi-node setting was first introduced in \cite{beznosikov2023biased}. In the meantime, other error compensation techniques are being actively developed, \cite{lin2022differentially, richtarik2021ef21}. The last approach called \texttt{EF21} was later studied in \cite{fatkhullin2021ef21}, \cite{gruntkowska2023ef21p}.

\textbf{Markovian stochasticity.}
Another recent trend in the literature is to design algorithms that use Markovian stochastic processes instead of $i.i.d.$ random variables in various ways. For instance, \cite{duchi2012ergodic} introduced a version of the \texttt{Mirror Descent} algorithm that yields optimal convergence rates for non-smooth and convex problems. Later, \cite{doan2020convergence, dorfman2023adapting, beznosikov2023first} studied first-order methods in the Markovian noise setting. Alternatively, token algorithms \cite{hendrikx2022principled,ayache2022walk} are also a popular area of research in Markovian stochasticity. In particular, \cite{even2023stochastic} obtained optimal rates of convergence, and \cite{pmlr-v162-sun22b, mao2019walkman,doan2020finitetime} looked at the token algorithm from the angle of the Lagrangian duality and from variants of the \texttt{ADMM} method. At the same time, there exist the results, e.g., \cite{bresler2020squares}, which provide a lower bound for the particular finite sum problems in the Markovian setting.

Despite all of the above, to the best of our knowledge, there are currently no works that combine compressed communications and Markovian stochasticity of the compressors.

\subsection{Technical preliminaries}

    \textbf{Notations.}~~ We use $\dotprod{x}{y} := \sum_{i = 1}^d x_i y_i$ to denote standard inner product of vectors $x, y \in \mathbb{R}^d$ and $\left( x \odot y \right)_i = x_i y_i$ to denote Hadamard product of vectors $x, y \in \mathbb{R}^d$. We introduce $l_2$-norm of vector $x \in \mathbb{R}^d$ as $\norm{x} := \sqrt{\dotprod{x}{x}}$. We define $x^* \in \mathbb{R}^d$ as a point, where we reach the minimum in the problem \eqref{eq:problem}. We also denote $f^* > - \infty$ as a global (potentially not unique) minimum of $f$. We use a standard notation for $(d-1)$-dimensional simplex $\Delta_d := \textstyle{\left\{p \in \mathbb{R}^d ~|~ p_j \geq 0 ~\text{ and }~ \sum_{j=1}^d p_j = 1 \right\}}$ and for a set of natural numbers $\overline{1, n} := \{ 1, 2, \ldots, n \}$. 
    We denote $C_m^k$ as the binomial coefficient $ \binom{m}{k}$.

    Throughout the paper, we assume that the objective functions $f_i$ and the function $f$ from \eqref{eq:problem} satisfy the following assumptions.

    \begin{assumption}[$L_i$-smooth]
    \label{as:lip}
        Every function $f_i$ is $L_i$-smooth on $\mathbb{R}^d$ with $L_i > 0$, i.e. it is differentiable and there exists a constant $L_i > 0$ such that for all $x, y \in \mathbb{R}^d$ it holds that $
            \norms{\nabla f_i (x) - \nabla f_i (y)}^2 \leq L_i^2 \norm{x - y}^2 .
        $
        We define $L^2 := \frac{1}{n}\sum_{i=1}^nL_i^2$.
    \end{assumption}

    \begin{assumption}[$\mu$-strongly convex]
    \label{as:strconv}
        The function $f$ is $\mu$-strongly convex on $\mathbb{R}^d$, i.e., it is differentiable and there is a constant $\mu > 0$ such that for all $x, y \in \mathbb{R}^d$ it holds that $
            (\mu / 2)\norm{x - y}^2 \leq f(x) - f(y) - \dotprod{\nabla f(y)}{x - y} .
        $
    \end{assumption}

    \begin{assumption}[PL-condition]
        \label{as:PL}
            The function $f$ satisfies the PL-condition, i.e., it is differentiable and there is a constant $\mu > 0$ such that for all $x \in \mathbb{R}^d$ it holds that $
                \norm{\nabla f(x)}^2 \geq 2 \mu \left( f(x) - f^* \right) .
            $
    \end{assumption}

    \begin{assumption}[Data similarity]
    \label{as:sim}
        The functions $f_i$ are similar on $\mathbb{R}^d$, i.e., there are constants $\delta, \sigma \geq 0$, such that the following inequality holds for all $x \in \mathbb{R}^d$: $
            \norm{\nabla f_i(x) - \nabla f(x)}^2 \leq \delta^2 \norm{\nabla f(x)}^2 + \sigma^2 .
        $
    \end{assumption}

    The last inequality implies that the data stored at each device does not differ significantly. This Assumption is quite standard in the literature \cite{shamir2014communication, arjevani2015communication, khaled2020tighter, woodworth2020minibatch, gorbunov2021local, beznosikov2022decentralized, beznosikov2023first}.

    Now we introduce important definitions related to the theory of Markov processes.

    \begin{definition}[Markov chain]
        Markov chain with a finite state space $\{\nu_n\}_{n=0}^N$ is a stochastic process $\{ X_t \}_{t\geq0}$, that satisfies Markov property, i.e. $\mathbb{P} \{X_t = \nu_t ~|~ X_{t-1} = \nu_{t-1}, X_{t-2} = \nu_{t-2}, ...,  X_{0} = \nu_0\} = \mathbb{P} \{X_t = \nu_t ~|~ X_{t-1} = \nu_{t-1}\}$.
    \end{definition}

    \begin{definition}[Ergodicity of Markov chain]
    \label{def:tmp_def}
        Markov chain $\{ X_t \}_{t\geq0}$ with a finite state space $\{\nu_n\}_{n=0}^N$ is referred to be ergodic if for any $n \in \overline{1, N}$ there exists $\lim\limits_{t \to \infty} \prob{X_t = \nu_n ~|~ X_0 = \nu_0} = p_n$,
        where $0 \leq p_n \leq 1$ does not depend on the $\nu_0$. 
        If Markov chain is ergodic, then $\left\{ p_n \right\}_{n=0}^N \in \Delta_N$ and there exist $0 < \rho < 1, C > 0$, such that 
        $
            \left| \prob{X_t = \nu_n ~|~ X_0 = \nu_0} - p_n \right| \leq C \rho^{t} .
        $
    \end{definition}

    \begin{definition}[Mixing time of the discrete Markov chain]
        We say that $\tau_{\text{mix}}(\varepsilon)$ is the mixing time of the ergodic Markov chain $\{ X_t \}_{t\geq0}$ with stationary distribution $\{p_n\}_{n=0}^N$, if $\forall \varepsilon > 0, \,\forall t \geq \tau_{\text{mix}}(\varepsilon) \hookrightarrow \max\limits_{n \in \overline{0, N}}\{\| \mathbb{P}\{X_t = \nu_n ~|~ X_0 = \nu_0\} - p_n \| \} \leq \varepsilon \cdot p_{\min}$, where $p_{\min} := \min_{n \in \overline{0, N}}\{ p_n \}$. From Definition \ref{def:tmp_def}, it follows that $\tau_{\text{mix}}(\varepsilon) \geq \frac{\log(C / p_{\min} \varepsilon)}{\log(1/\rho)}.$
    \end{definition}

    These definitions are extremely important for further analysis of the Markovian compressors, which are presented in the next section.

\section{Main results}

\subsection{Markovian compressors} \label{subsec:compressors}

    In this section, we introduce Markovian compressors that take into account the information transmitted in previous $K$ iterations. It is assumed that these compressors function within an iterative algorithm aimed at minimizing the problem \eqref{eq:problem}, wherein a discrete variable $t$ denotes time step. Consequently, due to the dependence of the compressors on previous states, they exhibit a reliance on the step $t$. Let us narrow down the class of compressors to be discussed in this paper.

    \begin{definition}[Random sparsification]
    \label{def:rand_spars}

        $Q_t(x)$ is a random sparsification compressor, if it operates on the vector $x \in \mathbb{R}^d$ as $
            Q_t(x) = \frac{d}{m} x \odot \mathbbm{1}(\nu_t),$
        where $\nu_t$ is a set of $m$ coordinates: $\nu_t \subseteq \overline{1, d}$.
    \end{definition}

    The classical Rand$m$ operator fits Definition \ref{def:rand_spars}, in particular, subsets $\nu_t$ for this compressor are generated uniformly at each step $t$, therefore it is unbiased, i.e., $\mathbb{E}_{t}[Q_t(x)] = x$ for all $t$. In this paper, we do not generate $\nu_t$ independently, but according to some Markov chain, i.e., compressors start to take into account previous iterations. We formulate this idea as an assumption.

    \begin{assumption}[Asymptotic unbiasedness of Markovian compressors]
    \label{as:compressors}
        We assume that operator $Q_t$ is a random sparsification compressor (Definition \ref{def:rand_spars}) and $\left\{ \nu_t \right\}_{t \geq 0}$ are realizations of some ergodic Markov chain with uniform stationary distribution.
    \end{assumption}

    Assumption \ref{as:compressors} implies that in the limit as $t \to \infty$, the compressor $Q_t$ is unbiased, i.e., $\expect{Q_t(x)} \to x$ as $t \to \infty$, because the stationary distribution of the Markov chain is uniform. We are now ready to introduce two examples that adhere to Assumption \ref{as:compressors}.     
    The first compressor is called \texttt{BanLast}($K, m$), it prohibits sending coordinates that have been sent at least once in the last $K$ iterations.

    \begin{definition}[\texttt{BanLast}($K, m$) compressor]
    \label{def:banlastk}
        Let $Q_t(x)$ be a random sparsification compressor (Definition \ref{def:rand_spars}). The $j \in \nu_t$ are chosen according to the distribution $p^t \in \Delta_d$ and $p^t$ is  given by the formula: 
        \begin{equation*}
         p^t_j =
        \begin{cases}
            0, & \text{if }~ j \in \bigcup_{s = t-K}^{t-1}\nu_s, 
            \\ \frac{1}{d - K m}, & \text{otherwise}.
            \end{cases}
        \end{equation*}
    \end{definition}
    
    The \texttt{BanLast}($K, m$) compressor exhibits a limitation in its utility due to an application restriction: $d \geq (K+1) m$, since we need at least $m$ coordinates to have a non-zero probability at each step $t$. In order to avoid these limitations, we introduce a more flexible Markovian compressor \texttt{KAWASAKI}($K, b, \pi_\Delta, m$).
    \begin{definition}[\texttt{KAWASAKI}($K, b, \pi_\Delta, m$) compressor]
    \label{def:kawasaki}
        Let $Q_t(x)$ be a random sparsification compressor (Definition \ref{def:rand_spars}). The $j \in \nu_t$ are chosen according to the distribution $p^t \in \Delta_d$, which is given by the formula:
        \begin{equation*}
        \begin{split}
            &\widetilde{p}_j^{~t} = \frac{1/d}{b^{\text{\# of choices $j$ for the last $K$ iterations}}},~~~ j \in \overline{1, d};
            \\&p^t = \pi_\Delta \left( \widetilde{p}^{~t} \right),
        \end{split}
        \end{equation*}
        where $b > 1$ is a forgetting rate and $\pi_{\Delta} : \mathbb{R}^d \to \Delta_d$ is an activation function. 
    \end{definition}

    The \texttt{KAWASAKI}($K, b, \pi_\Delta, m$) compressor is now applicable for arbitrary values of $d \ge m$, and $K$. On top of that, it introduces two additional parameters in comparison with \texttt{BanLast}($K, m$), namely $b$ and $\pi_\Delta$. The parameter $b$ is responsible for the how strongly we penalize a coordinate if it was selected in previous iterations: the larger $b$ is, the less likely we are to select a coordinate in step $t$ if it was selected in steps from $t-K$ to $t-1$. The function $\pi_\Delta$ is required in order to obtain the probability vector $p^t$ from the vector $\widetilde{p}^{~t}$. The following examples illustrate potential selections for $\pi_\Delta$:
    \begin{equation*}
    \begin{split}
        (&\pi_{\Delta}(\widetilde{p}))_j = |\widetilde{p}_j| / \|\widetilde{p}\|_1,\\ \,\,&\pi_{\Delta}\left(\widetilde{p}\right) = \text{Softmax}\left( \widetilde{p} \right),\\
        &\pi_{\Delta}\left(\widetilde{p}\right) = \underset{p \in \Delta_d}{\arg\min} \{\|\widetilde{p} - p \|^2\} .
    \end{split}
    \end{equation*}
    We now provide an example where using the Markovian compressor \texttt{BanLast}($K, m$) speeds up the optimization process by a factor of three compared to the unbiased compressor Rand$m$.
    \begin{example}
    \label{example:1}
        Consider the \texttt{QSGD} algorithm (Algorithm \ref{alg:GD}), which solves the problem \eqref{eq:problem} in the case $n=1$, of the form $x^{t+1} = x^t - \gamma Q(\nabla f(x^t))$. Assume that at some step $t$ we observe gradient of the form $(1, 0, ... , 0)^T \in \mathbb{R}^d$. In the QSGD algorithm, we compress the gradient at each step, therefore, we do not always send the first coordinate to the server, i.e. we do not move from the point $x^t$.

        In the case of $m = 0.1 \cdot d$, i.e. we send 10\% of all coordinates at each step, if we use the \texttt{BanLast}($K, m$) compressor, then the mathematical expectation of the number of steps to leave the point $x^t$ is approximately $3.4$ in the case of $K = 7$. For Rand10\%, this number is equal to $10$, i.e. we speed up the optimization process by a factor of three. For arbitrary values of $d$ and $m$, the formula for calculating the number of steps to leave the point $x^t$ is provided in Appendix \ref{appendix:example_1}.
    \end{example}

    Moreover, in Appendix \ref{appendix:example_1}, we obtain the general results for an arbitrary value of $\alpha \in (0;1]$ with $d = \alpha \cdot m$. 
    For each fixed $\alpha$, we can find the optimal value of $K^*(\alpha)$. It turns out that empirically this dependence is close to a linear one of the form $K^*(\alpha) \approx 0.73 \cdot \alpha$. It is important to highlight that at this point, $K$ is not treated as a hyperparameter; rather, it can be selected automatically.


    We now present a theorem demonstrating that our Markovian compressors from Definitions \ref{def:banlastk} and \ref{def:kawasaki} satisfy the conditions outlined in Assumption \ref{as:compressors}.

    \begin{theorem}[Asymptotic unbiasedness of \texttt{BanLast}($K, m$) and \texttt{KAWASAKI}($K, b, \pi_\Delta, m$)]
    \label{theorem:compressors}
        Compressors from Definitions \ref{def:banlastk} and \ref{def:kawasaki} can be described using Markov chains with states $\left\{ \nu_1, \nu_2, ..., \nu_K \right\}_{\nu_1, ..., \nu_K \in M}$, where $M$ is the set of all subsets of $\overline{1, d}$ of size $m$. Moreover,

        $\bullet$~~
            \texttt{BanLast}($K, m$) (Definition \ref{def:banlastk}) is ergodic with a uniform stationary distribution, if $d > (K+1) m$. 
        
        $\bullet$~~
            If $d > (2K+1)m$, then for \texttt{BanLast}($K, m$) we get
            \begin{equation*}
                \rho = \sqrt{1 - \left(\frac{C_{d - 2Km}^{m}}{(C_{d - Km}^{m})^2}\right)^K} \text{ and } C = \rho^{-2}.
            \end{equation*}
        $\bullet$~~
            If for all permutations $\phi$ of the set $\overline{1, d}$ it holds that $\pi_\Delta \left(\phi\left( \widetilde{p} \right) \right) = \phi \left(\pi_\Delta\left( \widetilde{p} \right) \right)$, then \texttt{KAWASAKI}($K, b, \pi_\Delta, m$) (Definition \ref{def:kawasaki}) is ergodic with a uniform stationary distribution.

        $\bullet$~~ 
            %
        If $\left(\pi_{\Delta}\left(\widetilde{p}\right)\right)_j = |\widetilde{p}_j| / \|\widetilde{p}\|_1$, then 
        \begin{equation}
        \rho = 1 - \left[ db^K - m(b^K - 1) \right]^{-mK} \text{ and } C = \rho^{-1}.
        \end{equation}
    \end{theorem}
    
    The proof of Theorem \ref{theorem:compressors} is provided in Appendix \ref{appendix:comprssors}. The outcomes of Theorem \ref{theorem:compressors} hold significant importance for the subsequent investigation of algorithms aimed at solving problem \eqref{eq:problem} employing Markovian compressors. Note that the examples of activation functions $\pi_\Delta$ provided above satisfy the conditions of Theorem \ref{theorem:compressors}.

    \subsection{Distributed gradient descent with Markovian compressors}

    In this section, we propose a new algorithm \texttt{Markovian QSGD} (Algorithm \ref{alg:GD}). This algorithm is similar to the vanilla \texttt{QSGD} \cite{alistarh2017qsgd}, but in line \ref{line:QSGD_Q} of Algorithm \ref{alg:GD} we use Markovian compressor $Q_t^i$ that we introduced in Section \ref{subsec:compressors}. That is, $Q_t^i$ can be either \texttt{BanLast}($K, m)$ or \texttt{KAWASAKI}($K, b, \pi_\Delta, m$), or any other Markovian compressor.

    \begin{theorem}[Convergence of \texttt{MQSGD} (Algorithm \ref{alg:GD})]
    \label{theorem:GD_odd}
        Consider Assumptions \ref{as:lip}, \ref{as:sim} and  \ref{as:compressors}. Let the problem \eqref{eq:problem} be solved by Algorithm \ref{alg:GD}.

        $\bullet$~~
            For any $\varepsilon, \gamma > 0$, $T > \tau > \tau_{\text{mix}}(\varepsilon)$ satisfying conditions, described in Appendix \ref{appendix_subsec:full_GD_odd}, it holds that
            \begin{equation*}
                \expect{\norm{\nabla f(\widehat{x}^{T})}^2} 
                =
                \mathcal{O} \left( \frac{F_\tau}{\gamma T}
                +
                \frac{\gamma L \tau d^2}{m^2} \sigma^2 \right),
            \end{equation*}
            where $\widehat{x}^{T}$ is chosen uniformly from $\left\{x^t \right\}_{t=0}^{T}$.

        $\bullet$~~
        If $f$ additionally verifies the PL-condition (Assumption \ref{as:PL}), then for any $\varepsilon > 0$, $\gamma > 0$, $\tau > \tau_{\text{mix}}(\varepsilon)$ and $T > \tau$ satisfying conditions, described in Appendix \ref{appendix_subsec:full_GD_odd}, it holds that
        \begin{equation*}
            \begin{split}
                F_T
                &=
                \mathcal{O} \left(
                \left( 1 - \frac{\mu \gamma}{12} \right)^{T-\tau} F_\tau
                +
                \frac{\gamma d^2 L \tau}{\mu m^2}\sigma^2 
                \right) .
            \end{split}
        \end{equation*}
        Here we use the notations $F_t := \expect{f(x^t) - f(x^*)}$  and $F_\tau := \expect{f(x^\tau) - f(x^*)}$.
    \end{theorem}
    The proof of Theorem \ref{theorem:GD_odd} is provided in Appendix \ref{appendix_subsec:non_conv}, \ref{appendix_subsec:GD_PL}. If Assumption \ref{as:sim} does not hold, we observe different results, which are provided in the Appendix \ref{appendix:proof_GD_sem}.
    
    \begin{wrapfigure}[10]{r}{7cm}
    \vspace{-0.35cm}
    \begin{minipage}{7cm}
    \begin{algorithm}{\texttt{Markovian QSGD}}\label{alg:GD}
        \begin{algorithmic}[1]
            \State {\bfseries Input:} starting point $x^0 \in 
               \mathbb{R}^d$, \\step size $\gamma>0$, \\number of iterations $T$
               \For{$t=0$ {\bfseries to} $T$}
               \State Broadcast $x^t$ to all workers
                   \For{$i=1$ {\bfseries to} $n$ in parallel}
                   \State Set $g_i^t$ = $Q_t^i\left( \nabla f_i(x^t) \right)$ \label{line:QSGD_Q}
                   \State Send $g_i^t$ to the server
                   \EndFor
                \State Aggregate $g^t = \frac{1}{n} \sum\limits_{i=1}^n g_i^t$
                \State Update $x^{t+1} = x^t - \gamma g^t$
               \EndFor
        \end{algorithmic}
    \end{algorithm}
    \end{minipage}
    \end{wrapfigure}    
    Usually in convergence evaluations of various methods, expressions with the term of $F_0$, i.e., something that depends on the initial guess, arise as constants, but in Theorem \ref{theorem:GD_odd}, a term of the form $F_\tau$ appears. This can be explained by the fact that, in iterations from $t = 0$ to $t = \tau$, the Markov chain has not been stabilized yet, and the initial state can be taken as $t = \tau$.

    \textbf{Sketch proof of Theorem \ref{theorem:GD_odd}}. 
    Let us write out a descent lemma of the form:
    \begin{equation}
        \begin{split}
        \label{eq:tmp_main}
            \mathbb{E} \left[ {\norm{x^{t+1} - x^*}^2} \right]
            &=
            \expect{\norm{x^t - x^*}^2}
            -
            2 \expect{\gamma \dotprod{\nabla f(x^t)}{x^t - x^*}}
            \\&-
            \underbrace{\frac{2 \gamma}{n}\sum\limits_{i=1}^n \expect{\dotprod{Q_t^i(\nabla f(x^t)) - \nabla f_i(x^t)}{x^t - x^*}}}_{\circledOne}
            \\&+
            \gamma^2 \expect{\norm{\frac{1}{n} \sum\limits_{i = 1}^n Q_t^i(\nabla f_i(x^t))}^2} .
        \end{split}
    \end{equation}
    The expression $\circledOne$ in \eqref{eq:tmp_main} is zero if $Q_t^i$ are unbiased and independent from iteration $t$, because $\mathbb{E}\langle Q_t^i(\nabla f(x^t)) - \nabla f_i(x^t), x^t - x^*\rangle = \mathbb{E} \langle \mathbb{E}_t [Q_t^i(\nabla f(x^t)) - \nabla f_i(x^t)], x^t - x^* \rangle = 0,$ where $\mathbb{E}_t \left[ \cdot \right]$ is the conditional expectation at a step $t$. Therefore, the theory for such compressors is highly developed. In our case, $Q_t^i(x^s)$ are unbiased only if $t - s \to \infty$, which follows from asymptotic unbiasedness of our Markovian compressors obtained from Assumption \ref{as:compressors}. However, we can use some coarsening rather than unbiasedness when $t - s = \tau$, where $\tau > \tau_{\text{mix}}(\varepsilon)$, using the technique of \textit{stepping back} as follows:
    \begin{equation*}
        \label{eq:dotprod_tau}
        \expect{\dotprod{Q_t^i( a^{t-\tau}) - a^{t-\tau}}{b^{t-\tau}}} 
        \leq 
        \frac{\varepsilon d}{m} \expect{\norm{a^{t-\tau}} \norm{b^{t-\tau}}} .
    \end{equation*}
    Importantly, we must apply the compressor $Q_t$ at step $t$ to the vector $a^{t-\tau}$ at step $t-\tau$, since if we apply it to the vector $a^t$ at step $t$, we will not be able to uncover the conditional expectation, since we will have randomness in $a^t$ (see details in Appendix \ref{section:main_lemmas}). As can be seen from \eqref{eq:tmp_main} we need to apply the last inequality with $a^{t-\tau} = \nabla f_i(x^{t-\tau})$ and $b^{t - \tau} = x^{t - \tau} - x^*$, but in \eqref{eq:tmp_main} we only obtain expression with variables at step $t$, therefore, it has to be handled in some way. In order to resolve this issue we use a straightforward algebra:
    \begin{align}
        \label{eq:meta_yeban}
        \mathbb{E}\langle Q_t^i( \nabla f_i(x^t)) - \nabla f_i(x^t), x^t - x^* \rangle \notag
        &=\mathbb{E}\langle Q_t^i( \nabla f_i(x^{t-\tau})) - \nabla f_i(x^{t-\tau}), x^{t- \tau} - x^*\rangle
        \\&-
        \mathbb{E}[ \langle Q_t^i( \nabla f_i(x^t) - \nabla f_i(x^{t-\tau})) - \nabla f_i(x^t) + \nabla f_i(x^{t-\tau}), x^t - x^{t-\tau} \rangle] 
        \\&+
        \mathbb{E}[\langle Q_t^i( \nabla f_i(x^t) - \nabla f_i(x^{t-\tau})) - \nabla f_i(x^t) + \nabla f_i(x^{t-\tau}), x^t - x^*\rangle ] \notag
        \\&+
        \mathbb{E}\langle Q_t^i\left( \nabla f_i(x^t) \right) - \nabla f_i(x^t), x^t - x^{t-\tau}\rangle. \notag
    \end{align}
        The first term in the last inequality \eqref{eq:meta_yeban} is solved with the \textit{stepping back}. The other scalar products are solved using the Fenchel-Young inequality. Terms with $\mathbb{E}\|x^t - x^{t-\tau}\|^2$ are evaluated using line 9 of Algorithm~\ref{alg:GD}: $x^t - x^{t - \tau} = -\gamma \sum_{s = t - \tau}^{t-1} g^s$. The terms with $\mathbb{E}\|Q_t^i (\nabla f_i(x^t) - \nabla f_i(x^{t-\tau}))\|^2$ are obtained from the following inequality (see details in Appendix \ref{appendix:GD_odd}):
        \begin{equation*}
        \begin{split}
            \norm{Q_t^i\left(\nabla f(x) - \nabla f(y)\right)}^2
            \leq
            \frac{d^2 L^2}{m^2} \norm{x - y}^2,
        \end{split}
        \end{equation*}
        Since the evaluation of $\mathbb{E}\norm{x^{t+1} - x^*}^2$ raises the terms of the form $\mathbb{E}\norm{x^{t-\tau} - x^*}^2$, we have to do a summation of $\mathbb{E}\norm{x^{t+1} - x^*}^2$ from $t = \tau$ to $t = T$. These terms greatly complicate the proof of Theorem \ref{theorem:GD_odd} compared to unbiased compressors. 
        
        The results of Theorem \ref{theorem:GD_odd} can be rewritten as an upper complexity bound on a number of iterations $T$ of the Algorithm \ref{alg:GD} by carefully tuning the step size $\gamma$.
    \begin{corollary}[Step tuning for Theorem \ref{theorem:GD_odd}]
    \label{corollary:GD_odd}
    $\\$
        $\bullet$~~ 
            Under the conditions of Theorem \ref{theorem:GD_odd} in the non-convex case, choosing $\gamma$ as in Appendix \ref{appenix:corollary_GD_odd_full},
            in order to achieve the $\epsilon$-approximate solution (in terms of $\expect{\norm{\nabla f(x^T)}^2} \leq \epsilon^2$), it takes
            \begin{equation*}
                \mathcal{O}\left( \frac{L \tau d^2}{m^2} F_\tau \left(
                \frac{\delta^2 + 1}{\epsilon^2} 
                +  
                \frac{\sigma^2}{\epsilon^4}
                \right)\right)
            \end{equation*}
        iterations of Algorithm \ref{alg:GD}.
        
        $\bullet$~~
            Under the conditions of Theorem \ref{theorem:GD_odd} in the PL-condition (Assumption \ref{as:PL}) case, choosing $\gamma$ as in Appendix \ref{appenix:corollary_GD_odd_full}
            in order to achieve the $\epsilon$-approximate solution (in terms of $\expect{f(x^t) - f(x^*)} \leq \epsilon$), it takes
            \begin{equation*}
            \begin{split}
                \mathcal{O}\left( \frac{d^2 L \tau}{m^2 \mu} \left(
                (\delta^2+1)\log\left(\frac{1}{\epsilon}\right) + \frac{\sigma^2}{\mu \epsilon} 
                \right)\right)
            \end{split}
            \end{equation*}
            iterations of Algorithm \ref{alg:GD}.
    \end{corollary}
\begin{algorithm}{\texttt{Accelerated Markovian QSGD}}
    \label{alg:acc}
    \begin{algorithmic}[1]
        \State {\bfseries Input:} starting point $x^0 \in 
           \mathbb{R}^d$, step size $\gamma > 0$, momentums $\theta, \eta, \beta, p$, number of iterations $T$
           \For{$t=0$ {\bfseries to} $T$}
           \State Update $x_g^t = \theta x_f^t + (1-\theta) x^t$ \label{line_x_g_acc}
           \State Broadcast $x_g^t$ to all workers
               \For{$i=1$ {\bfseries to} $n$ in parallel}
               \State Set $g_i^t$ = $Q_t^i\left( \nabla f_i(x_g^t) \right)$ \label{line_g_i_acc}
               \State Send $g_i^t$ to the server
               \EndFor
            \State Aggregate $g^t = \frac{1}{n} \sum\limits_{i=1}^n g_i^t$ \label{line_g_acc}
            \State Update $x_f^{t+1} = x_g^t - p\gamma g^t$ \label{line_x_f_acc}
            \State Update \label{line_x_acc} 
            $x^{t+1} = \eta x_f^{t+1} + (p - \eta)x_f^t$\\$\quad\quad\quad\quad\quad~+ (1-p)(1-\beta)x^t + (1-p)\beta x_g^t$ 
           \EndFor
    \end{algorithmic}
\end{algorithm}
\subsection{Accelerated method}

%
After giving the convergence result for the vanilla distributed \texttt{SGD} with Markovian compression operator, we now move on to the accelerated scheme. 

Since we do not assume boundedness of the gradient variance, the classical Nesterov acceleration \cite{10.5555/2670022} does not produce the expected effect, and therefore an additional momentum has to be introduced \cite{Nesterov2012EfficiencyOC, vaswani2019fast}. By applying a multi-step strategy partially similar to \cite{beznosikov2023first}, we obtain our Algorithm \ref{alg:acc}. 

%

        \begin{theorem}[Convergence of \texttt{AMQSGD} (Algorithm \ref{alg:acc})]
        \label{theorem:GD_acc}
        Consider Assumptions \ref{as:lip}, \ref{as:strconv}, \ref{as:sim}. Let the problem \eqref{eq:problem} be solved by Algorithm \ref{alg:acc}. Then for any $\gamma, \varepsilon > 0$, $T > \tau > \tau_{\text{mix}}(\varepsilon), \beta, \theta, \eta, p$ satisfying conditions, described in Appendix \ref{appendix_subsec:full_GD_acc},it holds that
        \begin{equation*}
            \begin{split}
                F_{T+1}
                =
                \mathcal{O} \Bigg(
                \exp\left[-(T - \tau)\sqrt{\frac{p^2\mu\gamma}{3}} \right]F_{\tau}
                + 
                \exp\left[-T\sqrt{\frac{p^2\mu\gamma}{3}} \right]\Delta_{\tau} 
                + 
                \frac{\gamma}{\mu}\sigma^2
                \Bigg).
            \end{split}
        \end{equation*}
        Here we use the notations: $F_t := \E[\|x^{t} - x^*\|^2 + 3/\mu ( f(x_f^{t}) - f(x^*))]$ and $\Delta_{\tau} \leq \gamma^{1/2}\tau^{-4/3}\mu^{-1/3}\sum_{t=0}^{\tau}\big(\E\|\nabla f(x_g^t)\|^2 + \E\|x^{t} - x^*\|^2 + \E[f(x_f^{t}) - f(x^*)] \big)$.
    \end{theorem}
    
    The above theorem shows that in the strongly convex case \texttt{Accelerated Markovian QSGD} with constant step-size can attain sublinear convergence. In terms of dealing with Markovian stochasticity, its proof follows quite similar ideas as the proof of Theorem \ref{theorem:GD_odd}: here again we use the technique of \textit{stepping back} for mixing time, which allows us to effectively deal with the bias of the gradient estimator. 
    The full proof is provided in Appendix \ref{appendix_subsec:full_GD_acc_proof}.
    The results of Theorem \ref{theorem:GD_acc} can be rewritten as an upper complexity bound on a number of iterations $T$ of the Algorithm \ref{alg:acc} by carefully tuning the step size $\gamma$.

        \begin{corollary}[Step tuning for Theorem \ref{theorem:GD_acc}]
        \label{corollary:GD_acc}
            Under the conditions of Theorem \ref{theorem:GD_acc}, choosing $\gamma$ as in Appendix \ref{appendix:corollary_GD_acc_full}
            in order to achieve the $\epsilon$-approximate solution (in terms of $\expect{\norm{x^T - x^*}^2} \leq \epsilon^2$), it takes
            \begin{equation*}
            \begin{split}
                &\mathcal{O}\left( \frac{d^2 L^{\frac{2}{3}} \tau^{\frac{4}{3}}}{m^2 \mu^{\frac{2}{3}}} \left(
                (\delta^2+1)\log\left(\frac{1}{\epsilon}\right) + \frac{\sigma^2}{\mu \epsilon} 
                \right)\right) 
            \end{split}
            \end{equation*}
            iterations of Algorithm \ref{alg:acc}.
    \end{corollary}

    \subsection{Discussion}
    \label{sec:discussion}
    
    Our Example \ref{example:1} and the numerical experiments in Section \ref{section:exp} show that the usage of Markovian compressors could lead to a better performance quite well, however, the theoretical guarantees turn out to be poorer than in the unbiased case. In particular, if we use Rand$m$ in the \texttt{QSGD} algorithm, then we observe the following estimates \cite{beznosikov2023biased}:
    $
        X_T = \mathcal{O}((1 - \mu \gamma)^T X_0 + \gamma\frac{d}{m}\frac{\sigma^2}{\mu n} ),
    $
    where $X_t = \expect{\norms{x^{t} - x^*}^2}\text{ and }\gamma  \lesssim \frac{1}{L (1 + d/mn)}$. Whereas Theorem \ref{theorem:GD_odd} gives us such estimates: 
    $
        F_T  = \mathcal{O} (
                \left( 1 - \frac{\mu \gamma}{12} \right)^{T} F_\tau
                +
                \gamma\frac{d^2}{m^2}\frac{\tau L \sigma^2}{\mu}
                ),
    $
    where $F_t := \expect{f(x^T) - f(x^*)}$ and $\gamma \lesssim \frac{m^2}{L d^2 \tau (\delta^2+1)}$.
    It is important to note that not only has the theory for Markovian compressors not yet been studied well, but also dealing with the Markovian stochasticity itself implies quite strict limitations. For instance,

    \begin{itemize}
        \item $\mathlarger{\mathbf{d/m}}$ \textbf{vs} $\mathlarger{\mathbf{d^2/m^2}}$\textbf{.} We are forced to uniformly bound the noise of the compressor (linearity in the compression constant is prevented by this) due to the impossibility of using the expectation trick, in contrast to the unbiased case \cite{beznosikov2023biased}, where the authors estimated the variance of the compressor noise. The assumption of uniformly bounded noise cannot be rejected by any authors who work with Markovian stochasticity \cite{beznosikov2023first, dorfman2023adapting, doan2020convergence, sun2018markov, even2023stochastic}, therefore, there is no possibility to achieve linearity in the compression rate in our theoretical guaranties, according to the current theoretical advances.
        
        \item \textbf{Mixing time.} Furthermore, it is imperative to emphasize that it follows from Theorems \ref{theorem:GD_odd} and \ref{theorem:GD_acc} that the convergence rate is improved as $\tau$ (and, consequently, $K$) diminishes. In other words, the distribution of the compressor's underlying Markov chain has to converge to a uniform distribution as fast as possible, but empirically one wants the choice of coordinates to depend on previous iterations rather than be random (e.g. for Rand$m$ compressor $\tau = 1, K = 0$). This causes a logical contradiction: while using a large $K$ will theoretically give poorer convergence, in practice algorithms with non-zero values of $K$ perform better (see Section \ref{section:exp}). It is also worth mentioning that when Markovian stochasticity is employed, we can never avoid $\tau$ in our estimates, since it appears in the lower bounds on the convergence rate of methods that involve Markovian properties \cite{bresler2020squares}. Thus, our Algorithms \ref{alg:GD} and \ref{alg:acc} have a reasonably good polynomial dependence on mixing time (Theorem \ref{theorem:GD_odd} shows an optimal estimation in terms of $\tau$), considering the fact there are several works \cite{doan2020finitetime} whose bounds include terms that are even \textit{exponential} in the mixing time.
        
        \item $\bm{L/\mu}$\textbf{.} In spite of the difficulties listed above, we still can observe that the momentums implementation in Algorithm \ref{alg:acc} gives an acceleration in terms of $L/\mu$ compared to vanilla \texttt{QSGD} (Algorithm \ref{alg:GD}). In the classical version of accelerated Gradient Descent, one can achieve an acceleration of the form $\sqrt{L/\mu}$ \cite{nesterov1983method}, but our analysis allows only to achieve $(L/\mu)^{2/3}$ in Theorem \ref{theorem:GD_acc}. When Markovian stochasticity is employed, it is also possible to achieve estimation of the form $\sqrt{L/\mu}$ \cite{beznosikov2023first}, but it is obtained by using batches with size scaled as $2^{j}$, where $j$ is drawn from a truncated geometric distribution. Unfortunately, this specific batching technique cannot be applied in our paper, as we consider compressors that act as random sparsification (Definition \ref{def:rand_spars}), which necessitates that the gradient be compressed only once at each iteration.
        
        \item \textbf{Variance reduction.} In our paper, we focus on the \texttt{QSGD} method and its accelerated version (Algorithms \ref{alg:GD} and \ref{alg:acc}). However, in modern studies on distributed optimization, techniques of variance reduction are of a great interest (\texttt{DIANA} \cite{mishchenko2019distributed}, \texttt{MARINA} \cite{gorbunov2021marina}, \texttt{DASHA} \cite{tyurin2022dasha}), because these methods converge linearly to the exact solution of the problem \eqref{eq:problem}, while \texttt{QSGD} (Algorithms \ref{alg:GD} and \ref{alg:acc}) converges only to the $\sigma^2$-neighborhood of the solution. We implement Markovian compressors (Definitions \ref{def:banlastk} and \ref{def:kawasaki}) in these methods in our experiments, but we do not provide theoretical guarantees for such algorithms since we have just developed a theoretical baseline for the study of Markovian compressors. This represents a promising direction for future research.
    \end{itemize}

    Even though it is not entirely clear whether it is possible to achieve significant improvements in the theoretical results, due to the peculiarities of dealing with Markovian randomness, for now we could only highlight a significantly better performance of Algorithms \ref{alg:GD} and \ref{alg:acc} compared to similar algorithms using a vanilla unbiased compressor Rand$m$ (see Section \ref{section:exp}).



\section{Experiments}
\label{section:exp}

In order to justify the practical usage of the proposed methods and analyze their behavior, we conduct a series of experiments using Markovian compression on distributed optimization problems. Despite the considerations mentioned in Section \ref{sec:discussion}, we observe that Markovian compressors, when used with \texttt{MQSGD} and \texttt{AMQSGD}, as well as with classical \texttt{SGD}, \texttt{DIANA} and \texttt{Adam} show better results compared to the baselines. Appendix~\ref{appx:experiments} provides a description of the technical setup, extended experiments with hyperparameters analysis, and an application of Markovian compressors to model-parallel neural network training.

\subsection{Logistic regression}
\label{section:logistic_regression}

Firstly, we experiment on a classification task using a logistic regression model with $L_2$ regularization of the form:
$
\min_{w \in \R^d} \{f(w) = \frac{1}{n} \sum_{i=1}^{n} \log( 1 + e^{ -y_s w^Tx_s }) + \lambda \|w\|^2 \},
$


with a regularization term $\lambda=0.05$. The dataset is divided among $n=10$ clients. We use \texttt{Mushrooms}, \texttt{A9A}, and \texttt{W8A} datasets from LibSVM~\cite{chang2011libsvm} and MNIST~\cite{deng2012mnist}. 
We experiment with \texttt{MQSGD}, \texttt{AMQSGD}, and \texttt{DIANA} optimizers, employing Rand10\% as a sparsification compressor. Markovian compressors were utilized independently on each client, with history size $K$ taken according to theory in Appendix \ref{appendix:example_1} and forgetting rate $b=50$.
Figure~\ref{fig:exp_logreg}  shows the convergence of the baselines and Markovian compressors on the \texttt{MQSGD} and \texttt{AMQSGD} algorithms on MNIST dataset. Note that improved convergence is achieved without compressor's hyperparameters tuning, making it easy to utilize in any practical optimization task. Moreover, the analysis in Appendix \ref{appendix:example_1} (see Figure \ref{figure:optimal_history_size}) shows that our theory suggests optimal values of history size $K$.

\begin{figure}[h]
    \centering
    \begin{subfigure}
        \centering
        \includegraphics[width=0.46\columnwidth]{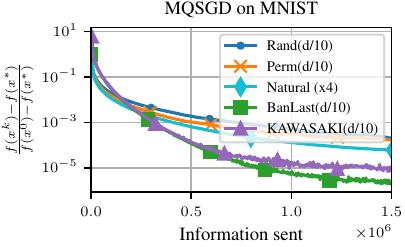}
    \end{subfigure}
    \begin{subfigure}
        \centering \includegraphics[width=0.46\columnwidth]{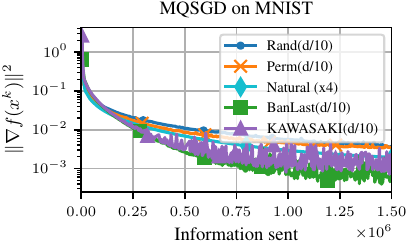}
    \end{subfigure}
    \\
    \begin{subfigure}
        \centering \includegraphics[width=0.46\columnwidth]{ACCGD-MNIST-VS-OTHERS-fdist_ratio.pdf}
    \end{subfigure}
    \begin{subfigure}
        \centering \includegraphics[width=0.46\columnwidth]{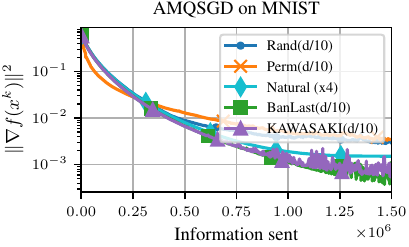}
    \end{subfigure}
    \caption{Logistic Regression on MNIST experiments results. Best runs for each method are displayed.}
    \label{fig:exp_logreg}
\end{figure}

Both Markovian compressors achieve faster convergence than the baseline and more complex compressors: PermK~\cite{szlendak2021permutation} and Natural compression ~\cite{horvath2022natural}. 
Additionally, as our compressors are fully compatible with classical compressors, we conduct experiments with combination with Natural compression in Appendix~\ref{appendix-comparison-combination}.

\subsection{Neural networks}
\subsubsection{Computer Vision}
We also apply Markovian compressors in more complex optimization tasks, such as image classification on CIFAR-10~\cite{krizhevsky2009cifar} dataset with \texttt{ResNet-18} convolutional neural network~\cite{he2016deepresnset}. Formally, we solve the optimization problem:
$
\min_{w \in \R^d} \{f(w) = \frac{1}{n} \sum_{i=1}^{n} l(\sigma(f(x_i, w)), y_i) \}, 
$
where $x_i$ is a training image, $y_i$ is its respective class, $\sigma(\cdot)$ is a Softmax function and $l(\cdot, \cdot)$ is a cross-entropy loss. Dataset is split equally between $n=5$ clients. 
\begin{table*}[htbp]
\small
\captionof{table}{Numerical results of training ResNet-18 on CIFAR-10 with different compressors. 
}
\centering
\begin{tabular}{lccc}
\mbox{} & Rand5\% & \texttt{Banlast} & \texttt{KAWASAKI} \\ \hline Train Loss & \makecell{0.0743} & \makecell{0.0734} & \textbf{0.0305} \\ \hline
Gradient Norm & 1.403 & 1.383 & \textbf{0.745} \\ \hline
Test Accuracy & 87.9 & 88.0 & \textbf{89.05} \\ \hline
\end{tabular}
\label{table:exp_nn_res}
\end{table*}
Figure~\ref{fig:exp_nn} depicts the training loss and gradient norm, with the aggregate values shown in Table~\ref{table:exp_nn_res}. As in the previous case, the application of the Markovian compressor favors faster convergence and better validation results. Hyperparameters, such as the learning rate and batch size, are fine-tuned.
\begin{figure}[h]    
    \centering \includegraphics[width=0.46\columnwidth]{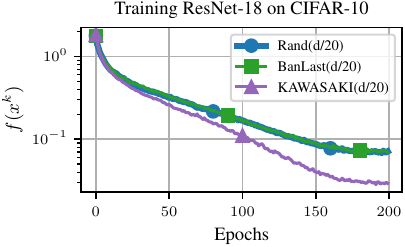}
    \includegraphics[width=0.46\columnwidth]{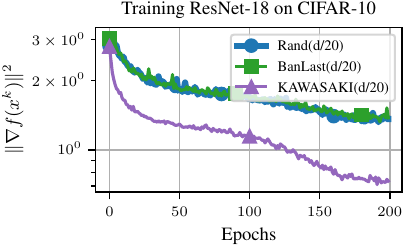}
    \caption {Image classification with ResNet-18 on CIFAR-10 experiments results. Best runs for each method are displayed.}
    \label{fig:exp_nn}
\end{figure}

Figure \ref{fig:exp_nn_2} presents comparison with Permutation and Natural compression, which confirm the practical usefulness of Markovian compressors on more complex and non-convex optimization problems. Note that our methods can be applied in combination with complex sparsification techniques like Natural compression, making our approach even more flexible.

\begin{figure}[h]
    \centering \includegraphics[width=0.46\columnwidth]{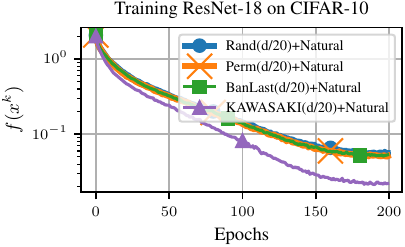}
    \includegraphics[width=0.46\columnwidth]{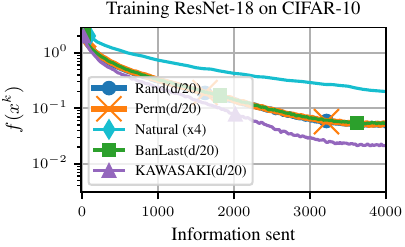}
    \caption {Comparison with other compressors on Resnet-18 training on CIFAR-10 dataset for Rand5\% sparsification on $n=20$ clients. Natural compression factor is 4. Left figure is sequential combination with Natural compression. Right figure is comparison against PermK and Natural compression independently.
    }
    \label{fig:exp_nn_2}
\end{figure}
\subsubsection{Natural Language Processing}
Additionally, we perform a series of real-life experiments of fine-tuning DeBERTaV3-base model \cite{he2021debertadecodingenhancedbertdisentangled} for a subset of the GLUE \cite{wang2019gluemultitaskbenchmarkanalysis} benchmark. Figure \ref{fig:exp_nn_debertav3} presents the convergence rate of \texttt{Adam} \cite{kingma2014adam} with our Markovian compressors for QNLI and SST2 datasets. A detailed description of this experiment is provided in Appendix \ref{appx:exp_nn_deberta}. 
Convergence is also notably enhanced in this case in comparison to the baseline.

\begin{figure}[h]
    \centering 
    \includegraphics[width=0.46\columnwidth]{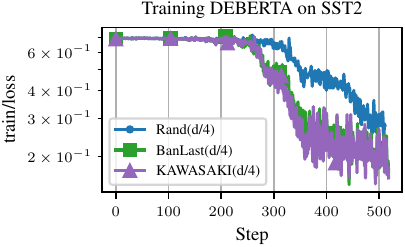}
    \includegraphics[width=0.46\columnwidth]{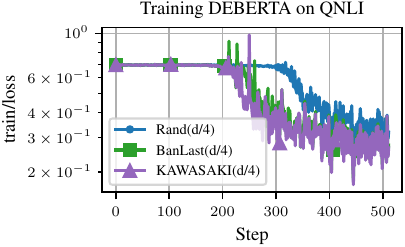}
    \caption {Comparison of the performance of \texttt{BanLast}, \texttt{KAWASAKI} and \texttt{Rand}m on the fine-tuning task on a subset of GLUE benchmark for $n=10$ clients.}
    \label{fig:exp_nn_debertav3}
\end{figure}
\section*{Acknowledgments}
The work was done in the Laboratory of Federated Learning Problems (Supported by Grant App. No. 2 to Agreement No. 075-03-2024-214).
\end{mainpart}

\begin{appendixpart}


\section{Auxiliary Lemmas and Facts}

    In this section we list auxiliary facts and our results that we use several times in our proofs.

    \subsection{Cauchy–Schwarz inequality}
    \label{axil:Cauchy-Schwarz}
        For all $x, y \in \mathbb{R}^d$

        \begin{equation*}
            \dotprod{x}{y} \leq \norms{x}\norms{y} .
        \end{equation*}

    \subsection{Fenchel-Young inequality}
    \label{axil:Fenchel-Young}
    For all $x, y \in \mathbb{R}^d$ and $\beta > 0$
    
    \begin{equation*}
        2 \dotprod{x}{y} \leq \beta^{-1} \|x\|^2 + \beta \|y\|^2 .
    \end{equation*}

\addtocontents{toc}{\protect\setcounter{tocdepth}{2}}

\section{Mathematical calculations from Example \ref{example:1}}
\label{appendix:example_1}

    By definition of the mathematical expectation of an integer positive random variable $Z$, we obtain that $\mathbb{E}[Z] = \sum_{s = 1}^\infty s \cdot \mathbb{P}\{Z = s\}$. 
    In our problem, $Z$ is the number of an iteration where we first selected the desired coordinate.
    For Rand$m$ compressor, we have $\mathbb{P}\{Z = s\} = \frac{m}{d} \cdot \left( 1 - \frac{m}{d}\right)^{s-1}$. The first term is the probability of picking the desired coordinate at iteration $s$ and the second term is the probability of not picking the desired coordinate at iterations from $1$ to $s-1$.
    Using this, the mathematical expectation of the number of steps to quit the point $x^t$ for Rand$m$ compressor is equal to 
    
    \begin{equation}
    \label{eq:rand_m}
     \sum\limits_{s = 1}^{\infty} s \left( 1 - \frac{m}{d}\right)^{s-1} \frac{m}{d} = \frac{d}{m}.
    \end{equation}

    Now we calculate the expectation for \texttt{BanLast}($K, m$) compressor (Definition \ref{def:banlastk}). If $s > K$, similarly to the Rand$m$ case, we obtain that $\mathbb{P}\{Z = s\} = \frac{m}{d - Km} \left( 1 - \frac{m}{d - Km}\right)^{s-1}$, because we cannot choose $K m$ coordinates. If $s \leq K$, then the formula of $\mathbb{P}\{Z = s\}$ becomes a bit more complicated, because the probability of not picking the desired coordinate at iterations from $1$ to $s-1$ is different at each iteration and is equal to $\prod_{h = 0}^{s - 2} \left(1 - \frac{m}{d - hm} \right)$. If $s = 1$, then this probability is equal to one. Using this, we can calculate the mathematical expectation of the number of steps to leave the point $x^t$ for \texttt{BanLast}($K, m$) compressor:
    \begin{equation}
    \begin{split}
    \label{eq:zalupa_K}
        &\sum\limits_{s=1}^K \frac{s m}{d - (s-1) m} \prod_{h = 0}^{s - 2} \left(1 - \frac{m}{d - hm} \right)
        +
        \sum\limits_{s = K+1}^{\infty} s \left( 1 - \frac{m}{d - Km}\right)^{s-1} \frac{m}{d - Km}
        \\&=
        \sum\limits_{s=1}^K \frac{s m}{d - (s-1) m} \prod_{h = 0}^{s - 2} \left(1 - \frac{m}{d - hm} \right)
        +
        \frac{d}{m} \left( 1 - \frac{m}{d - Km}\right)^K
        \\&=
        \sum\limits_{s=1}^K \frac{s}{\alpha - (s-1)} \prod_{h = 0}^{s - 2} \left(1 - \frac{1}{\alpha - h} \right)
        +
        \alpha \left( 1 - \frac{1}{\alpha - K}\right)^K,
    \end{split}
    \end{equation}
    
    where we used the notation $\alpha = d/m$ to show that \eqref{eq:zalupa_K} depends only on $d/m$, but not on $d$ and $m$ separately. We can consider \eqref{eq:zalupa_K} as an optimization problem with respect to $K$. Since $K$ is an integer and the objective function in \eqref{eq:zalupa_K} is complex, we numerically find the optimal $K$ for different $\alpha$. For the sake of clarity, we show the difference between formulas \eqref{eq:rand_m} and \eqref{eq:zalupa_K} on Figure \ref{figure:optimal_history_size}(c).

    We consider $\alpha \in [5.3, 6.7, 8.3, 10, 11.1, 12.5, 14.3, 16.7, 20]$ and find the optimal $K$ by a complete brute force search -- see Figure \ref{figure:optimal_history_size} (a). Then, we perform a linear approximation and obtain the formula $K^*(\alpha) \approx 0.7323\alpha$ -- see Figure \ref{figure:optimal_history_size} (b). Since the correlation coefficient between the points and the approximated line is equal to $0.73$, we can consider this formula to be accurate enough for practical applications.


    \begin{figure}[ht]
    \begin{minipage}{0.33\textwidth}
        \includegraphics[width =  \textwidth]{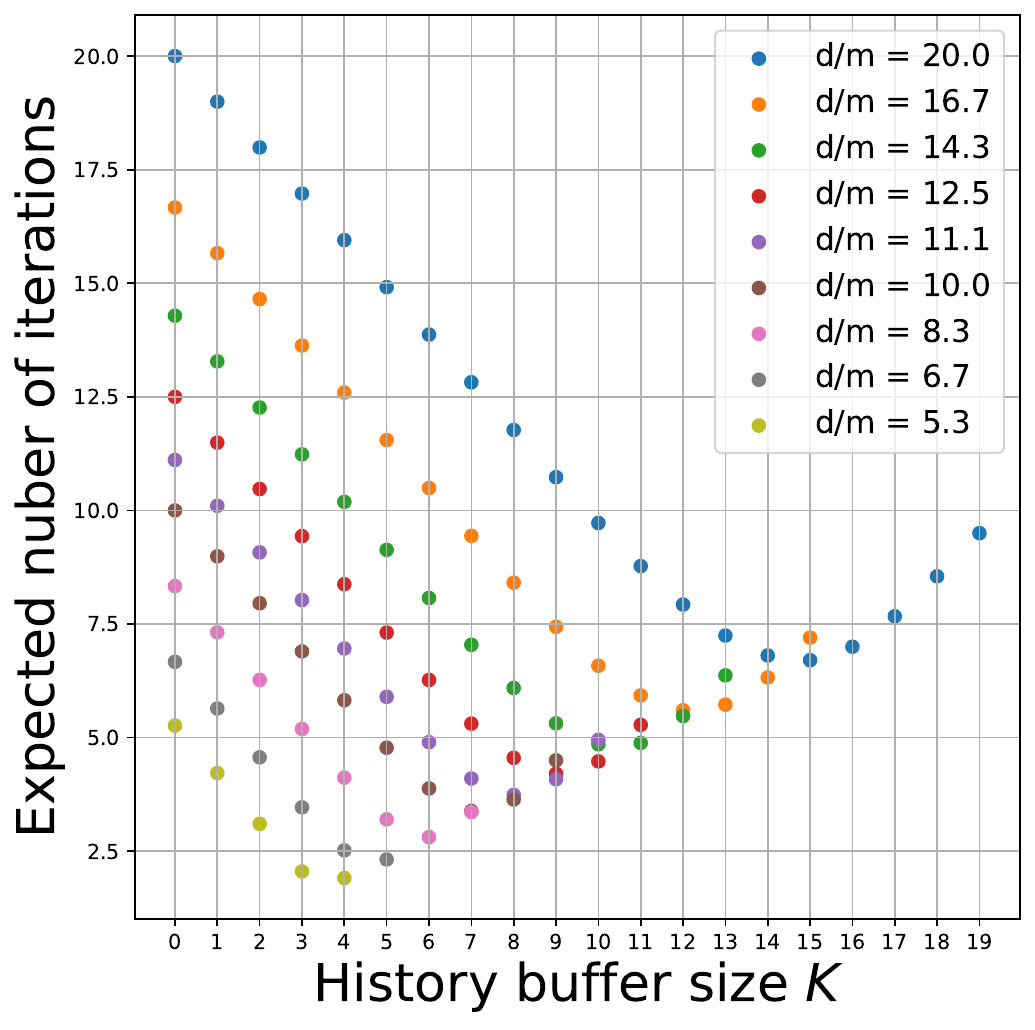}
    \end{minipage}
    \begin{minipage}{0.33\textwidth}
        \includegraphics[width =  \textwidth]{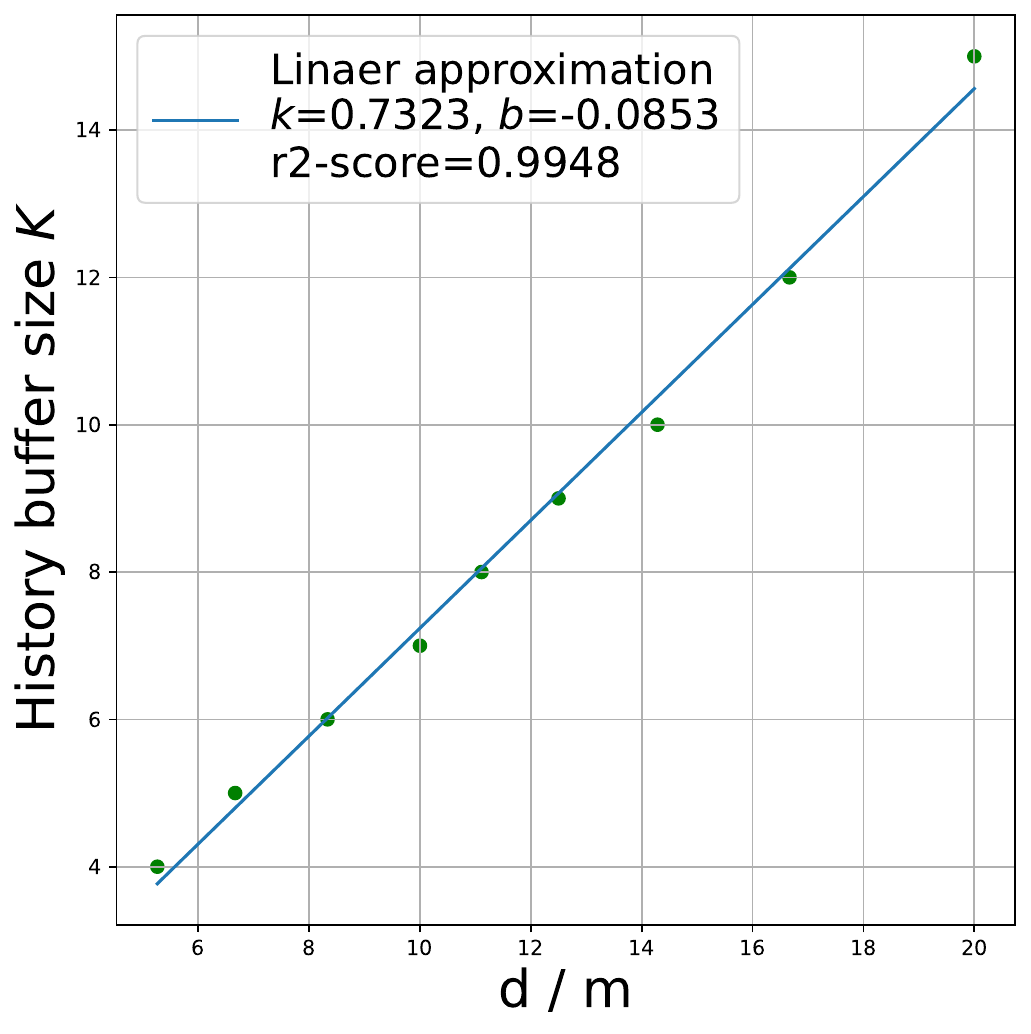}
    \end{minipage}
    \begin{minipage}{0.33\textwidth}
        \includegraphics[width =  \textwidth]{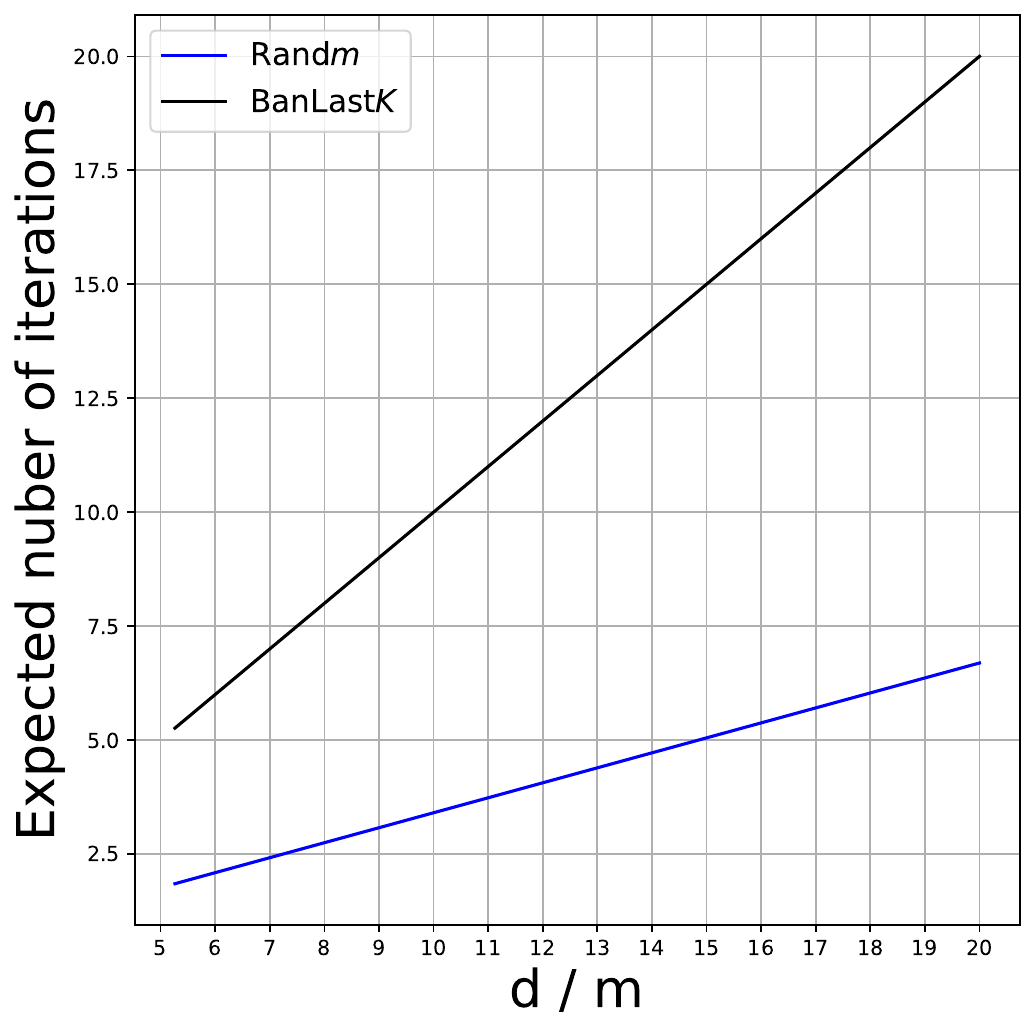}
    \end{minipage}
    \\
    \begin{minipage}{0.01\textwidth}
\quad
\end{minipage}%
\begin{minipage}{0.33\textwidth}
  \centering
(a)
\end{minipage}%
\begin{minipage}{0.34\textwidth}
  \centering
(b)
\end{minipage}%
\begin{minipage}{0.33\textwidth}
  \centering
(c)
\end{minipage}%
    \caption{Theoretical estimate on dependence of history buffer size $K$ on parameter $\alpha = d/m$: (a) represents expected number of iterations required to transfer all coordinates to server on history buffer size $K$ for different $\alpha$, (b) represents scaling of optimal history buffer size $K^*$ on $\alpha$. (c) represents comparison of expected number of iterations required to transfer all coordinates to server on problems parameter $\alpha$ for Rand$m$ and BanLast$K$.}
    \label{figure:optimal_history_size}
    \end{figure}

\section{Proof of Theorem \ref{theorem:compressors}}
\label{appendix:comprssors}
\begin{lemma}
\label{lemma:stationary}
    If $P$ is a transition matrix of a finite homogeneous Markov chain, i.e.
    $$P := (p_{ij})^{n}_{i, j = 1},$$
    where $p_{ij}$ is probability of moving from i to j in one time step.
    And the matrix $P$ is symmetric, i.e. $P^T = P$, then stationary distribution exists and it is uniformly distributed.
\end{lemma} 
\begin{proof}[Proof of Lemma \ref{lemma:stationary}.]
    Let us look at uniform distribution
    $$\pi := \left( \frac{1}{n}, \frac{1}{n}, \ldots, \frac{1}{n} \right).$$
    We can easily obtain that $\pi$ is a stationary distribution, using symmetry and stochastic property of matrix $P$:
    $$\pi P =  \frac{1}{n} \textbf{1}^T P = \frac{1}{n} (P \textbf{1})^T = \frac{1}{n} \textbf{1}^T = \pi .$$
\end{proof}

\begin{proof}[Proof of Theorem \ref{theorem:compressors}.]

    We consider states of Markov chain as $s := \left\{ \nu_1, \nu_2, ..., \nu_K \right\}_{\nu_1, ..., \nu_K \in M}$, where $M$ is the set of all subsets of $\overline{1, d}$ of size $m$. We define $p(s, s', i)$ as the probability to move from state $s$ to state $s'$ for the number of steps $i$.

        $\bullet$ For both compressors \texttt{BanLast}($K, m$) (Definition \ref{def:banlastk}) and \texttt{KAWASAKI}($K, b, \pi_\Delta, m$) (Definition \ref{def:kawasaki}) corresponding Markov chain is finite and indecomposable.
        
        The finiteness of the chain is apparent, as the number of states can be explicitly expressed as $|M| = (C^m_d)^K$.
        We show that both chains are indecomposable below. Then we deduce that the chain is ergodic based on the Ergodic Theorem \cite{neumann1932proof}. Thus, we know that a stationary distribution exists. Than we show that the statinary distribution is uniform over the set of states using Lemma \ref{lemma:stationary}.

        All that remains is to show that both chains are indecomposable and that transition matrixes for both chaines are symmetric. 
        
        We will start with \texttt{BanLast}($K, m$). Restriction on $K, m$ and $d$ is $d > (K+1) m$. That makes obvious that any two states are communicated, i.e. for any $s, s'$ there exists way from $s$ to $s'$. Thus, the Markov chain is indecomposable. 

        For the compressor probability to move from $s$ to $s'$ in one time step can be explicitly expressed as:
        $$p(s, s', 1) = \left( \frac{1}{C^{m}_{d - Km}} \right)^K,$$
        where $C^m_{d - K m} = \frac{(d - Km)!}{m! (d - (K+1)m)!}$ is a binomial coefficient. And all these states are equal in probability. If $d = (K+1) m$, then for $s$ there will be only one set $s'$, such that $p(s, s', 1) > 0$, in this case chain will not be ergodic. If $d > (K+1) m$, then there are more then one state $s'$, for witch $p(s, s', 1) > 0$, therefore chain will be ergodic.

        $\bullet$ According to the Ergodic Theorem, $\rho = (1 - \delta)^{1/N_0}$ and $C = (1 - \delta)^{-1}$, where $N_0$ is the minimal number of iterations through which is strictly greater then zero and $\delta := \min_{s, s'}\{ p(s, s', N_0) \} > 0$. For \texttt{BanLast}($K, m$) in case of $d > (2K +1) m$ it holds that 

        $$N_0 = 2 \text{ and } \delta = p(s, s, 2) = \left( \frac{C^m_{d - 2 K m}}{C^m_{d - K m}} \right)^K \cdot \left( \frac{1}{C^m_{d - Km }} \right)^K, $$
        because the smallest probability is to return to state $s$ in two steps.
        
        $\bullet$ For \texttt{KAWASAKI}($K, b, \pi_\Delta, m$) from any given state, there exists a path to any other state in just one iteration, because probabilities to choose any set of coordinates $\nu$ are non-zero. Thus, the corresponding markov chain is indecomposable.

        We focus on the case where $K = 1$ and that generalize analysis to accommodate larger values of $K$. Let us look at probabilities to move from $\nu_i$ to $\nu_j$ and from $\nu_j$ to $\nu_i$.
        We show that both these probabilities correspond to random choice of the same indexes with the same distribution vector $p$, defined in \ref{def:kawasaki}, i.e. the probabilities are equal.
        For this case let us define $\nu$ as operator
        $$\Psi_i(\overline{1, d}) := \nu_i,$$
        i.e. operator chooses indexes that are in $\nu_i$ from $\overline{1, d}$. And
        $$\Phi(p, \Psi_i) := \mathbb{P}\{ \text{choose $\nu_i$ with distribution vector $p$} \} .$$

        According to \ref{def:kawasaki}, probability to move from $\nu_i$ to $\nu_j$ equals a probability to choose indexes $\nu_j$ with distribution
        $$p_i = \pi_{\Delta}(\widetilde{p}_i), $$
        where
        $$\widetilde{p}_i^k = \begin{cases}
            1 / bd & \text{if } k \in \nu_i \\
            1 / d & \text{if } k \notin \nu_i
        \end{cases},$$
        i.e.
        $$p_{ij} = \Phi(p_i, \Psi_j).$$

        By the definition of $\Phi$, for arbitrary permutation $\phi$ and index choice $\Psi$ holds
        $$\Phi(\phi(p), \Psi \circ \phi) = \Phi(p, \Psi).$$
        
        Now we point out that for arbitrary $\nu_i$ and $\nu_j$ exists permutation $\phi_{ij}$, such that
        $$\Psi_j \circ \phi_{ij} = \Psi_i.$$
        For such permutation holds $\phi_{ij}(\widetilde{p}_i) = \widetilde{p}_j$, i.e. the permutations moves indexes from $\nu_i$ to indexes from $\nu_j$. Then we need to use the property of $\pi_{\Delta}$ to get the same equality for $p_i, p_j$:
        $$\phi_{ij}(p_i) = \phi_{ij}(\pi_{\Delta}(\widetilde{p}_i)) = \pi_{\Delta}{\phi_{ij}((\widetilde{p}_i))} = \pi_{\Delta}(p_j) .$$
        
        This allows us to write
        $$p_{ij} = \Phi(p_i, \Psi_j) = \Phi(\phi_{ij}(p_i), \Psi_j \circ \phi_{ij}) = \Phi(p_j, \Psi_i) = p_{ji} .$$
        
        Thus we get equality of probabilities to move from $\nu_j$ to $\nu_i$ and to opposite way. 
    
        Now we can easily generalize the proof for arbitrary $K$. All that is required is to consider, instead of the sets of indices $\nu$, combinations of sets of indices that were chosen for transmission over the previous $K$ steps. In this way, the number of states is increased, but the logic of reasoning remains unchanged.

        $\bullet$ As was mentioned above, for \texttt{KAWASAKI}($K, b, \pi_\Delta, m$) $N_0 = 1$. We now compute $\delta := p(s, s, 1)$, where $s = \{ \nu, ... , \nu \}$, where $\nu$ occurs $K$ times. In this case probability to choose $\nu$ another $K$ times is equal to $\mathbb{P}\{ j \in \nu \}^{m K}$. And 
        $$\mathbb{P}\{ j \in \nu \} = \min \Bigg\{ \pi_{\Delta} \Bigg[ \widetilde{p} :=  \Big( \underbrace{\frac{1 / d}{b^K} , ... , \frac{1 / d}{b^K}}_{m} , \underbrace{\frac{1 / d}{1} , ... , \frac{1 / d}{1}}_{d - m}\Big)^T \Bigg] \Bigg\}.$$
        If we consider $\left(\pi_{\Delta}\left(\widetilde{p}\right)\right)_j = |\widetilde{p}_j| / \|\widetilde{p}\|_1$, then, since $\|\widetilde{p}\|_1 = \frac{1}{d b^k} (d b^K - m(b^K - 1))$, it hold that $\delta = (d b^K - m(b^K - 1))^{-mK}$. This finishes the proof.
\end{proof}

\section{Main lemmas}
\label{section:main_lemmas}
    \begin{lemma}
        \label{lemma:main_dotprod}
            For any $i \in \overline{1, n}$, $\varepsilon > 0$, $\tau >  \tau_{\text{mix}}(\varepsilon)$, $t > \tau$, for any $a^{t-\tau}, b^{t-\tau} \in \mathbb{R}^d$, such that if we fix all randomness up to step $t-\tau$, $a^{t-\tau}$ and $b^{t-\tau}$ become non-random, it holds that
    
            \begin{equation*}
                \expect{\dotprod{Q_t^i\left( a^{t-\tau} \right) - a^{t-\tau}}{b^{t-\tau}}} 
                \leq 
                \frac{\varepsilon d}{m} \expect{\norm{a^{t-\tau}} \cdot \norm{b^{t-\tau}}} .
            \end{equation*}
        \end{lemma} 
    \begin{proof}
        We begin by using tower property:

        \begin{equation}
        \label{eq:tmp_lemma_GD_1}
            \expect{\dotprod{Q_t^i\left(  a^{t-\tau} \right) - a^{t-\tau}}{b^{t-\tau}}}
            =
            \expect{\dotprod{\mathbb{E}_{t-\tau} \left[Q_t^i\left(  a^{t-\tau} \right) - a^{t-\tau} \right]}{b^{t-\tau}}},
        \end{equation}

        where $\mathbb{E}_{t-\tau} \left[ \cdot \right]$ is the conditional expectation with fixed randomness of all steps up to $t - \tau$. Since on a step $t$ we compress vector $a^{t-\tau}$ according to distribution $\pi_t^i$ by the formula $Q_t^i\left(  a^{t-\tau} \right) = d/m a^{t-\tau} \odot \mathbbm{1}(\nu_t^i)$, where $\nu_t^i$ is some set of $m$ coordinates : $\nu_t^i \subset \overline{1, d}$ and $\mathbbm{1}(\nu_t^i)$ is vector with $1$ on coordinates $\nu_t^i$ on $0$ otherwise. Using this  we can obtain:

        \begin{equation*}
            \mathbb{E}_{t-\tau} \left[Q_t^i\left(  a^{t-\tau} \right)  - a^{t-\tau} \right] = \sum\limits_{\widetilde{\nu
            }_i \in M} \left(\mathbb{P}_{t-\tau}\left\{\nu_t^i = \widetilde{\nu}_i\right\} - \frac{1}{C_d^m} \right) a^{t-\tau} \odot \mathbbm{1}(\widetilde{\nu}_i) \frac{d}{m},
        \end{equation*}

        where $M$ is set of all subsets of $\overline{1, d}$ of size $m$.
        This equality follows from the fact that $\sum_{\widetilde{\nu}_i \in M} a^{t-\tau} \odot \mathbbm{1}(\widetilde{\nu}_i) = C_{d-1}^{m-1} a^{t-\tau}$ and $C_{d-1}^{m-1} / C_{d}^{m} = m/d$. Now with the help of Cauchy–Schwarz inequality \ref{axil:Cauchy-Schwarz} we can estimate \eqref{eq:tmp_lemma_GD_1}:

        \begin{equation}
        \label{eq:tmp_lemma_GD_2}
            \eqref{eq:tmp_lemma_GD_1}
            \leq
            \expect{\sum\limits_{\widetilde{\nu}_i \in M}
            \left| \mathbb{P}_{t-\tau}\left\{\nu_t^i = \widetilde{\nu}_i\right\} - \frac{1}{C_d^m} \right|
            \norm{a^{t-\tau} \odot \mathbbm{1}(\widetilde{\nu}_i)} \frac{d}{m} \norm{b^{t-\tau}}}.
        \end{equation}

        Since $t > \tau$ and $\tau >  \tau_{\text{mix}}(\varepsilon)$ it holds that $\left| \mathbb{P}_{t-\tau}\left\{\nu_t^i = \widetilde{\nu}_i\right\} - 1 / C_d^m \right| \leq \varepsilon \cdot 1 / C_d^m$, because stationary distribution of our Markov chain is uniform. Using the fact that $\norm{a^{t-\tau} \odot \mathbbm{1}(\widetilde{\nu}_i)} \leq \norm{a^{t-\tau}}$ we can obtain:

        \begin{equation*}
            \eqref{eq:tmp_lemma_GD_2} 
            \leq 
            \expect{\sum\limits_{\widetilde{\nu}_i \in M} \varepsilon \frac{1}{C_d^m} \norm{a^{t-\tau}} \frac{d}{m}\norm{b^{t-\tau}}} 
            = 
            \frac{\varepsilon d}{m} \expect{\norm{ a^{t-\tau} } \cdot \norm{b^{t-\tau}}}.
        \end{equation*}

        This finishes the proof.
        
    \end{proof}

    \begin{lemma}
        \label{lemma:main_even}
            For any $i \in \overline{1, n}$, $\varepsilon > 0$, $\tau >  \tau_{\text{mix}}(\varepsilon)$, $t > \tau$, for any $a^{t-\tau} \in \mathbb{R}^d$, such that if we fix all randomness up to step $t-\tau$, $a^{t-\tau}$ becomes non-random, it holds that
    
            \begin{equation*}
                \expect{\norms{\mathbb{E}_{t - \tau} \left[ Q_t^i(a^{t-\tau}) \right] - a^{t-\tau}}^2}
                \leq
                \frac{\varepsilon^2 d^2}{m^2}\expect{\norm{a^{t-\tau}}^2} .
            \end{equation*}
    \end{lemma}
    
    \begin{proof}
        Using same notation as in the proof of Lemma \ref{cor:main_dotprod} we obtain 

        \begin{equation*}
        \begin{split}
            \E\Big[\|\mathbb{E}_{t - \tau} \left[ Q_t^i(a^{t-\tau}) \right] - &a^{t-\tau}\|^2\Big]
            =
            \expect{\norms{ \sum\limits_{\widetilde{\nu}_i \in M} \left( \mathbb{P}_{t-\tau} \left\{ \nu_t^i = \widetilde{\nu}_i\right\} - \frac{1}{C_d^m} \right) \frac{d}{m} a^{t-\tau} \odot \mathbbm{1}(\widetilde{\nu}_i)}^2}
            \\&\leq
            \expect{\frac{d^2}{m^2}C_d^m \sum\limits_{\widetilde{\nu}_i \in M} \left(\left| \mathbb{P}_{t-\tau} \left\{ \nu_t^i = \widetilde{\nu}_i\right\} - \frac{1}{C_d^m} \right|^2 \norm{a^{t-\tau} \odot \mathbbm{1}(\widetilde{\nu}_i)}^2 \right)} .
        \end{split}
        \end{equation*}

        Since $t > \tau$ and $\tau >  \tau_{\text{mix}}(\varepsilon)$ it holds that $\left| \mathbb{P}_{t-\tau}\left\{\nu_t^i = \widetilde{\nu}_i\right\} - 1 / C_d^m \right| \leq \varepsilon \cdot 1 / C_d^m$, because stationary distribution of our Markov chain is uniform. Using the fact that $\norm{a^{t-\tau} \odot \mathbbm{1}(\widetilde{\nu}_i)} \leq \norm{a^{t-\tau}}$ we can obtain:

        \begin{equation*}
            \expect{\norms{\mathbb{E}_{t - \tau} \left[ Q_t^i(a^{t-\tau}) \right] - a^{t-\tau}}^2}
            \leq
            \frac{\varepsilon^2 d^2}{m^2}\expect{\norm{a^{t-\tau}}^2} .
        \end{equation*}

        This finishes the proof.
        
    \end{proof}


    \begin{lemma}
    \label{lemma:main_norms}
        For any $i \in \overline{1, n}$ and $a \in \mathbb{R}^d$ it holds that

        \begin{equation*}
            \norm{Q^i(a)}^2 \leq \frac{d^2}{m^2} \norm{a}^2 
            ~~\text{ and }~~
            \norm{Q^i(a) - a}^2 \leq 4 \frac{d^2}{m^2} \norm{a}^2 .
        \end{equation*}
    \end{lemma}

    \begin{proof}
        Consider the first inequality. Since $Q^i\left(  a \right) = d/m a \odot \mathbbm{1}(\nu^i)$, then $\norm{Q^i(a)} \leq d/m \norm{a}$, therefore 

        $$\norm{Q^i(a)}^2 \leq \frac{d^2}{m^2} \norm{a}^2.$$

        Consider the second inequality. Using Fenchel-Young inequality  \ref{axil:Fenchel-Young} with $\beta = 1$ we can estimate

        \begin{equation*}
            \norms{Q^i(a) - a}^2 
            \leq
            2 \norm{Q^i(a)}^2 + 2 \norm{a}^2 
            \leq
            2 \left( \frac{d^2}{m^2} + 1\right) \norm{a}^2 
            \leq
            4 \frac{d^2}{m^2} \norm{a}^2 .
        \end{equation*}

        This finishes the proof.
    \end{proof}

    \begin{corollary}
            \label{cor:main_dotprod}
            For any $i \in \overline{1, n}$, $\varepsilon > 0$, $\tau >  \tau_{\text{mix}}(\varepsilon)$, $t > \tau$, for any $a^{t}, b^{t} \in \mathbb{R}^d$, such that if we fix all randomness up to step $t$, $a^{t}$ and $b^{t}$ become non-random. And for any $\hat{a}^{t-\tau}, \hat{b}^{t-\tau}$, such that if we fix all randomness up to step $t - \tau$, $\hat{a}^{t-\tau}$ and $\hat{b}^{t-\tau}$ become non-random, it holds that
    
            \begin{equation*}
            \begin{split}
                &2 \left|\expect{\dotprod{Q_t^i\left( a^{t} \right) - a^{t}} {b^{t}}} \right|
                \leq
                \frac{\varepsilon d}{m \beta_0} \expect{\norm{\hat{a}^{t-\tau}}^2}
                +
                \frac{\varepsilon d \beta_0}{m} \expect{\norm{\hat{b}^{t-\tau}}^2}
                +
                \frac{1}{\beta_2} \expect{\norm{b^t}^2},
                \\&+
                \left( \frac{1}{\beta_1} + \frac{1}{\beta_3} \right)
                \expect{\norm{b^t - \hat{b}^{t-\tau}}^2}
                +
                4 \frac{d^2}{m^2} \beta_3 \expect{\norm{a^t}^2} 
                +
                4 \frac{d^2 \left( \beta_1 + \beta_2 \right)}{m^2}  \expect{\norm{a^t - \hat{a}^{t-\tau}}^2}
            \end{split}
            \end{equation*}

            where $\beta_0, \beta_1, \beta_2, \beta_3 > 0$.
        \end{corollary}

    \begin{proof}
        Using straightforward algebra we obtain 

        \begin{equation*}
        \begin{split}
            \expect{\dotprod{Q_t^i\left( a^{t} \right) - a^{t}}{b^{t}}}
            &=
            \expect{\dotprod{Q_t^i\left( \hat{a}^{t-\tau} \right) - \hat{a}^{t-\tau}}{\hat{b}^{t-\tau}}} 
            \\&-
            \expect{\dotprod{Q_t^i\left( a^t - \hat{a}^{t-\tau} \right) - a^t + \hat{a}^{t-\tau}}{b^t - \hat{b}^{t-\tau}}}
            \\&+
            \expect{\dotprod{Q_t^i\left( a^t - \hat{a}^{t-\tau} \right) - a^t + \hat{a}^{t-\tau}}{b^t}}
            \\&+
            \expect{\dotprod{Q_t^i\left( a^t \right) - a^t}{b^{t} - \hat{b}^{t-\tau}}} .
        \end{split}
        \end{equation*}

        Using Lemma \ref{lemma:main_dotprod} with $a^{t-\tau} = \hat{a}^{t-\tau}, b^{t-\tau} = \hat{b}^{t-\tau}$ and Fenchel-Young inequality \ref{axil:Fenchel-Young} with $\beta_1, \beta_2, \beta_3 > 0$ we obtain:

        \begin{equation*}
        \begin{split}
            2 \left|\expect{\dotprod{Q_t^i\left( a^{t} \right) - a^{t}} {b^{t}}} \right|
            &\leq
            2 \frac{\varepsilon d}{m} \expect{\norm{\hat{a}^{t-\tau}} \cdot \norm{\hat{b}^{t-\tau}}}
            \\&+
            \beta_1 \expect{\norm{Q_t^i\left( a^t - \hat{a}^{t-\tau} \right) - a^t + \hat{a}^{t-\tau}}^2}
            +
            \frac{1}{\beta_1} \expect{\norm{b^t - \hat{b}^{t-\tau}}^2}
            \\&+
            \beta_2 \expect{\norm{Q_t^i\left( a^t - \hat{a}^{t-\tau} \right) - a^t + \hat{a}^{t-\tau}}^2}
            +
            \frac{1}{\beta_2} \expect{\norm{b^t}^2}
            \\&+
            \beta_3 \expect{\norm{Q_t^i\left( a^t \right) - a^t}^2}
            +
            \frac{1}{\beta_3} \expect{\norm{b^t - \hat{b}^{t - \tau}}^2} .
        \end{split}
        \end{equation*}

        Using Lemma \ref{lemma:main_norms} and Fenchel-Young inequality \ref{axil:Fenchel-Young} with $\beta_0 > 0$ we obtain

        \begin{equation*}
        \begin{split}
            2 \left|\expect{\dotprod{Q_t^i\left( a^{t} \right) - a^{t}} {b^{t}}} \right|
            &\leq
            \frac{\varepsilon d}{m \beta_0} \expect{\norm{\hat{a}^{t-\tau}}^2}
            +
            \frac{\varepsilon d \beta_0}{m} \expect{\norm{\hat{b}^{t-\tau}}^2}
            \\&+
            4 \frac{d^2}{m^2} \left( \beta_1 + \beta_2 \right) \expect{\norm{a^t - \hat{a}^{t-\tau}}^2}
            +
            \left( \frac{1}{\beta_1} + \frac{1}{\beta_3} \right)
            \expect{\norm{b^t - \hat{b}^{t-\tau}}^2}
            \\&+
            4 \frac{d^2}{m^2} \beta_3 \expect{\norm{a^t}^2} 
            +
            \frac{1}{\beta_2} \expect{\norm{b^t}^2} . 
        \end{split}
        \end{equation*}

        This finishes the proof.
    \end{proof}

    \begin{lemma}
        
        Assume \ref{as:sim}, then for any $x \in \mathbb{R}^d$ it holds that

        \begin{equation*}
            \label{lemma:main_sim}
            \frac{1}{n} \sum\limits_{i=1}^n \norm{\nabla f_i(x)}^2 
            \leq
            2 (\delta^2 + 1) \norm{\nabla f(x)}^2 + 2 \sigma^2 .
        \end{equation*}
    \end{lemma}

    \begin{proof}
        Using straightforward algebra and Fenchel-Young inequality \ref{axil:Fenchel-Young} with $\beta = 1$  we obtain 

        \begin{equation*}
        \begin{split}
            \frac{1}{n} \sum\limits_{i=1}^n \norm{\nabla f_i(x)}^2 
            &\leq
            \frac{2}{n} \sum\limits_{i=1}^n \norm{\nabla f_i(x) - \nabla f(x)}^2 
            +
            2 \norm{\nabla f(x)}^2 
            \\&\leq
            2(\delta^2 + 1) \norm{\nabla f(x)}^2 + 2 \sigma^2 .
        \end{split}
        \end{equation*}

        The last inequity follows from \ref{as:sim}. This finishes the proof.
    \end{proof}

\section{Extensions for Theorem \ref{theorem:GD_odd}}
\label{appendix:GD_odd}

\subsection{Full version of Theorem \ref{theorem:GD_odd}}
\label{appendix_subsec:full_GD_odd}

\begin{theorem}[Convergence of \texttt{MQSGD} (Algorithm \ref{alg:GD}), extension of \ref{theorem:GD_odd}]
    \label{theorem:full_GD_odd}
        Consider Assumptions \ref{as:lip}, \ref{as:sim} and  \ref{as:compressors}. Let problem \eqref{eq:problem} be solved by Algorithm \ref{alg:GD}.
        
        $\bullet$~~
            For any $\varepsilon > 0$, $\gamma > 0$, $\tau > \tau_{\text{mix}}(\varepsilon)$ and $T > \tau$ satisfying
            $$\gamma \lesssim \frac{m^2}{d^2 L (\delta^2 + 1) \tau} 
            ~~\text{ and }~~
            \varepsilon \lesssim \frac{m^2}{d^2 (\delta^2 + 1)}, $$
            it holds that
            \begin{equation*}
                \expect{\norm{\nabla f(\widehat{x}^{T})}^2} 
                =
                \mathcal{O} \left( \frac{F_\tau}{\gamma T}
                +
                \frac{\gamma L \tau d^2}{m^2} \sigma^2 \right),
            \end{equation*}
               
            where $\widehat{x}^{T}$ is chosen uniformly from $\left\{x^t \right\}_{t=0}^{T}$.

        $\bullet$~~
            If $f$ additionally verifies the PL-condition (Assumption \ref{as:PL}), then for any $\varepsilon > 0$, $\gamma > 0$, $\tau > \tau_{\text{mix}}(\varepsilon)$ and $T > \tau$ satisfying
            
            $$\gamma \lesssim \frac{m^2}{L d^2 \tau (\delta^2+1)}
            ~~\text{ and }~~
            \varepsilon =  \sqrt{\gamma L \tau} \lesssim \frac{m}{d \sqrt{\delta^2+1}}, $$
            
            it holds that
    
            \begin{equation*}
            \begin{split}
                F_T
                &=
                \mathcal{O} \left(
                \left( 1 - \frac{\mu \gamma}{12} \right)^{T-\tau} F_\tau
                +
                \frac{\gamma d^2 L \tau}{\mu m^2}\sigma^2 
                \right) .
            \end{split}
            \end{equation*}

            Here we use a notation $F_t := \expect{f(x^t) - f(x^*)}$ .
    \end{theorem}

\subsection{Full version of Corollary \ref{corollary:GD_odd}}
\label{appenix:corollary_GD_odd_full}

\begin{corollary}[Step tuning for Theorem \ref{theorem:GD_odd}, extension of Corollary \ref{corollary:GD_odd}]
    $\\$
        $\bullet$~~ 
            Under the conditions of Theorem \ref{theorem:GD_odd} in the non-convex case, choosing $\gamma$ as 
            
            \begin{equation*}
            \gamma \lesssim \frac{m}{d \sqrt{L \tau}}\min\left\{ 
                                        \frac{m}{d (\delta^2 + 1) \sqrt{L \tau}} ~;~
                                        \sqrt{\frac{F_\tau}{T \sigma^2}},
                                    \right\},
            \end{equation*}
            
            in order to achieve $\epsilon$-approximate solution (in terms of $\expect{\norm{\nabla f(x^T)}^2} \leq \epsilon^2$) it takes
            
            \begin{equation*}
                \mathcal{O}\left( \frac{L \tau d^2}{m^2} F_\tau \left(
                \frac{\delta^2 + 1}{\epsilon^2} 
                +  
                \frac{\sigma^2}{\epsilon^4}
                \right)\right) \text{ iterations of Algorithm \ref{alg:GD}} .
            \end{equation*}
        
        $\bullet$~~
            Under the conditions of Theorem \ref{theorem:GD_odd} in the PL-condition (Assumption \ref{as:PL}) case, choosing $\gamma$ as
            
            \begin{equation*}\gamma \lesssim \min\left\{
                                    \frac{m^2}{L d^2 \tau (\delta^2+1)}
                                    ~;~
                                    \frac{\log\left( \max\left\{ 2 ;
                                    \frac{\mu^2 m^2 F_\tau T}{d^2 L \tau \sigma^2}\right\} \right)}{\mu T}
                                    \right\},
                                    \end{equation*}

            in order to achieve $\epsilon$-approximate solution (in terms of $\expect{f(x^t) - f(x^*)} \leq \epsilon$) it takes
            
            \begin{equation*}
                \mathcal{O}\left( \frac{d^2 L \tau}{m^2 \mu} \left(
                (\delta^2+1)\log\left(\frac{1}{\epsilon}\right) + \frac{\sigma^2}{\mu \epsilon} 
                \right)\right) \text{ iterations of Algorithm \ref{alg:GD}} .
            \end{equation*}
            
    \end{corollary}

\subsection{Proof of Theorem \ref{theorem:GD_odd}, non-convex case}
\label{appendix_subsec:non_conv}

    \begin{proof}
        Denoting $F_t := \expect{f(x^t) - f(x^*)}$, we have using $L$-smoothness:

        \begin{equation}
        \begin{split}
        \label{eq:tmp_even_1_1}
            F_{t+1} - F_t 
            &\leq
            - \gamma \expect{\dotprod{\frac{1}{n} \sum\limits_{i=1}^n Q_t^i (\nabla f_i(x^t))}{\nabla f(x^t)}}
            +
            \frac{\gamma^2 L}{2} \expect{\norm{\frac{1}{n} \sum\limits_{i=1}^n Q_t^i(\nabla f_i(x^t))}^2} .
        \end{split}
        \end{equation}

        Consider $- \gamma \expect{\dotprod{\frac{1}{n} \sum\limits_{i=1}^n Q_t^i (\nabla f_i(x^t))}{\nabla f(x^t)}}$. Using straightforward algebra: $\pm \nabla f_i(x^{t-\tau})$ and $\pm \nabla f(x^{t-\tau})$ we can re-write this term:

        \begin{equation*}
        \begin{split}
            &- \gamma \E\Bigg[\Bigg\langle\frac{1}{n} \sum\limits_{i=1}^n Q_t^i (\nabla f_i(x^t)), \nabla f(x^t)\Bigg\rangle\Bigg]
            \\&=
            \underbrace{ - \gamma \expect{ \dotprod{ \frac{1}{n} \sum\limits_{i=1}^n Q_t^i (\nabla f_i(x^{t-\tau})) }{ \nabla f(x^{t-\tau}) } } }_{\circledOne}
            \\&\quad
            \underbrace{ - \gamma \expect{ \dotprod{ \frac{1}{n} \sum\limits_{i=1}^n Q_t^i (\nabla f_i(x^t)) }{ \nabla f(x^t) - \nabla f(x^{t-\tau}) } } }_{\circledTwo}
            \\&\quad
            \underbrace{ - \gamma \expect{ \dotprod{ \frac{1}{n} \sum\limits_{i=1}^n Q_t^i (\nabla f_i(x^t) - \nabla f_i(x^{t-\tau})) }{ \nabla f(x^{t-\tau}) } } }_{\circledThree} .
        \end{split}
        \end{equation*}

        Consider $\circledOne$. Using straightforward algebra, tower property, Lemmas \ref{lemma:main_even} and \ref{lemma:main_sim} we obtain

        \begin{equation}
        \begin{split}
        \label{eq:tmp_even_1_2}
            \circledOne
            &=
             - \gamma \expect{ \dotprod{ \frac{1}{n} \sum\limits_{i=1}^n \mathbb{E}_{t-\tau} \left[ Q_t^i (\nabla f_i(x^{t-\tau}))\right] }{ \nabla f(x^{t-\tau}) }}
             \\&=
             -\frac{\gamma}{2} \expect{\norm{\frac{1}{n} \sum\limits_{i=1}^n \mathbb{E}_{t-\tau} \left[ Q_t^i (\nabla f_i(x^{t-\tau}))\right]}^2}
             \\&~~~~\,+
             \frac{\gamma}{2} \expect{\norm{\nabla f(x^{t-\tau}) - \frac{1}{n} \sum\limits_{i=1}^n \mathbb{E}_{t-\tau} \left[ Q_t^i (\nabla f_i(x^{t-\tau}))\right]}^2}
             -
             \frac{\gamma}{2} \expect{\norms{\nabla f(x^{t-\tau})}^2}
             \\&\leq
             \frac{\gamma}{2} \varepsilon^2 \frac{d^2}{m^2} \frac{1}{n} \sum\limits_{i=1}^n \expect{\norm{\nabla f_i(x^{t-\tau})}^2}
             - 
             \frac{\gamma}{2} \expect{\norms{\nabla f(x^{t-\tau})}^2}
             \\&\leq
             \gamma\left( \varepsilon^2 \frac{d^2}{m^2} (\delta^2 +1) - \frac{1}{2} \right) \expect{\norms{\nabla f(x^{t-\tau})}^2}
             +
             \gamma \varepsilon^2 \frac{d^2}{m^2} \sigma^2
             \\&\leq
             - \frac{\gamma}{4} \expect{\norms{\nabla f(x^{t-\tau})}^2}
             +
             \gamma \varepsilon^2 \frac{d^2}{m^2} \sigma^2 .
        \end{split}
        \end{equation}

        The last inequality follows from the fact, that $$\varepsilon \leq \frac{m}{2 d \sqrt{\delta^2 + 1}} .$$

        Consider $\circledTwo$. Using Cauchy-Schwarz \ref{axil:Cauchy-Schwarz} and Fenchel-Young \ref{axil:Fenchel-Young} with $\beta = 1$ inequalities we obtain

        \begin{equation}
        \begin{split}
        \label{eq:tmp_even_1_3}
            \circledTwo
            &\leq
            \expect{\norm{-\frac{\gamma}{n} \sum\limits_{i=1}^n Q_t^i (\nabla f_i(x^t))} \norm{\nabla f(x^t) - \nabla f(x^{t - \tau})}}
            \\&\leq
            \gamma L \expect{\norm{\frac{1}{n} \sum\limits_{i=1}^n Q_t^i (\nabla f_i(x^t))} \norm{x^t - x^{t - \tau}}}
            \\&=
            \gamma^2L \expect{\norm{\frac{1}{n} \sum\limits_{i=1}^n Q_t^i (\nabla f_i(x^t))} \norm{\sum\limits_{s = t - \tau}^{t-1} \frac{1}{n} \sum\limits_{i=1}^n Q_s^i(\nabla f_i(x^s))}}
            \\&\leq
            \frac{\gamma^2 L}{2} \left( \tau \expect{\norm{\frac{1}{n} \sum\limits_{i=1}^n Q_t^i (\nabla f_i(x^t))}^2 }
            + 
            \sum\limits_{s = t - \tau}^{t-1} \expect{\norm{\frac{1}{n} \sum\limits_{i=1}^n Q_s^i (\nabla f_i(x^s))}^2}
            \right) .
        \end{split}
        \end{equation}

        Third equality holds since $x^t - x^{t-\tau} = \gamma \sum_{s = t - \tau}^{t-1} \frac{1}{n} \sum\limits_{i=1}^n Q_s^i(\nabla f_i(x^s))$. Consider $\circledThree$. Using Cauchy-Schwarz \ref{axil:Cauchy-Schwarz} and Fenchel-Young \ref{axil:Fenchel-Young} with $\beta = m/d$ inequalities we obtain

        \begin{equation}
        \begin{split}
        \label{eq:tmp_even_1_4}
            \circledThree
            &\leq
            \expect{\norm{-\frac{\gamma}{n} \sum\limits_{i=1}^n Q_t^i (\nabla f_i(x^t) - \nabla f_i(x^{t-\tau}))} \norm{\nabla f(x^{t - \tau})}}
            \\&\leq
            \gamma L \expect{\norm{\frac{1}{n} \sum\limits_{i=1}^n Q_t^i (\nabla f_i(x^t) - \nabla f_i(x^{t-\tau})} \norm{\nabla f(x^{t- \tau})}}
            \\&\leq
            \gamma^2L \frac{d}{m} \expect{\norm{\frac{1}{n} \sum\limits_{i=1}^n Q_t^i (\nabla f_i(x^t))} \norm{\sum\limits_{s = t - \tau}^{t-1} \frac{1}{n} \sum\limits_{i=1}^n Q_s^i(\nabla f_i(x^s))}}
            \\&\leq
            \frac{\gamma^2 L}{2} \left(\sum\limits_{s = t - \tau}^{t-1} \expect{\norm{\frac{1}{n} \sum\limits_{i=1}^n Q_s^i (\nabla f_i(x^s))}^2}
            + 
            \frac{d^2 \tau}{m^2} \expect{\norm{\nabla f(x^{t - \tau})}^2}
            \right) .
        \end{split}
        \end{equation}

        Wrapping \eqref{eq:tmp_even_1_1} - \eqref{eq:tmp_even_1_4} up we obtain

        \begin{equation*}
        \begin{split}
            F_{t+1} - F_t 
            &\leq 
            \frac{\gamma^2 L}{2} \expect{\norm{\frac{1}{n} \sum\limits_{i=1}^n Q_t^i(\nabla f_i(x^t))}^2}
            - \frac{\gamma}{4} \expect{\norms{\nabla f(x^{t-\tau})}^2}
            +
            \gamma \varepsilon^2 \frac{d^2}{m^2} \sigma^2
            \\&\quad+
            \frac{\gamma^2 L}{2} \left( \tau \expect{\norm{\frac{1}{n} \sum\limits_{i=1}^n Q_t^i (\nabla f_i(x^t))}^2 }
            + 
            \sum\limits_{s = t - \tau}^{t-1} \expect{\norm{\frac{1}{n} \sum\limits_{i=1}^n Q_s^i (\nabla f_i(x^s))}^2}
            \right)
            \\&\quad+
            \frac{\gamma^2 L}{2} \left(\sum\limits_{s = t - \tau}^{t-1} \expect{\norm{\frac{1}{n} \sum\limits_{i=1}^n Q_s^i (\nabla f_i(x^s))}^2}
            + 
            \frac{d^2 \tau}{m^2} \expect{\norm{\nabla f(x^{t - \tau})}^2}
            \right)
            \\&\leq
            \gamma^2 L \tau \expect{\norm{\frac{1}{n} \sum\limits_{i=1}^n Q_t^i(\nabla f_i(x^t))}^2}
            +
            \gamma \varepsilon^2 \frac{d^2}{m^2} \sigma^2
            \\&\quad+
            \gamma^2 L \sum\limits_{s = t - \tau}^{t-1} \expect{\norm{\frac{1}{n} \sum\limits_{i=1}^n Q_s^i (\nabla f_i(x^s))}^2}
            + 
            \left( \frac{\gamma^2 L \tau d^2}{2 m^2} - \frac{\gamma}{4} \right) \expect{\norm{\nabla f(x^{t - \tau})}^2} .
        \end{split}
        \end{equation*}

        Using Lemma \ref{lemma:main_sim} we obtain

        \begin{equation}
        \begin{split}
        \label{eq:tmp_need_for_PL}
            F_{t+1} - F_t 
            &\leq 
            \frac{2 d^2 \gamma^2 L \tau}{m^2} \left( (\delta^2 + 1) \expect{\norm{\nabla f(x^t)}^2} + \sigma^2 \right)
            + 
            \left( \frac{\gamma^2 L \tau d^2}{2 m^2} - \frac{\gamma}{4} \right) \expect{\norm{\nabla f(x^{t - \tau})}^2}
            \\&\quad+
            \frac{2 d^2 \gamma^2 L}{m^2} \sum\limits_{s = t - \tau}^{t-1} \left( (\delta^2 + 1) \expect{\norm{\nabla f(x^s)}^2} + \sigma^2 \right)
            +
            \gamma \varepsilon^2 \frac{d^2}{m^2} \sigma^2
            \\&=
            \frac{2 d^2 \gamma^2 L (\delta^2+1) \tau}{m^2} \expect{\norm{\nabla f(x^t)}^2}
            +
            \frac{2 d^2 \gamma^2 L (\delta^2 + 1)}{m^2} \sum\limits_{s = t - \tau}^{t-1} \expect{\norm{\nabla f(x^s)}^2}
            \\&\quad+
            \left( \frac{\gamma^2 L \tau d^2}{2 m^2} - \frac{\gamma}{4} \right) \expect{\norm{\nabla f(x^{t - \tau})}^2}
            +
            \frac{\gamma d^2}{m^2}\left(4 \gamma L \tau + \varepsilon^2 \right)\sigma^2 .
        \end{split}
        \end{equation}

        Summing \eqref{eq:tmp_need_for_PL} from $t = \tau$ to $t = T$ and using the fact that $\varepsilon^2 \leq \gamma L \tau$ and $1 + \delta^2 \geq 1$ we obtain

        \begin{equation*}
        \begin{split}
            \sum\limits_{t = \tau}^T \frac{\gamma}{4} \expect{\norm{\nabla f(x^{t - \tau})}^2}
            &\leq
            F_{\tau}
            +
            \frac{2 d^2 \gamma^2 L (\delta^2 + 1)}{m^2} 
            \Bigg( 
                \tau \sum\limits_{t = \tau}^T \expect{\norm{\nabla f(x^{t})}^2}
                \\&\hspace{-10mm}+
                \sum\limits_{t = \tau}^T \sum\limits_{s = t - \tau}^{t-1} \expect{\norm{\nabla f(x^s)}^2}
                +
                \tau \sum\limits_{t = \tau}^T\expect{\norm{\nabla f(x^{t - \tau})}^2}
            \Bigg)
            +
            \sum\limits_{t = \tau}^T 5 \frac{\gamma^2 L \tau d^2}{m^2} \sigma^2 .
        \end{split}
        \end{equation*}

        Since $ \sum_{t = \tau}^T \sum_{s = t - \tau}^{t-1} \expect{\norm{\nabla f(x^s)}^2} \leq \tau \sum_{t = 0}^T \expect{\norm{\nabla f(x^t)}^2}$, we get

        \begin{equation*}
            \gamma \sum\limits_{t = 0}^{T-\tau} \expect{\norm{\nabla f(x^{t})}^2}
            \leq
            4 F_\tau
            +
            \frac{24 d^2 \gamma^2 L (\delta^2 + 1) \tau}{m^2} \sum\limits_{t = 0}^T \expect{\norm{\nabla f(x^t)}^2} 
            +
            20 \sum\limits_{t = \tau}^T \frac{\gamma^2 L \tau d^2}{m^2} \sigma^2 .
        \end{equation*}

        Taking  $$\gamma \leq \frac{m^2}{48 d^2 L (\delta^2 + 1) \tau}, $$ we obtain

        \begin{equation}
        \label{eq:tmp_GD_even_5}
            \gamma \sum\limits_{t = 0}^{T-\tau} \expect{\norm{\nabla f(x^{t})}^2}
            \leq 
            8 F_\tau
            +
            \frac{48 d^2 \gamma^2 L (\delta^2 + 1) \tau}{m^2} \sum\limits_{t = T-\tau}^T \expect{\norm{\nabla f(x^t)}^2} 
            +
            40 \sum\limits_{t = \tau}^T \frac{\gamma^2 L \tau d^2}{m^2} \sigma^2 .
        \end{equation}

        We now prove that for any $t \geq 0$, we have $$\sup_{t \leq s \leq t + \tau} \left\{ \expect{\norm{\nabla f(x^s)}^2} \right\} \leq 4 \expect{\norm{\nabla f(x^t)}^2} + 8 L^2 \gamma^2 \tau^2 \frac{d^2}{m^2} \sigma^2 .$$

        For $t \leq s \leq t + \tau$ it holds that

        \begin{equation*}
        \begin{split}
            \expect{\norm{\nabla f(x^s)}^2}
            &\leq
            2 \expect{\norm{\nabla f(x^t)}^2}
            +
            \expect{\norm{\nabla f(x^s) - \nabla f(x^t)}^2}
            \\&\leq
            2 \expect{\norm{\nabla f(x^t)}^2} 
            +
            2 L^2 \gamma^2 \expect{\norm{\sum\limits_{r = t}^{s-1} \frac{1}{n} \sum\limits_{i=1}^n Q_r^i(\nabla f_i(x^r))}^2}
            \\&\leq
            2 \expect{\norm{\nabla f(x^t)}^2} 
            +
            2 L^2 \gamma^2 \tau \frac{d^2}{m^2} \sum\limits_{r = t}^{s-1} \frac{1}{n} \sum\limits_{i=1}^n \expect{\norm{\nabla f_i(x^r)}^2}
            \\&\leq
            2 \expect{\norm{\nabla f(x^t)}^2} 
            +
            4 L^2 \gamma^2 \tau \frac{d^2}{m^2} \sum\limits_{r = t}^{s-1}\left( (\delta^2 + 1) \expect{\norm{\nabla f(x^r)}^2} + \sigma^2\right)
            \\&\leq
            2 \expect{\norm{\nabla f(x^t)}^2} 
            +
            4 L^2 \gamma^2 \tau^2 \frac{d^2}{m^2} \left( (\delta^2 + 1) \hspace{-1mm}\sup_{t \leq s \leq t + \tau}\hspace{-1mm} \left\{ \expect{\norm{\nabla f(x^s)}^2} \right\} + \sigma^2\right) .
        \end{split}
        \end{equation*}

        Since $$\gamma \leq \frac{m}{\sqrt{8} d L \sqrt{\delta^2 + 1} \tau},$$

        it holds that 

        \begin{equation*}
            \sup_{t \leq s \leq t + \tau} \left\{ \expect{\norm{\nabla f(x^s)}^2} \right\} \leq
            4 \expect{\norm{\nabla f(x^t)}^2} 
            +
            8 L^2 \gamma^2 \tau^2 \frac{d^2}{m^2} \sigma^2 .
        \end{equation*}

        Using this \eqref{eq:tmp_GD_even_5} takes form 

        \begin{equation*}
        \begin{split}
            \gamma \sum\limits_{t = 0}^{T-\tau} \expect{\norm{\nabla f(x^{t})}^2}
            &\leq 
            8 F_\tau
            +
            \frac{192 d^2 \gamma^2 L (\delta^2 + 1) \tau}{m^2} \sum\limits_{t = T-2 \tau}^{T-\tau} \expect{\norm{\nabla f(x^t)}^2} 
            \\&+
            384 L^3 \gamma^4 \tau^3 \frac{d^4}{m^4} (\delta^2+1) \sigma^2
            +
            40 \sum\limits_{t = \tau}^T \frac{\gamma^2 L \tau d^2}{m^2} \sigma^2 .
        \end{split}
        \end{equation*}

        Taking $$\gamma \leq \frac{m}{384 d L \sqrt{\delta^2 + 1} \tau}, $$ and dividing both sides of the inequality by $T - \tau$, we obtain

        \begin{equation*}
            \frac{1}{T-\tau}\sum\limits_{t = 0}^{T-\tau} \expect{\norm{\nabla f(x^{t})}^2}
            \leq 
            16 \frac{F_\tau}{\gamma (T-\tau)}
            +
            80 \frac{\gamma^2 L \tau d^2}{m^2} \sigma^2 .
        \end{equation*}

        Therefore, if $\widehat{x}^T$ is chosen uniformly from $\left\{x^t \right\}_{t=0}^{T-1}$, then it holds that

        \begin{equation*}
            \expect{\norm{\nabla f(\widehat{x}^T)}^2} 
            \leq
            16 \frac{F_\tau}{\gamma T}
            +
            80 \frac{\gamma^2 L \tau d^2}{m^2} \sigma^2 .
        \end{equation*}
        
        This finishes the proof.
    \end{proof}


\subsection{Proof of Theorem \ref{theorem:GD_odd}, Under PL-condition}
\label{appendix_subsec:GD_PL}

    \begin{proof}
        We start from \eqref{eq:tmp_need_for_PL}:

        \begin{equation*}
        \begin{split}
            F_{t+1} - F_t 
            &=
            \frac{2 d^2 \gamma^2 L (\delta^2+1) \tau}{m^2} \expect{\norm{\nabla f(x^t)}^2}
            +
            \frac{2 d^2 \gamma^2 L (\delta^2 + 1)}{m^2} \sum\limits_{s = t - \tau}^{t-1} \expect{\norm{\nabla f(x^s)}^2}
            \\&+
            \left( \frac{\gamma^2 L \tau d^2}{2 m^2} - \frac{\gamma}{4} \right) \expect{\norm{\nabla f(x^{t - \tau})}^2}
            +
            \frac{\gamma d^2}{m^2}\left(4 \gamma L \tau + \varepsilon^2 \right)\sigma^2 .
        \end{split}
        \end{equation*}

        If $f$ satisfies PL-inequality (Assumption \ref{as:PL}), then $-\expect{\norm{\nabla f(x^{t-\tau}}^2} \leq -2 \mu F_{t-\tau}$, so that, for some $0 < \alpha < 1$ we obtain 

        \begin{equation}
        \begin{split}
        \label{eq:tmp_PL_1}
            F_{t+1} - F_t 
            &=
            \frac{2 d^2 \gamma^2 L (\delta^2+1) \tau}{m^2} \expect{\norm{\nabla f(x^t)}^2}
            +
            \frac{2 d^2 \gamma^2 L (\delta^2 + 1)}{m^2} \sum\limits_{s = t - \tau}^{t-1} \expect{\norm{\nabla f(x^s)}^2}
            \\&+
            \left( \frac{\gamma^2 L \tau d^2}{2 m^2} - \frac{(1 - \alpha) \gamma}{4} \right) \expect{\norm{\nabla f(x^{t - \tau})}^2}
            \\&-
            \frac{\alpha \gamma \mu}{2} F_{t-\tau}
            +
            \frac{\gamma d^2}{m^2}\left(4 \gamma L \tau + \varepsilon^2 \right)\sigma^2 .
        \end{split}
        \end{equation}

        For $t \geq 0$, let $p_t = p^t$ and $p = (1 - \alpha \mu \gamma / 4)^{-1}$. We multiply the above expression by $p_t$ and sum for $t < T$, hoping for cancellations. Using PL-condition (Assumption \ref{as:PL}), for $T \geq \tau$ we obtain

        \begin{equation*}
        \begin{split}
            \sum\limits_{t = \tau}^{T-1} p_{t+1} \left( F_t - F_{t+1} - \frac{\alpha \gamma \mu}{4} F_{t-\tau} \right)
            &=
            \sum\limits_{t = \tau}^{T-1} p_{t+1} \left[ \left(1 - \frac{\alpha \gamma \mu}{4} \right) F_t - F_{t+1} + \frac{\alpha \gamma \mu}{4} (F_t - F_{t-\tau}) \right]
            \\&=
            \sum\limits_{t = \tau}^{T-1} p_{t} F_t 
            -
            \sum\limits_{t = \tau+1}^{T} p_{t} F_t
            +
            \frac{\alpha \gamma \mu}{4} \sum\limits_{t = \tau}^{T-1} p_{t+1} (F_t - F_{t-\tau})
            \\&\leq
            p_\tau F_\tau - p_T F_T 
            +
            \frac{\alpha \gamma \mu}{4} \sum\limits_{t = \tau}^{T-1} p_{t+1} F_t
            \\&\quad-
            \frac{\alpha \gamma \mu p_\tau}{4} \sum\limits_{t = 0}^{T-1 - \tau} p_{t+1} F_t
            \\&\leq 
            p_\tau F_\tau - p_T F_T 
            +
            \frac{\alpha \gamma \mu}{4} \sum\limits_{t = T - \tau}^{T-1} p_{t+1} F_t
            \\&\leq
            p_\tau F_\tau - p_T F_T 
            +
            \frac{\alpha \gamma}{8} \sum\limits_{t = T - \tau}^{T-1} p_{t+1} \expect{\norm{\nabla f(x^t)}^2} .
        \end{split}
        \end{equation*}

        For any $t \geq 0$ we use a notation $b_t := \expect{\norm{\nabla f(x^t)}^2}$. We now handle $b_t$ terms from \eqref{eq:tmp_PL_1}.

        \begin{equation}
        \begin{split}
        \label{eq:tmp_PL_2}
            - \sum\limits_{t = \tau}^{T-1} \frac{(1-\alpha) \gamma}{4} p_{t+1} b_{t-\tau}
            +
            \gamma^2L\frac{d^2}{m^2} \sum\limits_{t = \tau}^{T-1} p_{t+1} \left(2 \tau (\delta^2 + 1) b_t + 2 (\delta^2 + 1) \sum\limits_{s = t-\tau}^{t-1} b_s + \frac{\tau}{2} b_{t-\tau}\right) .
        \end{split}
        \end{equation}

        If $p_t = p^t$, $p = (1 - \alpha \mu \gamma /2)^{-1}$ and $\gamma = \gamma_1 / \tau$, then, using the fact that 
        $(1 - a/x)^{-x} \leq 2 e^{a} \leq 2e$ if $x \geq 2$ and $0 \leq a \leq 1$, we can get that 
        $1 \geq p_\tau = (1 - \mu \gamma_1 /(2 \tau))^{-\tau} \leq 2 e^{\mu \gamma_1 /2} \leq 2e \leq 6$. Then
        
        $$\sum\limits_{t = \tau}^T p_{t+1} \sum\limits_{s=t-\tau}^{t-1} b_s \leq 
        p^{\tau}\sum\limits_{t = \tau}^T \sum\limits_{s=t-\tau}^{t-1} p_{s+1} b_s 
        \leq 
        6 \tau \sum\limits_{t = 0}^T p_{t+1} b_t .$$

        Now we can estimate \eqref{eq:tmp_PL_2}:

        \begin{equation*}
        \begin{split}
            \eqref{eq:tmp_PL_2}
            &\leq
            - \sum\limits_{t = 0}^{T- \tau - 1} \frac{(1-\alpha) \gamma}{4} p_{t+1} b_{t}
            +
            \gamma^2L\frac{d^2 \tau}{m^2} \left(
                2 (\delta^2+1) \sum\limits_{t = \tau}^{T - 1} b_t 
                +
                12 (\delta^2+1) \sum\limits_{t = 0}^{T - 1} b_t 
                +
                3 \sum\limits_{t = 0}^{T - \tau} b_t 
            \right)
            \\&\leq
            - \sum\limits_{t = 0}^{T- \tau - 1} p_{t+1} \gamma b_t
            \left(
                \frac{1-\alpha}{4}
                -
                17 \gamma L\frac{d^2 \tau (\delta^2+1)}{m^2} 
            \right)
            +
            14 \gamma^2 L\frac{d^2 \tau (\delta^2+1)}{m^2} \sum\limits_{t = T - \tau}^{T - 1} p_{t+1} b_t .
        \end{split}
        \end{equation*}

        Taking $$\gamma \leq \frac{m^2(1-\alpha)}{136 L d^2 \tau (\delta^2+1) \beta}, $$

        where $\beta \geq 1$, we obtain

        \begin{equation*}
            \eqref{eq:tmp_PL_2}
            \leq 
            - \frac{(1-\alpha)\gamma}{8} \sum\limits_{t = 0}^{T- \tau - 1} p_{t+1} b_t
            +
            \frac{(1-\alpha) \gamma}{4 \beta} \sum\limits_{t = T - \tau}^{T - 1} p_{t+1} b_t .
        \end{equation*}

        Now we can estimate \eqref{eq:tmp_PL_1}:

        \begin{equation}
        \begin{split}
        \label{eq:tmp_PL_3}
            0
            \leq
            p_\tau F_\tau - p_T F_T 
            &+
            \left(\frac{\alpha \gamma}{8} + \frac{(1-\alpha) \gamma}{4 \beta} \right) \sum\limits_{t = T - \tau}^{T-1} p_{t+1} b_t
            -
            \frac{(1-\alpha)\gamma}{8} \sum\limits_{t = 0}^{T- \tau - 1} p_{t+1} b_t
            \\&+
            \sum\limits_{t = \tau}^{T - 1}p_{t+1} \frac{\gamma d^2}{m^2}\left(4 \gamma L \tau + \varepsilon^2 \right)\sigma^2 .
        \end{split}
        \end{equation}

        Using that we proved in \ref{appendix_subsec:non_conv} we have $b_t \leq 4 b_{t-\tau} + 8 L^2 \gamma^2 \tau^2 \frac{d^2}{m^2} \sigma^2$. Then, we can obtain 

        \begin{equation*}
        \begin{split}
            \gamma \left(\frac{\alpha}{8} + \frac{1-\alpha}{4 \beta} \right) \sum\limits_{t = T - \tau}^{T-1} p_{t+1} b_t
            &\leq
            24 \gamma \left(\frac{\alpha}{8} + \frac{1-\alpha}{4 \beta} \right) \sum\limits_{t = T - 2 \tau}^{T- \tau - 1} p_{t+1} b_t
            \\&+
            48 L^2 \gamma^3 \tau^3 \frac{d^2}{m^2} \left(\frac{\alpha}{8} + \frac{1-\alpha}{4 \beta} \right) \sigma^2 .
        \end{split}
        \end{equation*}

        Taking $\alpha = 1/6$ and $\beta = 4$, we obtain

        $$\frac{\alpha}{8} + \frac{1-\alpha}{4 \beta} = \frac{1 - \alpha}{8},$$

        and \eqref{eq:tmp_PL_3} takes form

        \begin{equation}
        \label{eq:tmp_PL_3.5}
        \begin{split}
            0
            \leq
            p_\tau F_\tau - p_T F_T 
            &+
            48 L^2 \gamma^3 \tau^3 \frac{d^2}{m^2} \sigma^2
            +
            \sum\limits_{t = \tau}^{T - 1}p_{t+1} \frac{\gamma d^2}{m^2}\left(4 \gamma L \tau + \varepsilon^2 \right)\sigma^2 .
        \end{split}
        \end{equation}

         Using the fact that 

        \begin{equation*}
            \sum\limits_{t = \tau}^T \left(1 - \frac{\alpha \mu \gamma}{2} \right)^{T -t} 
            =
            \sum\limits_{t = 0}^{T-\tau} \left(1 - \frac{\alpha \mu \gamma}{2} \right)^{t}
            \leq 
            \sum\limits_{t = 0}^{+\infty} \left(1 - \frac{\alpha \mu \gamma}{2} \right)^{t}
            =
            \frac{2}{\alpha \mu \gamma},
        \end{equation*}

        and  taking 

        $$\gamma \leq \frac{m^2}{625 L d^2 \tau (\delta^2+1)}
        ~~\text{ and }~~
        \varepsilon =  \sqrt{\gamma L \tau} \leq \frac{m}{25 d \sqrt{\delta^2+1}}, $$

        by dividing \eqref{eq:tmp_PL_3.5} by $p_\tau$, we obtain 

        \begin{equation*}
            \expect{f(x^T) - f(x^*)} 
            \leq
            \left( 1 - \frac{\mu \gamma}{12} \right)^{T-\tau} \expect{f(x^\tau) - f(x^*)} 
            +
            636 \frac{\gamma d^2 L \tau}{ \mu m^2}\sigma^2 .
        \end{equation*}

        This finishes the proof.

    \end{proof}

\section{Convergence of Algorithm \ref{alg:GD} without data similarity}
\label{appendix:proof_GD_sem}

    \begin{theorem}[Convergence of GD Algorithm \ref{alg:GD} without data similarity]
    \label{theorem:GD_sem}
        Consider Assumptions \ref{as:lip} and \ref{as:strconv}. Let problem \eqref{eq:problem} be solved by Algorithm \ref{alg:GD}. Then for any $\varepsilon > 0$, $\gamma > 0$, $\tau > \tau_{\text{mix}}(\varepsilon)$ and $T > \tau$ satisfying

        \begin{equation*}
            \gamma \leq \frac{m^2 \sqrt{\mu}}{24 d^2 L^{3/2} \tau}
            ~~\text{ and }~~
            \varepsilon
            \leq 
            \frac{m \sqrt{\mu}}{24 d} \min\left\{ \frac{1}{L^{3/2}} ; \sqrt{\mu} \right\},
        \end{equation*}

        it holds that
        
        $$\expect{\norms{x^{T+1} - x^*}^2}
            \leq
            \left(1 - \frac{\mu \gamma}{2} \right)^{T -\tau} \expect{\norms{x^{\tau} - x^*}^2}
            +
            \left(1 - \frac{\mu \gamma}{2} \right)^{T} \Delta_{\tau}
            +
            26 \frac{\gamma d^2 \tau}{\mu m^2}  \sigma_*^2 ,$$
    where
        $$\Delta_{\tau} = \mathcal{O} \Bigg( \frac{\gamma^2 d^2}{m^2} \sqrt{\frac{\mu}{L}} \sum\limits_{t = 0}^{\tau} \Big[\tau\expect{\norms{x^{t} - x^*}^2}+ 4L\expect{f(x^t) - f(x^*)} \Big]\Bigg).$$
    
    \end{theorem}
    
    \begin{proof}[Proof of Theorem \ref{theorem:GD_sem}]
        We start by writing out step of the Algorithm \ref{alg:GD}:
    
        \begin{equation}
        \begin{split}
        \label{eq:tmp_GD_01}
            \expect{\norms{x^{t+1} - x^*}^2} 
            &=
            \expect{\norms{x^t - x^*}^2}
            -
            2 \gamma \expect{\frac{1}{n} \sum\limits_{i=1}^d \dotprod{Q_t^i\left(\nabla f_i(x^{t})\right) - \nabla f_i(x^{t})}{x^{t} - x^*}}
            \\&-
            2 \gamma \expect{\dotprod{\nabla f(x^t)}{x^t - x^*}}
            + \gamma^2 \expect{\norms{\frac{1}{n} \sum\limits_{i=1}^n Q_t^i\left(\nabla f_i(x^t)\right)}^2} .
        \end{split}
        \end{equation}

        Consider $\expect{ \dotprod{Q_t^i\left(\nabla f_i(x^{t})\right) - \nabla f_i(x^{t})}{x^{t} - x^*}}$. Using Corollary \ref{cor:main_dotprod} with $a^t = \nabla f_i(x^t)$, $b^t = x^t - x^*$, $\hat{a}^{t-\tau} = \nabla f_i(x^{t - \tau})$ and $\hat{b}^{t-tau} = x^{t-\tau} - x^*$ we obtain
        \begin{equation}
        \begin{split}
        \label{eq:tmp_GD_0}
            &2 \E[\frac{1}{n} \sum\limits_{i=1}^n |\langle Q_t^i\left(\nabla f_i(x^{t})\right) - \nabla f_i(x^{t}), x^{t} - x^*\rangle|]
            \leq
            \frac{\varepsilon d}{m \beta_0} \frac{1}{n} \sum\limits_{i=1}^n \expect{\norm{\nabla f_i(x^{t-\tau})}^2} 
            \\&+
            \frac{\varepsilon d \beta_0}{m} \expect{\norms{x^{t-\tau} - x^*}^2}
            +
            4 \frac{d^2L^2}{m^2}(\beta_1 + \beta_2) \expect{\norms{x^t - x^{\tau}}^2} 
            +
            \left( \frac{1}{\beta_1} + \frac{1}{\beta_3}  \right) \expect{\norms{x^t - x^{\tau}}^2} 
            \\&+
            4 \frac{d^2}{m^2} \beta_3 \frac{1}{n} \sum\limits_{i=1}^d \expect{\norm{\nabla f_i(x^{t})}^2}
            +
            \frac{1}{\beta_2} \expect{\norms{x^{t} - x^*}^2} .
        \end{split}
        \end{equation}
        Using the fact that $f_i$ are $L$-smooth, we can obtain:
        \begin{equation}
        \begin{split}
        \label{eq:super_lemma}
            \frac{1}{n} \sum\limits_{i = 1}^{n} \norm{\nabla f_i(x^t)}^2
            &=
            \frac{1}{n} \sum\limits_{i = 1}^{n} \norm{\nabla f_i(x^t) - \nabla f_i(x^*) + \nabla f_i(x^*)}^2
            \\&\leq
            \frac{2}{n} \sum\limits_{i = 1}^{n} \norm{\nabla f_i(x^t) - \nabla f_i(x^*)}^2
            +
            \frac{2}{n} \sum\limits_{i = 1}^{n} \norm{\nabla f_i(x^*)}^2
            \\&\leq
            \frac{4L}{n} \sum\limits_{i = 1}^{n} \left(f_i(x^t) -f_i(x^*) - \dotprod{\nabla f_i(x^*)}{x^t - x^*} \right) + 2 \sigma_*^2
            \\&=
            4L (f(x^t) - f(x^*)) + 2 \sigma_*^2,
        \end{split}
        \end{equation}

        where we use a notation $\sigma_*^2 := \frac{1}{n}\sum_{i = 1}^{n}\norm{\nabla f_i(x^*)}^2$. Now we can estimate \eqref{eq:tmp_GD_0}:

        \begin{equation*}
        \begin{split}
            \eqref{eq:tmp_GD_0}
            &\leq
            \frac{2 \varepsilon d}{m \beta_0} (2L \expect{f(x^{t-\tau}) - f(x^*)} +  \sigma_*^2)
            +
            \frac{\varepsilon d \beta_0}{m} \expect{\norms{x^{t-\tau} - x^*}^2}
            \\&+
            \left(4 \frac{d^2L^2}{m^2}(\beta_1 + \beta_2) 
            +
            \frac{1}{\beta_1} + \frac{1}{\beta_3}  \right) \expect{\norms{- \gamma \sum\limits_{s = t-\tau}^{t-1} \frac{1}{n} \sum\limits_{i = 1}^n Q_s^i\left( \nabla f_i(x^s)\right)}^2} 
            \\&+
            8 \frac{d^2}{m^2} \beta_3 (2L \expect{f(x^t) - f(x^*)} + \sigma_*^2)
            +
            \frac{1}{\beta_2} \expect{\norms{x^{t} - x^*}^2} .
        \end{split}
        \end{equation*}

        Now we can estimate \eqref{eq:tmp_GD_01}. Using Lemma \ref{lemma:main_norms} and Assumption \ref{as:strconv} we can obtain
        \begin{equation}
        \begin{split}
        \label{eq:tmp_GD_02}
            \mathbb{E}\Big[\Big\|x^{t+1} &- x^*\Big\|^2\Big]
            \leq
            \left(1 - \mu \gamma + \frac{\gamma}{\beta_2} \right) \expect{\norms{x^{t} - x^*}^2} 
            +
            \frac{\varepsilon d \beta_0 \gamma}{m} \expect{\norms{x^{t-\tau} - x^*}^2} 
            \\&+
            4L \mathbb{E}\Bigg[
                \frac{\varepsilon d \gamma}{m \beta_0} (f(x^{t-\tau}) - f(x^*))
                +
                4 \frac{d^2 \beta_3 \gamma}{m^2} (f(x^t) - f(x^*))
                \\&+
                \left(4 \frac{d^2L^2}{m^2}(\beta_1 + \beta_2) +
                \frac{1}{\beta_1} + \frac{1}{\beta_3}  \right) \frac{\gamma^3 \tau d^2}{m^2} \sum\limits_{s = t-\tau}^{t-1}(f(x^s) - f(x^*))
                \\&+
                \frac{\gamma^2d^2}{m^2} (f(x^t) - f(x^*)) - \frac{\gamma}{2L} (f(x^t) - f(x^*))
            \Bigg]
            \\&+
            2 \Bigg[
                \frac{\varepsilon d \gamma}{m \beta_0}
                +
                4 \frac{d^2 \beta_3 \gamma}{m^2}  
                +
                \left(4 \frac{d^2L^2}{m^2}(\beta_1 + \beta_2) +
                \frac{1}{\beta_1} + \frac{1}{\beta_3}  \right) \frac{\gamma^3 \tau^2 d^2}{m^2}
                +
                \frac{\gamma^2d^2}{m^2}
            \Bigg] \sigma_*^2 .
        \end{split}
        \end{equation}

        Taking  $\beta_0 = \beta_1 = 1, \beta_3 = \gamma, \beta_2 = 4/\mu$ and using fact, that $\varepsilon \leq \gamma \tau d/m$ inequality \eqref{eq:tmp_GD_02} takes form
        \begin{equation}
        \begin{split}
        \label{eq:tmp_GD_3}
            \expect{\norms{x^{t+1} - x^*}^2} 
            &\leq
            \left(1 - \frac{3}{4 \mu \gamma} \right) \expect{\norms{x^{t} - x^*}^2} 
            +
            \frac{\varepsilon d \beta_0 \gamma}{m} \expect{\norms{x^{t-\tau} - x^*}^2} 
            \\&+
            4L \mathbb{E}\Bigg[
                \frac{\varepsilon d \gamma}{m \beta_0} (f(x^{t-\tau}) - f(x^*))
                +
                5 \frac{d^2 \gamma^2}{m^2} (f(x^t) - f(x^*))
                \\&+
                20 \frac{d^4L^2}{m^4} \frac{\gamma^3 \tau}{\mu} \sum\limits_{s = t-\tau}^{t-1}(f(x^s) - f(x^*))
                - 
                \frac{\gamma}{2L} (f(x^t) - f(x^*))
            \Bigg]
            \\&+
            4 \frac{d^2 \gamma^2 \tau}{m^2} \Bigg[
                3
                +
                10 \frac{d^2L^2}{m^2} \frac{\gamma}{\mu}
            \Bigg] \sigma_*^2 .
        \end{split}
        \end{equation}
        
        Let us perform the summation from $t = \tau$ to $t = T > \tau$ of equations \eqref{eq:tmp_GD_3} with coefficients $p_k$:

        \begin{align}
        \label{eq:tmp_GD_4}
            \sum\limits_{t = \tau}^T p_t \expect{\norms{x^{t+1} - x^*}^2}
            &\leq
            \sum\limits_{t = \tau}^T p_t (1 - \frac{3 \mu \gamma}{4}) \expect{\norms{x^t - x^*}^2}
            \\&+
            \sum\limits_{t = \tau}^T p_t \frac{\gamma \varepsilon d}{m} \expect{\norms{x^{t-\tau} - x^*}^2}
            \\&+
            \sum\limits_{t = \tau}^T p_t 4L \left( \frac{\gamma \varepsilon d}{m} + 5 \frac{\gamma^2 d^2 \tau}{m^2} - \frac{\gamma}{2 L} \right) \expect{f(x^t) - f(x^*)} 
            \\&+
            20 \sum\limits_{t = \tau}^T p_t 4 L \frac{d^4L^2}{m^4} \frac{\gamma^3 \tau}{\mu} \sum\limits_{s=t-\tau}^{t-1} \expect{f(x^s) - f(x^*)}
            \\&+
            \sum\limits_{t = \tau}^T p_t 4 \frac{d^2 \gamma^2 \tau}{m^2} \Bigg[
                3
                +
                10 \frac{d^2L^2}{m^2} \frac{\gamma}{\mu}
            \Bigg] \sigma_*^2 .
        \end{align}

        If $p_t = p^t$, $p = (1 - \mu \gamma /2)^{-1}$ and $\gamma = \gamma_1 / \tau$, then, using the fact that 
        $(1 - a/x)^{-x} \leq 2 e^{a} \leq 2e$ if $x \geq 2$ and $0 \leq a \leq 1$, we can get that 
        $p_\tau = (1 - \mu \gamma_1 /(2 \tau))^{-\tau} \leq 2 e^{\mu \gamma_1 /2} \leq 2e \leq 6$.
        
        $$\sum\limits_{t = \tau}^T p_t \sum\limits_{s=t-\tau}^{t-1} a_s \leq 
        p^{\tau}\sum\limits_{t = \tau}^T \sum\limits_{s=t-\tau}^{t-1} p_s a_s 
        \leq 
        6 \tau \sum\limits_{t = 0}^T p_t a_t .$$

        Using this we can estimate \eqref{eq:tmp_GD_4}:
        \begin{equation}
        \begin{split}
        \label{eq:tmp_GD_5}
            &\sum\limits_{t = \tau}^T p_t \expect{\norms{x^{t+1} - x^*}^2}
            \leq
            \sum\limits_{t = \tau}^T p_t \left(1 - \mu \gamma + 6 \frac{\gamma \varepsilon d}{m}\right) \expect{\norms{x^t - x^*}^2}
            \\&+
            \sum\limits_{t = \tau}^T 4 p_t L \left( \frac{\gamma \varepsilon d}{m} + 5 \frac{\gamma^2 d^2 \tau}{m^2}
            +
            120 \frac{d^4L^2}{m^4} \frac{\gamma^3 \tau^2}{\mu}
            - 
            \frac{\gamma}{2 L} \right) \expect{f(x^t) - f(x^*)}
            \\&+
            4 \sum\limits_{t = \tau}^T p_t \Bigg[
                3
                +
                10 \frac{d^2L^2}{m^2} \frac{\gamma}{\mu}
            \Bigg] \sigma_*^2 
            +
            \sum\limits_{t = 0}^{\tau} p_{t+\tau} \frac{\gamma \varepsilon d}{m} \expect{\norms{x^{t} - x^*}^2}
            \\&+
            80 \sum\limits_{t = 0}^{\tau} p_{t+\tau} L \frac{d^4L^2}{m^4} \frac{\gamma^3 \tau}{\mu} \expect{f(x^t) - f(x^*)} .
        \end{split}
        \end{equation}

        Taking

        \begin{equation*}
            \gamma \leq \frac{m^2 \sqrt{\mu}}{24 d^2 L^{3/2} \tau}
            ~~\text{ and }~~
            \varepsilon = \min\left\{\frac{\gamma d \tau}{m} ; \frac{\mu m }{24 d } \right\}
            \leq 
            \frac{m \sqrt{\mu}}{24 d} \min\left\{ \frac{1}{L^{3/2}} ; \sqrt{\mu} \right\}.
        \end{equation*}

        We get

        \begin{equation*}
            \frac{\gamma \varepsilon d}{m} + 5 \frac{\gamma^2 d^2 \tau}{m^2}
            +
            120 \frac{d^4L^2}{m^4} \frac{\gamma^3 \tau^2}{\mu}
            - 
            \frac{\gamma}{2 L} \leq 0
            ~~\text{ and }~~
            1 - \frac{3 \mu \gamma}{4} + 6 \frac{\gamma \varepsilon d}{m} = 1 - \frac{\mu \gamma}{2} .
        \end{equation*}

        Assume a notation 
        
        \begin{equation*}
        \begin{split}
            \Delta_{\tau} 
            &:=
            \sum\limits_{t = 0}^{\tau} p_{t+\tau} \frac{\gamma \varepsilon d}{m} \expect{\norms{x^{t} - x^*}^2}
            +
            80 \sum\limits_{t = 0}^{\tau} p_{t+\tau} L \frac{d^4L^2}{m^4} \frac{\gamma^3 \tau}{\mu} \expect{f(x^t) - f(x^*)}
            \\&~\leq
            120 \frac{\gamma^2 d^2}{m^2} \sqrt{\frac{\mu}{L}}
            \sum\limits_{t = 0}^{\tau} \left( \tau \expect{\norms{x^{t} - x^*}^2} + 4 L  \expect{f(x^t) - f(x^*)} \right). 
        \end{split}
        \end{equation*}

        Using the notation of $\Delta_{\tau}$, \eqref{eq:tmp_GD_5} takes form

        \begin{equation*}
            \sum\limits_{t = \tau}^T p_t \expect{\norms{x^{t+1} - x^*}^2}
            \leq
            \sum\limits_{t = \tau}^T p_t \left(1 - \frac{\mu \gamma}{2} \right) \expect{\norms{x^t - x^*}^2}
            +
            \sum\limits_{t = \tau}^T 13 p_t \frac{\gamma^2 d^2 \tau}{m^2}  \sigma_*^2 
            +
            \Delta_{\tau} .
        \end{equation*}

        Using $p_t = p^t$ and $p = \left(1 - \mu \gamma / 2 \right)^{-1}$ we can obtain:

        \begin{equation*}
        \begin{split}
            \sum\limits_{t = \tau}^T \left(1 - \frac{\mu \gamma}{2} \right)^{-t} \expect{\norms{x^{t+1} - x^*}^2}
            &\leq
            \sum\limits_{t = \tau}^T \left(1 - \frac{\mu \gamma}{2} \right)^{-t+1} \expect{\norms{x^t - x^*}^2}
            \\&+
            \sum\limits_{t = \tau}^T 13 \left(1 - \frac{\mu \gamma}{2} \right)^{-t} \frac{\gamma^2 d^2 \tau}{m^2}  \sigma_*^2 
            +
            \Delta_{\tau} .
        \end{split}
        \end{equation*}

        The summed terms on the left and right sides are reduced, therefore this expression takes the form:

        \begin{equation*}
        \begin{split}
            \left(1 - \frac{\mu \gamma}{2} \right)^{-T} \expect{\norms{x^{T+1} - x^*}^2}
            &\leq
            \left(1 - \frac{\mu \gamma}{2} \right)^{-\tau} \expect{\norms{x^{\tau} - x^*}^2}
            \\&+
            \sum\limits_{t = \tau}^T 13 \left(1 - \frac{\mu \gamma}{2} \right)^{-t} \frac{\gamma^2 d^2 \tau}{m^2}  \sigma_*^2 
            +
            \Delta_{\tau} .
        \end{split}
        \end{equation*}

        We can re-arrange this inequality:

        \begin{equation*}
        \begin{split}
            \expect{\norms{x^{T+1} - x^*}^2}
            &\leq
            \left(1 - \frac{\mu \gamma}{2} \right)^{T -\tau} \expect{\norms{x^{\tau} - x^*}^2}
            \\&+
            \sum\limits_{t = \tau}^T 13 \left(1 - \frac{\mu \gamma}{2} \right)^{T-t} \frac{\gamma^2 d^2 \tau}{m^2}  \sigma_*^2 
            +
            \left(1 - \frac{\mu \gamma}{2} \right)^{T} \Delta_{\tau} .
        \end{split}
        \end{equation*}

        Using the fact that 

        \begin{equation*}
            \sum\limits_{t = \tau}^T \left(1 - \frac{\mu \gamma}{2} \right)^{T -t} 
            =
            \sum\limits_{t = 0}^{T-\tau} \left(1 - \frac{\mu \gamma}{2} \right)^{t}
            \leq 
            \sum\limits_{t = 0}^{+\infty} \left(1 - \frac{\mu \gamma}{2} \right)^{t}
            =
            \frac{2}{\mu \gamma} .
        \end{equation*}

        We can estimate:

        \begin{equation*}
            \expect{\norms{x^{T+1} - x^*}^2}
            \leq
            \left(1 - \frac{\mu \gamma}{2} \right)^{T -\tau} \expect{\norms{x^{\tau} - x^*}^2}
            +
            \left(1 - \frac{\mu \gamma}{2} \right)^{T} \Delta_{\tau}
            +
            26 \frac{\gamma d^2 \tau}{\mu m^2}  \sigma_*^2 .
        \end{equation*}

        This finishes the proof.

    \end{proof}

\section{Extensions for Theorem \ref{theorem:GD_acc}}
\label{appendix_subsec:GD_ACC}

\subsection{Full version of Theorem \ref{theorem:GD_acc}}
\label{appendix_subsec:full_GD_acc}

\begin{theorem}[Convergence of AMQSGD Algorithm \ref{alg:acc}, full version]
        \label{theorem:full_GD_acc}
        Consider Assumptions \ref{as:lip}, \ref{as:strconv} and \ref{as:sim}. Let problem \eqref{eq:problem} be solved by Algorithm \ref{alg:acc}. Then for  any $\gamma > 0, \varepsilon > 0$, $\tau > \tau_{\text{mix}}(\varepsilon)$, $T > \tau$ and $\beta, \theta, \eta, p$ satisfying
        
        $$
        \gamma\lesssim\frac{\mu^{\frac{1}{3}}m^{\frac{1}{2}}}{\tau L^{\frac{4}{3}}d^{\frac{1}{2}}}, 
        ~~~~
        p\lesssim\frac{m^2}{\tau^2d^2(\delta^2+1)},
        ~~~~
        \e\lesssim\min\Big\{\frac{m^{\frac{7}{4}}}{d^{\frac{7}{4}}\tau^{\frac{5}{4}}L(\delta^2 + 1)}; \frac{m^{\frac{15}{4}}}{d^{\frac{15}{4}}\tau^{\frac{13}{4}}(\delta^2 + 1)^2}\Big\}
        $$
        
        $$
        \beta = \sqrt{\frac{2p^2\mu\gamma}{3}},
        ~~~~
        \eta = \sqrt{\frac{3}{2\mu\gamma}},
        ~~~~
        \theta = \frac{p\eta^{-1} - 1}{\beta p \eta^{-1} - 1}
        $$
        
        it holds that
        \begin{equation*}
            \begin{split}
                F_{T+1}
                &=
                \mathcal{O} \Bigg(
                \exp\left[-(T - \tau)\sqrt{\frac{p^2\mu\gamma}{3}} \right]F_{\tau}
                + 
                \exp\left[-T\sqrt{\frac{p^2\mu\gamma}{3}} \right]\Delta_{\tau} 
                + 
                \frac{\gamma}{\mu}\sigma^2
                \Bigg).
            \end{split}
        \end{equation*}
        
        Here we use notations: $F_t := \E[\|x^{t} - x^*\|^2 + \frac{3}{\mu}( f(x_f^{t}) - f(x^*))]$ and $\Delta_{\tau} \leq \frac{\sqrt{\gamma}}{\tau^{\frac{4}{3}}\mu^{\frac{1}{3}}}\sum\limits_{t=0}^{\tau}\big(\E\norm{\nabla f(x_g^t)} + \E\norm{x^{t} - x^*}^2 + \E[f(x_f^{t}) - f(x^*)] \big)$.
    \end{theorem}

\subsection{Full version of Corollary \ref{corollary:GD_acc}}
\label{appendix:corollary_GD_acc_full}
\begin{corollary}[Step tuning for Theorem \ref{theorem:GD_acc}, full version of Corollary \ref{corollary:GD_acc}]
            Under the conditions of Theorem \ref{theorem:GD_acc}, choosing $\gamma$ as
            
            $$\gamma \lesssim \min\left\{
                                    \frac{\mu^{\frac{1}{3}}}{L^{\frac{4}{3}} \tau^{\frac{8}{3}}}
                                    ~;~
                                    \frac{\log\left( \max\left\{ 2 ;
                                    \frac{\mu^\frac{2}{3} (F_\tau + \Delta_\tau) T}{ \tau^\frac{4}{3}L^\frac{2}{3} \sigma^2}\right\} \right)}{\mu p^2T^2}
                                    \right\},$$
            
            in order to achieve $\epsilon$-approximate solution (in terms of $\expect{\norm{x^T - x^*}^2} \leq \epsilon^2$) it takes

            \begin{equation*}
                \mathcal{O}\left( \frac{d^2 L^{\frac{2}{3}} \tau^{\frac{4}{3}}}{m^2 \mu^{\frac{2}{3}}} \left(
                (\delta^2+1)\log\left(\frac{1}{\epsilon}\right) + \frac{\sigma^2}{\mu \epsilon} 
                \right)\right) \text{ iterations} .
            \end{equation*}
            
    \end{corollary}

\subsection{Proof of Theorem \ref{theorem:full_GD_acc}}
\label{appendix_subsec:full_GD_acc_proof}


    \begin{lemma}
    \label{lemma:GD_acc_tilde}
        Consider Algorithm \ref{alg:acc} with $\theta = (p \eta^{-1} - 1) / (\beta \eta^-1 - 1) < 1$. Then for any $y^t = \kappa x_f^t + (1 - \kappa) x^t \in \text{conv}\left\{x_f^t, x^t\right\}$ for any $s < t$ exist constants $\alpha_f^s, \alpha^s \geq 0$ and $c_r \geq 0$ such that 

        \begin{equation*}
            y^t 
            =
            \widetilde{y}^s - p\gamma \sum\limits_{r = s}^{t-1} c_r g^r
            = \alpha_f^s x_f^s + \alpha^s x^s - p\gamma \sum\limits_{r = s}^{t-1} c_r g^r.
        \end{equation*}

        And $\alpha_f^s + \alpha^s = 1$ for any $s < t$. If $(1 - \kappa) \eta \leq 1$, then $c_r \leq t - s + 2$, otherwise we can only use the estimate $c_r \leq \eta$.
    \end{lemma}

    \begin{proof}
        We start by writing out lines 3 and 10 of Algorithm \ref{alg:acc}:

        \begin{equation}
        \label{eq:tmp_acc_tilde_1}
            x_f^s 
            = 
            x_g^{s-1} - p\gamma g^{s-1} 
            =
            \theta x_f^{s-1} + (1 - \theta) x^{s-1} - p \gamma g^{s-1} .
        \end{equation}

        Now let us handle expression $\eta x^k_g + (p - \eta)x^k_f + (1- p)(1 - \beta) x^k +(1 - p) \beta x^k_g  - x^*$ for a while. Taking into account the choice of $\theta$ such that $\theta = (p \eta^{-1} - 1) / (\beta p \eta^{-1} - 1)$ (in particular, $(p \eta^{-1} - 1) = (\beta p \eta^{-1} - 1)\theta$ and $ \eta (1 - \beta p \eta^{-1}) (1 - \theta) = p (1 - \beta)$), we get
        
        \begin{align*}
        \eta x^k_g + (p - \eta) & x^k_f + (1- p)(1 - \beta) x^k +(1 - p) \beta x^k_g
        \notag\\
        &=
        (\eta + (1 - p) \beta) x^k_g + (p - \eta)x^k_f + (1- p)(1 - \beta) x^k
        \notag\\
        &=
        (\eta + (1 - p) \beta) x^k_g + \eta (p \eta^{-1} - 1)x^k_f + (1- p)(1 - \beta) x^k
        \notag\\
        &=
        (\eta + (1 - p) \beta) x^k_g + \eta (\beta p \eta^{-1} - 1)\theta x^k_f + (1- p)(1 - \beta) x^k
        \notag\\
        &=
        (\eta + (1 - p) \beta) x^k_g + \eta (\beta p \eta^{-1} - 1)(x^k_g - (1 - \theta) x^k) + (1- p)(1 - \beta) x^k
        \notag\\
        &=
        (\eta + (1 - p) \beta) x^k_g + \eta (\beta p \eta^{-1} - 1)(x^k_g - (1 - \theta) x^k) + (1- p)(1 - \beta) x^k
        \notag\\
        &=
        \beta x^k_g - \eta (\beta p \eta^{-1} - 1)(1 - \theta) x^k + (1- p)(1 - \beta) x^k
        \notag\\
        &=
        \beta x^k_g + p(1 - \beta) x^k + (1- p)(1 - \beta) x^k
        \notag\\
        &=
        \beta x^k_g + (1 - \beta) x^k
        \,.
        \end{align*}

        Now we write out line \ref{line_x_acc} of Algorithm \ref{alg:acc}:

        \begin{equation}
        \begin{split}
        \label{eq:tmp_acc_tilde_2}
            x^s 
            &= 
            \beta x_g^{s-1} + (1 - \beta) x^{s-1} - \eta x_g^{s-1} + \eta x_f^s 
            = 
            \beta x_g^{s-1} + (1 - \beta) x^{s-1} - \eta p \gamma g^{s-1} 
            \\&=
            \beta (\theta x_f^{s-1} + (1 - \theta) x^{s-1}) + (1 - \beta) x^{s-1} - \eta p \gamma g^{s-1} 
            \\&=
            \beta \theta x_f^{s-1} + (1 - \beta \theta) x^{s-1} - \eta p \gamma g^{s-1} .
        \end{split}
        \end{equation}

        Now we use induction. $x_f^t = \theta x_f^{s-1} + (1 - \theta) x^{s-1} - p \gamma g^{s-1}$, then $\alpha_f^{t-1} = \theta \geq 0$, $\alpha^{t-1} = 1 - \theta \geq 0$, $c_r = 1 \leq \eta$ and $\alpha_f^{t-1} + \alpha^{t-1} = 1$, therefore base step is fulfilled. If $x_f^t = \alpha_f^{s} x_f^s + \alpha^s x^s - p \gamma \sum_{r = s}^{t-1} c_r g^r$ for some $s < t$, when with help of \eqref{eq:tmp_acc_tilde_1} and \eqref{eq:tmp_acc_tilde_2} we can write out

        \begin{equation*}
        \begin{split}
            x_f^t 
            &= \alpha_f^s \left(\theta x_f^{s-1} + (1 - \theta) x^{s-1} - p \gamma g^{s-1} \right)
            \\&+
            \alpha^s \left( \beta \theta x_f^{s-1} + (1 - \beta \theta) x^{s-1} - \eta p \gamma g^{s-1} \right)
            -
            p \gamma \sum\limits_{r = s}^{t-1} c_r g^r .
        \end{split}
        \end{equation*}

        Therefore $\alpha_f^{s-1} = \alpha_f^s \theta + \alpha^s \beta \theta \geq 0$, $\alpha^{s-1} = \alpha_f^s (1 - \theta) + \alpha^s (1 - \beta \theta) \geq 0$ and $c_{s-1} = \alpha_f^s + \eta \alpha^s \leq \eta$. Then, the step of the induction is fulfilled, since $\alpha_f^{s-1} + \alpha^{s-1} = 1$. Therefore results of this Lemma are true for $y^t = x_f^t \in \text{conv}\left\{ x_f^t, x^t\right\}$.

        Consider $y^t = x^t \in \text{conv}\left\{ x_f^t, x^t\right\}$. Form \eqref{eq:tmp_acc_tilde_2} follows that $\alpha_f^{t-1} = \beta \theta$ and $\alpha^{t-1} = 1 - \beta \theta$, therefore base step is fulfilled. The step of the induction will be the same as in $y^t = x_f^t$. Therefore results of this Lemma are true for $y^t = x^t$. Then, they are true for any $y^t \in \text{conv}\left\{ x_f^t, x^t\right\}$.

        If $y^t = \kappa x_f^t + (1 - \kappa) x^t$, then $\alpha^s(y) = \kappa \alpha^s(x_f^t) + (1 - \kappa) \alpha^s(x^t)$. Since $(1 - \theta) \eta \leq 1$, then $\alpha^{t-1}(x_f^t) \eta \leq 1 = t - (t - 1)$. Therefore $\alpha^{s}(x_f^t) \eta \leq t - s$ by induction, since $\alpha^{s-1}(x_f^t) \eta = \alpha_f^s(x_f^t) (1-\theta) \eta + (1 - \beta \theta) \alpha^s(x_f^t) \eta \leq \alpha_f^s(x_f^t) + (1 - \beta \theta) (t - s) \leq t - s + 1$. 
        
        Then, if $(1 - \kappa)\eta \leq 1$, then $\alpha^s(y^t) \eta = \kappa \alpha^s(x_f^t) \eta + (1 - \kappa) \eta \alpha^s(x^t) \leq \kappa (t - s) + \alpha^s(x^t) \leq t - s + 1$. Now we consider $c_s(y^t)$. $c_s(y^t) = \alpha_f^s(y^t) + \alpha^s(y^t) \eta \leq \alpha_f^s(y^t) + t - s + 1 \leq t - s + 2$.

    \end{proof}
\begin{lemma}
        Assume \ref{as:lip}, \ref{as:strconv} and \ref{as:sim}. Then for iterates of Algorithm \ref{alg:acc} with $\theta = (p\eta^-1 - 1)/(\beta p\eta^-1 - 1), \theta > 0, \eta \geq 1$, it holds that
        \begin{align}
            \label{lemma:GD_acc_lemma5}
                \E\| &x^{t+1} - x^*\|^2 \notag\\&\leq (1-\beta)(1 + \frac{\beta}{4})\E\norm{x^{t} - x^*}^2 + \beta(1 + \frac{\beta}{4})\E\norm{x^t_g - x^*}^2 + (\beta^2 - \beta)\E\norm{x^t - x^t_g}^2\notag
                \\&+  
                10\frac{d^2}{m^2}(\delta^2+1)p^2\gamma^2\eta^2\E\norm{\nabla f(x^t_g)}^2
                +
                p^2\gamma^2\eta^2\tau\Big(32\frac{\tau^2d^2L^2p^2\gamma^2}{m^2\beta}+\frac{5}{4} \Big)\sum\limits_{r=t-\tau}^{t-1}\norm{g^r}^2\notag
                \\&+
                3\e p\gamma\eta L\frac{d}{m}\sqrt{\delta^2+1}\expect{\norm{x^{t-\tau} - x^*}^2}
                +
                3\e p\gamma\eta L\frac{d}{m}\sqrt{\delta^2+1}\expect{\norm{x_f^{t-\tau} - x^*}^2}
                \\&-
                2\gamma\eta^2\E\dotprod{\nabla f(x^t_g)}{x^t_g + (p\eta^{-1} - 1)x^t_f-p\eta^{-1}x^*} 
                + 
                2p\gamma\eta\Bigg(\frac{\e d}{m\sqrt{\delta^2+1}L} + 4p\gamma\eta\frac{d^2}{m^2}\Bigg)\sigma^2 .\notag
        \end{align}
        \begin{proof}
            Using lines 10 and 11 of Algorithm \ref{alg:acc}, we get
            \begin{align*}
                    \E\|x^{t+1} &- x^*\|^2 = \E\norm{\eta x_f^{t+1} + (p-\eta)x_f^t + (1-p)(1-\beta)x^t + (1-p)\beta x_g^t - x^*}^2  \\&=\E\norm{\eta x_g^{t} - p\gamma\eta g^t + (p-\eta)x_f^t + (1-p)(1-\beta)x^t + (1-p)\beta x_g^t - x^*}^2
                    \\&=\E\norm{\eta x_g^{t} + (p-\eta)x_f^t + (1-p)(1-\beta)x^t + (1-p)\beta x_g^t - x^*}^2 + p^2\gamma^2\eta^2\E\norm{g^t}^2 
                    \\&\quad- 2p\gamma\eta\E\dotprod{g^t}{\eta x_g^{t} + (p-\eta)x_f^t + (1-p)(1-\beta)x^t + (1-p)\beta x_g^t - x^*}
                    \\&= \underbrace{\E\norm{\eta x_g^{t} + (p-\eta)x_f^t + (1-p)(1-\beta)x^t + (1-p)\beta x_g^t - x^*}^2}_{\circledOne} + \underbrace{p^2\gamma^2\eta^2\E\norm{g^t}^2}_{\circledTwo} 
                    \\&\quad\underbrace{- 2p\gamma\eta\E\dotprod{g^t - \nabla f(x^t_g)}{\eta x_g^{t} + (p-\eta)x_f^t + (1-p)(1-\beta)x^t + (1-p)\beta x_g^t - x^*}}_{\circledThree}
                    \\&\quad\underbrace{- 2p\gamma\eta\E\dotprod{\nabla f(x^t_g)}{\eta x_g^{t} + (p-\eta)x_f^t + (1-p)(1-\beta)x^t + (1-p)\beta x_g^t - x^*}}_{\circledFour} .
            \end{align*}
            Consider \circledOne. From Lemma \ref{lemma:GD_acc_tilde}, we know that 
            \begin{equation*}
                \eta x_g^{t} + (p-\eta)x_f^t + (1-p)(1-\beta)x^t + (1-p)\beta x_g^t  = \beta x_g^t + (1-\beta)x^t.
            \end{equation*}
            It implies 
            \begin{equation}
                \label{eq:part1}
                \begin{split}
                    \|\eta x_g^{t} &+ (p-\eta)x_f^t + (1-p)(1-\beta)x^t + (1-p)\beta x_g^t - x^*\|^2
                    \\&= \norm{\beta x_g^t + (1-\beta)x^t - x^*}^2
                    \\&= \norm{\beta(x_g^t - x^t) + x^t - x^*}^2
                    \\&= \norm{x^t - x^*}^2 + 2\beta\dotprod{x^t - x^*}{x_g^t - x^t} + \beta^2\norm{x_g^t - x^t}^2
                    \\&= \norm{x^t - x^*}^2 + \beta(\norm{x^t_g - x^*}^2-\norm{x^t - x^*}^2-\norm{x_g^t - x^t}^2) + \beta^2\norm{x_g^t - x^t}^2
                    \\&= (1-\beta)\norm{x^t - x^*}^2 + \beta\norm{x_g^t - x^*}^2 + (\beta^2-\beta)\norm{x^t - x^t_g}^2.
                \end{split}
            \end{equation}
            Consider \circledTwo. Using convexity of squared Euclidean norm and Lemma \ref{lemma:main_norms}, one can obtain 
            \begin{equation}
                \label{eq:part2}
                \begin{split}
                    p^2\gamma^2\eta^2\E\norm{g^t}^2 &= p^2\gamma^2\eta^2\E\norm{\frac{1}{n}\sum\limits^{n}_{i = 1} Q_t^i(\nabla f_i(x^t_g))}^2 
                    \\&\leq p^2\gamma^2\eta^2\frac{1}{n}\sum\limits^{n}_{i = 1}\E\norm{ Q_t^i(\nabla f_i(x^t_g))}^2
                    \\&\overset{\eqref{lemma:main_norms}}{\leq} p^2\gamma^2\eta^2\frac{d^2}{m^2}\frac{1}{n}\sum\limits^{n}_{i = 1}\E\norm{ \nabla f_i(x^t_g)}^2 
                    \\&\overset{\eqref{lemma:main_sim}}{\leq} 2p^2\gamma^2\eta^2\frac{d^2}{m^2}(\delta^2 + 1)\E\norm{\nabla f (x^t_g)}^2 + 2p^2\gamma^2\eta^2\frac{d^2}{m^2}\sigma^2,
                \end{split}
            \end{equation}
            where in the last inequality we used Lemma \ref{lemma:main_sim}. 
            \\Consider \circledThree. We first use Lemma \ref{lemma:GD_acc_tilde} twice
            $$x^t_g = \theta x^t_f + (1- \theta)x^t = \alpha_f^{t - \tau}x_f^{t-\tau} + \alpha^{t-\tau}x^{t-\tau} - p\gamma\sum\limits^{t-1}_{r = t - \tau}c_rg^r$$
            \begin{equation*}
                \begin{split}
                    \eta x_g^{t} + (p-\eta)x_f^t + (1-p)(1-\beta)x^t &+ (1-p)\beta x_g^t  = \beta x_g^t + (1-\beta)x^t
                    \\&= \beta\theta x_f^t + (1-\beta\theta)x^t 
                    \\&= \hat{\alpha}_f^{t - \tau}x_f^{t-\tau} + \hat{\alpha}^{t-\tau}x^{t-\tau} - p\gamma\sum\limits^{t-1}_{r = t - \tau}\hat{c}_rg^r .
                \end{split}
            \end{equation*}
            Next, we apply Corollary \ref{cor:main_dotprod} with $\hat{a}^{t-\tau} = \nabla f_i(\widetilde{x}^{t-\tau}_g)$, where $\widetilde{x}^{t-\tau}_g = \alpha_f^{t - \tau}x_f^{t-\tau} + \alpha^{t-\tau}x^{t-\tau}$, and $\hat{b}^{t-\tau} = \hat{\alpha}_f^{t - \tau}x_f^{t-\tau} + \hat{\alpha}^{t-\tau}x^{t-\tau} - x^*$, leading us to 
            \begin{align*}
                    &-2p\gamma\eta\E\Big\langle g^t - \nabla f(x^t_g), \eta x_g^{t} + (p-\eta)x_f^t + (1-p)(1-\beta)x^t + (1-p)\beta x_g^t - x^*\Big\rangle
                    \\&= -2p\gamma\eta\frac{1}{n}\sum\limits_{i = 1}^{n}\E\Big\langle Q_t^i(\nabla f_i(x^t_g)) - \nabla f_i(x^t_g), \eta x_g^{t} + (p-\eta)x_f^t + (1-p)(1-\beta)x^t + (1-p)\beta x_g^t - x^*\Big\rangle
                    \\&\leq
                    \frac{\varepsilon d}{m \beta_0}p\gamma\eta\frac{1}{n}\sum\limits_{i = 1}^{n}\expect{\norm{\nabla f_i(\widetilde{x}^{t-\tau}_g)}^2}
                    +
                    \frac{\varepsilon d \beta_0}{m}p\gamma\eta \expect{\norm{\hat{\alpha}_f^{t - \tau}x_f^{t-\tau} + \hat{\alpha}^{t-\tau}x^{t-\tau} - x^*}^2}
                    \\&\quad+
                    4 \frac{d^2}{m^2}p\gamma\eta \left( \beta_1 + \beta_2 \right)\frac{1}{n}\sum\limits_{i = 1}^{n}\expect{\norm{\nabla f_i(x_g^t) - \nabla f_i(\widetilde{x}^{t-\tau}_g)}^2} + p\gamma\eta\left(\frac{1}{\beta_1} + \frac{1}{\beta_3} \right)
                    \expect{\norm{-p\gamma\sum\limits^{t-1}_{r = t - \tau}\hat{c}_rg^r}^2}
                    \\&\quad+
                    4 \frac{d^2}{m^2}p\gamma\eta \beta_3 \frac{1}{n}\sum\limits_{i = 1}^{n}\expect{\norm{\nabla f_i(x_g^t)}^2} 
                    +
                    \frac{p\gamma\eta}{\beta_2} \expect{\norm{\beta x_g^t + (1-\beta)x^t - x^*}^2}.
            \end{align*}
            Using Assumption \ref{as:lip} and Lemma \ref{lemma:main_sim} with $c_r \leq \tau \leq 2 \tau$ and $\hat{c}_r \leq \eta$ one might obtain 
            \begin{align}
                \label{eq:part3}
                    &- 2p\gamma\eta\E\Big\langle g^t - \nabla f(x^t_g), \eta x_g^{t} + (p-\eta)x_f^t + (1-p)(1-\beta)x^t + (1-p)\beta x_g^t - x^*\Big\rangle\notag
                    \\&\leq
                    \frac{2\varepsilon d}{m \beta_0}p\gamma\eta(\delta^2+1)\expect{\norm{\nabla f(\widetilde{x}^{t-\tau}_g)}^2}
                    +
                    \frac{\varepsilon d \beta_0}{m}p\gamma\eta \expect{\norm{\hat{\alpha}_f^{t - \tau}x_f^{t-\tau} + \hat{\alpha}^{t-\tau}x^{t-\tau} - x^*}^2}\notag
                    \\&\quad+
                    4 \frac{d^2L^2}{m^2}p\gamma\eta \left( \beta_1 + \beta_2 \right)\expect{\norm{-p\gamma\sum\limits^{t-1}_{r = t - \tau}c_rg^r}^2}
                    +
                    p\gamma\eta\left(\frac{1}{\beta_1} + \frac{1}{\beta_3} \right)
                    \expect{\norm{-p\gamma\sum\limits^{t-1}_{r = t - \tau}\hat{c}_rg^r}^2}\notag
                    \\&\quad+
                    8 \frac{d^2}{m^2}(\delta^2+1)p\gamma\eta \beta_3\expect{\norm{\nabla f(x_g^t)}^2} 
                    +
                    \frac{p\gamma\eta}{\beta_2} \expect{\norm{\beta x_g^t + (1-\beta)x^t - x^*}^2}\notag
                    \\&\quad+
                    2p\gamma\eta(\frac{\e d}{m\beta_0} + 4\frac{d^2\beta_3}{m^2})\sigma^2
                    \\&\leq
                    \frac{\varepsilon d}{m}p\gamma\eta\Bigg(2(\delta^2+1)L^2\alpha_f^{t-\tau}\frac{1}{\beta_0} + \beta_0\hat{\alpha}_f^{t - \tau}\Bigg)\expect{\norm{x_f^{t-\tau} - x^*}^2}\notag
                    \\&\quad+
                    \frac{\varepsilon d}{m}p\gamma\eta\Bigg(2(\delta^2+1)L^2\alpha^{t-\tau}\frac{1}{\beta_0} + \beta_0\hat{\alpha}^{t - \tau}\Bigg)\expect{\norm{x^{t-\tau} - x^*}^2}\notag
                    \\&\quad+
                    p^3\gamma^3\eta\tau\Bigg(4\frac{\tau^2d^2L^2}{m^2}(\beta_1+\beta_2)+\eta^2\Big(\frac{1}{\beta_1}+\frac{1}{\beta_3}\Big)\Bigg)\sum\limits_{r=t-\tau}^{t-1}\norm{g^r}^2\notag
                    \\&\quad+
                    8 \frac{d^2}{m^2}(\delta^2+1)p\gamma\eta \beta_3\expect{\norm{\nabla f(x_g^t)}^2}\notag
                    \\&\quad+
                    \frac{p\gamma\eta}{\beta_2}\beta \expect{\norm{x_g^t - x^*}^2}
                    +
                    \frac{p\gamma\eta}{\beta_2}(1 -\beta) \expect{\norm{x^t - x^*}^2}
                    +
                    2p\gamma\eta(\frac{\e d}{m\beta_0} + 4\frac{d^2\beta_3}{m^2})\sigma^2 .\notag
            \end{align}
    
            Consider \circledFour. Taking into account line 4 and the choice of $\theta$ such that $\theta = (p \eta^{-1} - 1) / (\beta p \eta^{-1} - 1)$, one can note
            \begin{align}
            \label{eq:acc_temp8}
                    &\eta x^k_g + (p - \eta)x^k_f + (1- p)(1 - \beta) x^k +(1 - p) \beta x^k_g - x^*
                    \notag\\
                    &= (\eta + (1 - p) \beta)x^k_g + (p - \eta)x^k_f + (1- p)(1 - \beta) x^k - x^*
                    \notag\\
                    &= \eta p^{-1}\left( (p + (1 - p) p^{-1}\eta \beta)x^k_g + (p \eta^{-1} - 1) p x^k_f + (1- p)(1 - \beta) p \eta^{-1} x^k - \eta^{-1} p x^* \right)
                    \notag\\
                    &= \eta p^{-1}\left( (p + (1 - p) p^{-1}\eta \beta)x^k_g + (p \eta^{-1} - 1) p x^k_f + (1- p)(1 - \beta p \eta^{-1}) (1 - \theta) x^k - \eta^{-1} p x^* \right)
                    \notag\\
                    &= \eta p^{-1}\left( (p + (1 - p) p^{-1}\eta \beta)x^k_g + (p \eta^{-1} - 1) p x^k_f + (1- p)(1 - \beta p \eta^{-1}) (x^k_g - \theta x^k_f) - \eta^{-1} p x^* \right)
                    \notag\\
                    &= \eta p^{-1}\left( x^k_g + (p \eta^{-1} - 1) p x^k_f - (1- p)(1 - \beta p \eta^{-1}) \theta x^k_f - \eta^{-1} p x^* \right)
                    \notag\\
                    &= \eta p^{-1}\left( x^k_g + (p \eta^{-1} - 1) p x^k_f - (1- p)(p \eta^{-1} - 1) x^k_f - \eta^{-1} p x^* \right)
                    \notag\\
                    &= \eta p^{-1}\left( x^k_g + (p \eta^{-1} - 1)  x^k_f - \eta^{-1} p x^* \right).
            \end{align}
            Using that, we get 
            \begin{equation}
            \label{eq:part4}
                \begin{split}
                    - 2p\gamma\eta\E\Big\langle\nabla f(x^t_g), \eta x_g^{t} + (p-\eta)x_f^t &+ (1-p)(1-\beta)x^t + (1-p)\beta x_g^t - x^*\Big\rangle 
                    \\&= -2\gamma\eta^2\E\dotprod{\nabla f(x^t_g)}{x^t_g + (p\eta^{-1} - 1)x^t_f-p\eta^{-1}x^*} .
                \end{split}
            \end{equation}
        Summing \eqref{eq:part1}, \eqref{eq:part2}, \eqref{eq:part3} and \eqref{eq:part4} with $\beta_0 = \sqrt{\delta^2 + 1}L,\; \beta_1= \beta_2 = \frac{4p\gamma\eta}{\beta}$ and $\beta_3 = p\gamma\eta$ we finish the proof.
        \end{proof}
    \end{lemma}

    \begin{lemma}
    \label{lemma:GD_acc_lemma6}
        Assume \ref{as:lip}, \ref{as:strconv} and \ref{as:sim}. Then for iterates of Algorithm \ref{alg:acc} and for any $u \in \mathbb{R}^d$ it holds that
        \begin{equation*}
        \begin{split}
            &\expect{f(x_f^{t+1})}
            \leq
            \expect{f(u)}
            -
            \expect{\dotprod{\nabla f(x_g^t)}{u - x_g^t}}
            -
            \frac{\mu}{2} \norm{u - x_g^t}
            -
            \frac{p\gamma}{2} \expect{\norm{\nabla f(x_g^t)}^2}
            \\&+
            2 \varepsilon \gamma \expect{\norm{\nabla f(\widetilde{x}_g^{t - \tau})}^2}
            +
            20 \frac{L^2 d^3 \gamma^3 p^2 \tau^3 (\delta^2 + 1)}{m^3} \sum\limits_{s = t - \tau}^{t-1} \expect{\norm{\nabla f (x_g^s)}^2}
            +
            23 \frac{L^2 d^3 \gamma^3 p^2 \tau^4}{m^3} \sigma^2,
        \end{split}
        \end{equation*}
        where 
        $$\gamma \leq \frac{1}{L} ~~\text{ and }~~ p \leq \frac{m^2}{12 (\delta^2 + 1) d^2}. $$
    \end{lemma}
    
    \begin{proof}
        Using \ref{as:lip} with $x = x_f^{t+1}$, $y=x_g^t$ and line 3 of Algorithm \ref{alg:acc} we get

        \begin{equation}
        \begin{split}
        \label{eq:tmp_acc_lemma6_1}
            \expect{f(x_f^{t+1})}
            &\leq
            \expect{f(x_g^t)}
            +
            \expect{\dotprod{\nabla f(x_g^t)}{x_f^{t+1} - x_g^t}}
            +
            \frac{L}{2} \expect{\norm{x_f^{t+1} - x_g^t}^2}
            \\&=
            \expect{f(x_g^t)}
            -
            p \gamma \expect{\dotprod{\nabla f(x_g^t)}{g^t}}
            +
            \frac{Lp^2\gamma^2}{2} \expect{\norm{g^t}^2}
            \\&=
            \expect{f(x_g^t)}
            -
            p\gamma \expect{\dotprod{\nabla f(x_g^t)}{\nabla f(x_g^t)}}
            -
            p\gamma \expect{\dotprod{\nabla f(x_g^t)}{g^k - \nabla f(x_g^t)}}
            \\&\quad+
            \frac{Lp^2\gamma^2}{2} \expect{\norm{g^t}^2} .
        \end{split}
        \end{equation}

        Consider $\expect{\dotprod{\nabla f(x_g^t)}{g^k - \nabla f(x_g^t)}}$. Using Corollary \ref{cor:main_dotprod} with $a^t = \nabla f_i(x^t_g), b^t = \nabla f(x^t_g), \hat{a}^{t-\tau} = \nabla f_i(\widetilde{x}_g^{t-\tau}), \hat{b}^{t-\tau} = \nabla f(\widetilde{x}_g^{t-\tau})$, where $x_g^t \in \text{conv}\left\{ x_f^t, x^t\right\} = \widetilde{x}_g^{t-\tau} -p \gamma \sum_{s = t - \tau}^{t-1} c_s g^s$ from Lemma \ref{lemma:GD_acc_tilde}. Using Assumption \ref{as:lip} we obtain 
        \begin{align*}
            &2 \left|\expect{\dotprod{\nabla f(x_g^t)}{g^k - \nabla f(x_g^t)}} \right|
            \leq
            \frac{\varepsilon d}{m \beta_0} \expect{\frac{1}{n} \sum\limits_{i=1}^n
            \norm{\nabla f_i(\widetilde{x}_g^{t - \tau})}^2}
            +
            \frac{\varepsilon d \beta_0}{m} \expect{\norm{\nabla f(\widetilde{x}_g^{t - \tau})}^2 }
            \\&+
            4 \frac{d^2 L^2}{m^2}(\beta_1 + \beta_2) \expect{\norm{x_g^t - \widetilde{x}_g^{t-\tau}}^2}
            +
            L^2 \left(\frac{1}{\beta_1} + \frac{1}{\beta_3}\right) \expect{\norm{x_g^t - \widetilde{x}_g^{t-\tau}}^2}
            \\&+
            4 \frac{d^2}{m^2} \beta_3 \expect{\frac{1}{n} \sum\limits_{i=1}^n
            \norm{\nabla f_i(x_g^t)}^2}
            +
            \frac{1}{\beta_2} \expect{\norm{\nabla f(x_g^t)}^2} .
        \end{align*}

        Taking $\beta_0 = \sqrt{\delta^2 + 1}$, $\beta_1 = m/d$, $\beta_2 = m/(d p)$, $\beta_3 = p m/d$ and using results from Lemma \ref{lemma:main_sim} we obtain 

        \begin{align*}
            &2 \left|\expect{\dotprod{\nabla f(x_g^t)}{g^k - \nabla f(x_g^t)}} \right|
            \leq
            \frac{2 \varepsilon d}{m} \left(\sqrt{\delta^2 +1} \expect{\norm{\nabla f(\widetilde{x}_g^{t - \tau})}^2} + \frac{\sigma^2}{\sqrt{\delta^2 +1}} \right)
            \\&+
            \frac{d p}{m} \expect{\norm{\nabla f(x_g^t)}^2}
            +
            10 \frac{L^2 d}{p m} \expect{\norm{-p \gamma \sum\limits_{s = t - \tau}^{t-1} c_s \frac{1}{n} \sum\limits_{i=1}^n Q_s^i(\nabla f_i(x_g^s))}^2}
            \\&+
            \frac{8 d p}{m} \left( (\delta^2 + 1) \expect{\norm{\nabla f(x_g^t)}^2} + \sigma^2\right)
            +
            \frac{\varepsilon d \sqrt{\delta^2 +1}}{m} \expect{\norm{\nabla f(\widetilde{x}_g^{t - \tau})}^2} .
        \end{align*}

        Using Lemma \ref{lemma:main_norms} and \ref{lemma:main_sim}, convexity of the squared norm and the fact that $c_s \leq t - s + 2 \leq \tau +2 \leq 2\tau$ we obtain
        \begin{equation*}
        \begin{split}
            2 \left|\expect{\dotprod{\nabla f(x_g^t)}{g^k - \nabla f(x_g^t)}} \right|
            &\leq
            \frac{3 \varepsilon d \sqrt{\delta^2 +1}}{m}\expect{\norm{\nabla f(\widetilde{x}_g^{t - \tau})}^2}
            +
            \\&+
            40 \frac{L^2 d^3 \gamma^2 p \tau^3}{m^3} \sum\limits_{s = t - \tau}^{t-1} \expect{(\delta^2 +1 ) \norm{\nabla f (x_g^s)}^2 + \sigma^2}
            \\&+
            \frac{9 d p (\delta^2 + 1) }{m} \expect{\norm{\nabla f(x_g^t)}^2}
            +
            \frac{2 d}{m} \left( \frac{\varepsilon}{\sqrt{\delta^2 +1}} + p \right) \sigma^2 .
        \end{split}
        \end{equation*}

        Using the fact that $L^2 \gamma^2 d^2/m^2 \tau^4 \eta^2 \geq 1$ and $\varepsilon \leq \sqrt{\delta^2 + 1} p$ we obtain

        \begin{align*}
            &2 \left|\expect{\dotprod{\nabla f(x_g^t)}{g^k - \nabla f(x_g^t)}} \right|
            \leq
            \frac{3 \varepsilon d \sqrt{\delta^2 +1}}{m}\expect{\norm{\nabla f(\widetilde{x}_g^{t - \tau})}^2}
            +
            44 \frac{L^2 d^3 \gamma^2 p \eta^2 \tau^4}{m^3} \sigma^2
            \\&+
            40 \frac{L^2 d^3 \gamma^2 p \tau^3 (\delta^2 + 1)}{m^3} \sum\limits_{s = t - \tau}^{t-1} \expect{\norm{\nabla f (x_g^s)}^2}
            +
            \frac{9 d p (\delta^2 + 1) }{m} \expect{\norm{\nabla f(x_g^t)}^2} .
        \end{align*}

        Using this result, Lemmas \ref{lemma:main_norms} and \ref{lemma:main_sim} we can estimate \eqref{eq:tmp_acc_lemma6_1}:

        \begin{align*}
            \expect{f(x_f^{t+1})}
            &=
            \expect{f(x_g^t)}
            -
            p\gamma \expect{\norm{\nabla f(x_g^t)}^2}
            \\&-
            p\gamma \expect{\dotprod{\nabla f(x_g^t)}{g^k - \nabla f(x_g^t)}}
            +
            \frac{L}{2} \expect{\norm{g^t}^2} 
            \\&\leq
            \expect{f(x_g^t)}
            -
            p\gamma \expect{\norm{\nabla f(x_g^t)}^2}
            +
            \frac{2 \varepsilon p \gamma d \sqrt{\delta^2 +1}}{m}\expect{\norm{\nabla f(\widetilde{x}_g^{t - \tau})}^2}
            +
            \\&+
            20 \frac{L^2 d^3 \gamma^3 p^2 \tau^3 (\delta^2 + 1)}{m^3} \sum\limits_{s = t - \tau}^{t-1} \expect{\norm{\nabla f (x_g^s)}^2}
            +
            \frac{5 d \gamma p^2 (\delta^2 + 1) }{m} \expect{\norm{\nabla f(x_g^t)}^2}
            \\&+
            22 \frac{L^2 d^3 \gamma^3 p^2 \tau^4}{m^3} \sigma^2 
            +
            \frac{L p^2 \gamma^2 d^2}{m^2} (\delta^2 + 1) \expect{\norm{\nabla f(x_g^t)}^2} +  \frac{L p^2 \gamma^2 d^2}{m^2} \sigma^2 .
        \end{align*}

        Taking 
        
        $$\gamma \leq \frac{1}{L} ~~\text{ and }~~ p \leq \frac{m^2}{12 (\delta^2 + 1) d^2}, $$

        we obtain

        \begin{align*}
            \expect{f(x_f^{t+1})}
            &\leq
            \expect{f(x_g^t)}
            -
            \frac{p\gamma}{2} \expect{\norm{\nabla f(x_g^t)}^2}
            +
            2 \varepsilon \gamma \expect{\norm{\nabla f(\widetilde{x}_g^{t - \tau})}^2}
            +
            \\&+
            20 \frac{L^2 d^3 \gamma^3 p^2 \tau^3 (\delta^2 + 1)}{m^3} \sum\limits_{s = t - \tau}^{t-1} \expect{\norm{\nabla f (x_g^s)}^2}
            +
            23 \frac{L^2 d^3 \gamma^3 p^2 \tau^4}{m^3} \sigma^2 .
        \end{align*}

        Using \ref{as:strconv} with $x = u$ and $y = x_g^t$, one can conclude that for any $u \in \mathbb{R}^d$ it holds

        \begin{equation*}
        \begin{split}
            \expect{f(x_f^{t+1})}
            &\leq
            \expect{f(u)}
            -
            \expect{\dotprod{\nabla f(x_g^t)}{u - x_g^t}}
            -
            \frac{\mu}{2} \norm{u - x_g^t}
            \\&\quad-
            \frac{p\gamma}{2} \expect{\norm{\nabla f(x_g^t)}^2}
            +
            2 \varepsilon \gamma \expect{\norm{\nabla f(\widetilde{x}_g^{t - \tau})}^2}
            +
            \\&\quad+
            20 \frac{L^2 d^3 \gamma^3 p^2 \tau^3 (\delta^2 + 1)}{m^3} \sum\limits_{s = t - \tau}^{t-1} \expect{\norm{\nabla f (x_g^s)}^2}
            +
            23 \frac{L^2 d^3 \gamma^3 p^2 \tau^4}{m^3} \sigma^2 .
        \end{split}
        \end{equation*}

        This finishes the proof.
    \end{proof}
    \begin{theorem}[Theorem \ref{theorem:GD_acc}]
        Consider Assumptions \ref{as:lip}, \ref{as:strconv} and \ref{as:sim}. Let problem \eqref{eq:problem} be solved by Algorithm \ref{alg:acc}. Then for  any $\gamma > 0, \varepsilon > 0$, $\tau > \tau_{\text{mix}}(\varepsilon)$, $T > \tau$ and $\beta, \theta, \eta, p$ satisfying 
        $$
        \gamma\leq\frac{\mu^{\frac{1}{3}}m^{\frac{1}{2}}}{2\tau L^{\frac{4}{3}}d^{\frac{1}{2}}},
        ~~~~
        \e\leq\min\Big\{\frac{m^{\frac{7}{4}}}{6d^{\frac{7}{4}}\tau^{\frac{5}{4}}L(\delta^2 + 1)}; \frac{m^{\frac{5}{4}}}{\sqrt{2}\tau^{\frac{3}{4}}\mu^{\frac{1}{3}}L^{\frac{2}{3}}d^{\frac{5}{4}}}; \frac{m^{\frac{15}{4}}}{6d^{\frac{15}{4}}\tau^{\frac{13}{4}}(\delta^2 + 1)^2}\Big\}, 
        $$
        $$
        p\leq\frac{m^2}{13d^2(\delta^2+1)\tau^2},
        ~~~~
        \beta = \sqrt{\frac{2p^2\mu\gamma}{3}},
        ~~~~
        \eta = \sqrt{\frac{3}{2\mu\gamma}},
        ~~~~
        \theta = \frac{p\eta^{-1} - 1}{\beta p \eta^{-1} - 1} .
        $$
        it holds that
        \begin{align*}
                \E[\|x^{T+1} - x^*\|^2 + \frac{3}{\mu}( f(x_f^{T+1}) - f(x^*))] 
                &\leq 
                \exp\Bigg(-(T - \tau)\sqrt{\frac{2p^2\mu\gamma}{3}} \Bigg)F_{\tau}
                \\&+ 
                \exp\Bigg(-T\sqrt{\frac{2p^2\mu\gamma}{3}} \Bigg)\Delta_{\tau} + \frac{45\gamma}{\mu}\sigma^2, 
        \end{align*}
        where $F_\tau := \E[\|x^{\tau} - x^*\|^2 + \frac{3}{\mu}( f(x_f^{\tau}) - f(x^*))]$ and $\Delta_{\tau} \leq \frac{\sqrt{\gamma}}{\tau^{\frac{4}{3}}\mu^{\frac{1}{3}}}\sum\limits_{t=0}^{\tau}\Big(\E\norm{\nabla f(x_g^t)} + \E\norm{x^{t} - x^*}^2 + \E[f(x_f^{t}) - f(x^*)] \Big)$.
    \end{theorem}
    \begin{proof}
        We start by using Lemma \ref{lemma:GD_acc_lemma6} with $u = x^*$ and $u = x_f^t$
        \begin{align*}
            \expect{f(x_f^{t+1})}
            &\leq
            \expect{f(x^*)}
            -
            \expect{\dotprod{\nabla f(x_g^t)}{x^* - x_g^t}}
            -
            \frac{\mu}{2} \norm{x^* - x_g^t}
            -
            \frac{p\gamma}{2} \expect{\norm{\nabla f(x_g^t)}^2}
            \\&\hspace{-18mm}+
            2 \varepsilon \gamma \expect{\norm{\nabla f(\widetilde{x}_g^{t - \tau})}^2}
            +
            20 \frac{L^2 d^3 \gamma^3 p^2 \tau^3 (\delta^2 + 1)}{m^3} \sum\limits_{s = t - \tau}^{t-1} \expect{\norm{\nabla f (x_g^s)}^2}
            +
            23 \frac{L^2 d^3 \gamma^3 p^2 \tau^4}{m^3} \sigma^2,
            \\
            \expect{f(x_f^{t+1})}
            &\leq
            \expect{f(x^t_f)}
            -
            \expect{\dotprod{\nabla f(x_g^t)}{x^t_f - x_g^t}}
            -
            \frac{\mu}{2} \norm{x^t_f - x_g^t}
            -
            \frac{p\gamma}{2} \expect{\norm{\nabla f(x_g^t)}^2}
            \\&\hspace{-18mm}+
            2 \varepsilon \gamma \expect{\norm{\nabla f(\widetilde{x}_g^{t - \tau})}^2}
            +
            20 \frac{L^2 d^3 \gamma^3 p^2 \tau^3 (\delta^2 + 1)}{m^3} \sum\limits_{s = t - \tau}^{t-1} \expect{\norm{\nabla f (x_g^s)}^2}
            +
            23 \frac{L^2 d^3 \gamma^3 p^2 \tau^4}{m^3} \sigma^2 .
        \end{align*}
        Summing the first inequality with coefficient $2p\gamma\eta$, the second with coefficient $2p\gamma\eta(\eta-p)$ and \eqref{lemma:GD_acc_lemma5}, we get 
        \begin{align*}
                \E[\|x^{t+1} - x^*\|^2 &+ 2\gamma\eta^2f(x_f^{t+1})]
                \\&\leq 
                (1-\beta)(1 + \frac{\beta}{4})\E\norm{x^{t} - x^*}^2 
                +
                \beta(1 + \frac{\beta}{4})\E\norm{x^t_g - x^*}^2 
                +
                (\beta^2 - \beta)\E\norm{x^t - x^t_g}^2
                \\&\quad+
                10\frac{d^2}{m^2}(\delta^2+1)p^2\gamma^2\eta^2\E\norm{\nabla f(x^t_g)}^2
                +
                p^2\gamma^2\eta^2\tau\Big(32\frac{\tau^2d^2L^2p^2\gamma^2}{m^2\beta}+\frac{5}{4} \Big)\sum\limits_{r=t-\tau}^{t-1}\norm{g^r}^2
                \\&\quad+
                3\e p\gamma\eta L\frac{d}{m}\sqrt{\delta^2+1}\expect{\norm{x^{t-\tau} - x^*}^2}
                +
                3\e p\gamma\eta L\frac{d}{m}\sqrt{\delta^2+1}\expect{\norm{x_f^{t-\tau} - x^*}^2}
                \\&\quad-
                2\gamma\eta^2\E\dotprod{\nabla f(x^t_g)}{x^t_g + (p\eta^{-1} - 1)x^t_f-p\eta^{-1}x^*} 
                + 
                2p\gamma\eta\Bigg(\frac{\e d}{m\sqrt{\delta^2+1}L} + 4p\gamma\eta\frac{d^2}{m^2}\Bigg)\sigma^2 
                \\&\quad+
                2p\gamma\eta\Bigg(\expect{f(x^*)}
                -
                \expect{\dotprod{\nabla f(x_g^t)}{x^* - x_g^t}}
                -
                \frac{\mu}{2} \norm{x^* - x_g^t}
                -
                \frac{p\gamma}{2} \expect{\norm{\nabla f(x_g^t)}^2}
                \\&\quad+
                2 \varepsilon \gamma \expect{\norm{\nabla f(\widetilde{x}_g^{t - \tau})}^2}
                +
                20 \frac{L^2 d^3 \gamma^3 p^2 \tau^3 (\delta^2 + 1)}{m^3}\sum\limits_{s = t - \tau}^{t-1} \expect{\norm{\nabla f (x_g^s)}^2}
                \\&\quad+
                23 \frac{L^2 d^3 \gamma^3 p^2 \tau^4}{m^3} \sigma^2\Bigg)
                \\&\quad+
                2\gamma\eta(\eta - p)\Bigg(\expect{f(x_f^t)}
                -
                \expect{\dotprod{\nabla f(x_g^t)}{x_f^t - x_g^t}}
                -
                \frac{\mu}{2} \norm{x_f^t - x_g^t}
                -
                \frac{p\gamma}{2} \expect{\norm{\nabla f(x_g^t)}^2}
                \\&\quad+
                2 \varepsilon \gamma \expect{\norm{\nabla f(\widetilde{x}_g^{t - \tau})}^2}
                +
                20 \frac{L^2 d^3 \gamma^3 p^2 \tau^3 (\delta^2 + 1)}{m^3} \sum\limits_{s = t - \tau}^{t-1} \expect{\norm{\nabla f (x_g^s)}^2}
                \\&\quad+
                23 \frac{L^2 d^3 \gamma^3 p^2 \tau^4}{m^3} \sigma^2\Bigg)
                \\&\leq
                (1-\beta)(1 + \frac{\beta}{4})\E\norm{x^{t} - x^*}^2 
                +
                (\beta + \frac{\beta^2}{4} - p\gamma\eta\mu)\E\norm{x^t_g - x^*}^2 
                +
                (\beta^2 - \beta)\E\norm{x^t - x^t_g}^2
                \\&\quad+
                p^2\gamma^2\eta^2\Bigg(10\frac{d^2}{m^2}(\delta^2+1) - \frac{1}{p}\Bigg)\E\norm{\nabla f(x_g^t)} 
                +
                2p\gamma\eta\E f(x^*)
                +
                2\gamma\eta(\eta-p)\E f(x_f^t)
                \\&\quad+
                p^2\gamma^2\eta^2\tau(\delta^2 + 1)\frac{d^2}{m^2}\Big(32\frac{\tau^2d^2L^2p^2\gamma^2}{m^2\beta}+\frac{5}{4} \Big)\sum\limits_{r = t - \tau}^{t-1}\E\norm{\nabla f(x_g^r)}
                \\&\quad+
                \e\gamma\eta L(3p\frac{d}{m}\sqrt{\delta^2+1}+2\gamma\eta L)\expect{\norm{x^{t-\tau} - x^*}^2}
                \\&\quad+
                \e\gamma\eta L(3p\frac{d}{m}\sqrt{\delta^2+1}+2\gamma\eta L)\expect{\norm{x_f^{t-\tau} - x^*}^2}
                \\&\quad+
                2p\gamma\eta\Bigg(\frac{\e d}{m\sqrt{\delta^2+1}L} + 4p\gamma\eta\frac{d^2}{m^2} 
                \\&\qquad\qquad\quad+ 23p\gamma^3\eta\tau^4\frac{d^3}{m^3}L^2 + p\gamma\eta\tau^2\frac{d^2}{m^2}\Bigg(16\frac{\tau^2d^2L^2p^2\gamma^2}{m^2\beta}+\frac{5}{8} \Bigg)\Bigg)\sigma^2,
        \end{align*}
        where in the last inequality we used Lemma \ref{lemma:main_sim} and Assumption \ref{as:lip}. Since $\beta < 1$, the choice of $p\gamma\eta\mu = \frac{3\beta}{2}$ gives 
        \begin{align*}
            &(1-\beta)(1+\frac{\beta}{4}) \leq 1 - \frac{3\beta}{4},\\
            &\beta+\frac{\beta^2}{4} - p\gamma\eta\mu \leq \frac{3\beta}{2} - p\gamma\eta\mu \leq 0,\\
            &\beta^2-\beta \leq0.
        \end{align*}
        This lead us to 
        \begin{align}
            \label{eq:gd_acc_jopa}
                &\E[\|x^{t+1} - x^*\|^2 + 2\gamma\eta^2( f(x_f^{t+1}) - f(x^*))] \leq (1 - \frac{3\beta}{4})\E\norm{x^{t} - x^*}^2 
                +
                2p\gamma\eta^2(1-\frac{p}{\eta})\E[f(x_f^t) - f(x^*)]\notag
                \\&\quad+
                p^2\gamma^2\eta^2\Bigg(10\frac{d^2}{m^2}(\delta^2+1) - \frac{1}{p}\Bigg)\E\norm{\nabla f(x_g^t)}
                +
                p^2\gamma^2\eta^2\tau(\delta^2 + 1)\frac{d^2}{m^2}\Big(32\frac{\tau^2d^2L^2p^2\gamma^2}{m^2\beta}+\frac{5}{4} \Big)\sum\limits_{r = t - \tau}^{t-1}\E\norm{\nabla f(x_g^r)}
                \\&\quad+
                \e\gamma\eta L(3p\frac{d}{m}\sqrt{\delta^2+1}+2\gamma\eta L)\expect{\norm{x^{t-\tau} - x^*}^2}
                +
                \e\gamma\eta L(3p\frac{d}{m}\sqrt{\delta^2+1}+2\gamma\eta L)\frac{2}{\mu}\E[f(x_f^{t-\tau}) - f(x^*)]\notag
                \\&\quad+
                2p\gamma\eta\Bigg(\frac{\e d}{m\sqrt{\delta^2+1}L} + 4p\gamma\eta\frac{d^2}{m^2} + 23p\gamma^3\eta\tau^4\frac{d^3}{m^3}L^2 + p\gamma\eta\tau^2\frac{d^2}{m^2}\Bigg(16\frac{\tau^2d^2L^2p^2\gamma^2}{m^2\beta}+\frac{5}{8} \Bigg)\Bigg)\sigma^2,\notag
        \end{align}
        where we also used Assumption \ref{as:strconv} and subtracted $2\gamma\eta^2f(x^*)$ from both sides. Next, we perform the summation from $t = \tau$ to $t = T > \tau$ of equations \eqref{eq:gd_acc_jopa} with coefficients $p_t$:
        \begin{align*}
                &\sum\limits_{t=\tau}^{T}p_t\E[\|x^{t+1} - x^*\|^2 + 2\gamma\eta^2(f(x_f^{t+1}) - f(x^*))]
                \\&\leq
                \sum\limits_{t=\tau}^{T}p_t(1 - \frac{3\beta}{4})\E\norm{x^{t} - x^*}^2 
                \\&+
                \sum\limits_{t=\tau}^{T}p_t2p\gamma\eta^2(1-\frac{p}{\eta})\E[f(x_f^t) - f(x^*)]
                +
                \sum\limits_{t=\tau}^{T}p_tp^2\gamma^2\eta^2\Bigg(10\frac{d^2}{m^2}(\delta^2+1) - \frac{1}{p}\Bigg)\E\norm{\nabla f(x_g^t)} 
                \\&+
                \sum\limits_{t=\tau}^{T}p_tp^2\gamma^2\eta^2\tau(\delta^2 + 1)\frac{d^2}{m^2}\Big(32\frac{\tau^2d^2L^2p^2\gamma^2}{m^2\beta}+\frac{5}{4} \Big)\sum\limits_{r = t - \tau}^{t-1}\E\norm{\nabla f(x_g^r)}
                \\&+
                \sum\limits_{t=\tau}^{T}p_t\e\gamma\eta L(3p\frac{d}{m}\sqrt{\delta^2+1}+2\gamma\eta L)\expect{\norm{x^{t-\tau} - x^*}^2}
                \\&+
                \sum\limits_{t=\tau}^{T}p_t\e\gamma\eta L(3p\frac{d}{m}\sqrt{\delta^2+1}+2\gamma\eta L)\frac{2}{\mu}\E[f(x_f^{t-\tau}) - f(x^*)]
                \\&+
                \sum\limits_{t=\tau}^{T}p_t2p\gamma\eta\Bigg(\frac{\e d}{m\sqrt{\delta^2+1}L} + 4p\gamma\eta\frac{d^2}{m^2} 
                + 23p\gamma^3\eta\tau^4\frac{d^3}{m^3}L^2 + p\gamma\eta\tau^2\frac{d^2}{m^2}\Bigg(16\frac{\tau^2d^2L^2p^2\gamma^2}{m^2\beta}+\frac{5}{8} \Bigg)\Bigg)\sigma^2.
        \end{align*}
        Similar as in Theorem \ref{theorem:GD_sem} we take $p_t = p^t$, $p = (1-\frac{\beta}{2})^{-1}$, it implies $p_{\tau} \leq 6$ and therefore
        \begin{align*}
                &\sum\limits_{t=\tau}^{T}p_t\E[\|x^{t+1} - x^*\|^2 + 2\gamma\eta^2( f(x_f^{t+1}) - f(x^*))]
                \\&\leq
                \sum\limits_{t=\tau}^{T}p_t\Bigg(1 - \frac{3\beta}{4} + 6\e\gamma\eta L\Big(3p\frac{d}{m}\sqrt{\delta^2+1}+2\gamma\eta L\Big)\Bigg)\E\norm{x^{t} - x^*}^2 
                \\&+
                \sum\limits_{t=\tau}^{T}p_t\Bigg(2p\gamma\eta^2(1-\frac{p}{\eta}) + 12\frac{\e\gamma\eta L}{\mu}\Big(3p\frac{d}{m}\sqrt{\delta^2+1}+2\gamma\eta L\Big) \Bigg)\E[f(x_f^t) - f(x^*)]
                \\&+
                \sum\limits_{t=\tau}^{T}p_tp^2\gamma^2\eta^2\Bigg(10\frac{d^2}{m^2}(\delta^2+1) - \frac{1}{p} + \tau^2(\delta^2 + 1)\frac{d^2}{m^2}\Bigg(32\frac{\tau^2d^2L^2p^2\gamma^2}{m^2\beta}+\frac{5}{4} \Bigg)\Bigg)\E\norm{\nabla f(x_g^t)} 
                \\&+
                \sum\limits_{t=0}^{\tau}p_{t+\tau}8p^2\gamma^4\eta^2(\delta^2+1)\frac{d^3}{m^3}\tau^3L^2\Bigg(\frac{2p^2d}{m\beta} + 5\Bigg)\sum\limits_{r = t - \tau}^{t-1}\E\norm{\nabla f(x_g^r)}
                \\&+
                \sum\limits_{t=0}^{\tau}p_{t+\tau}\e\gamma\eta L(3p\frac{d}{m}\sqrt{\delta^2+1}+2\gamma\eta L)\expect{\norm{x^{t} - x^*}^2}
                \\&+
                \sum\limits_{t=0}^{\tau}p_{t+\tau}\e\gamma\eta L(3p\frac{d}{m}\sqrt{\delta^2+1}+2\gamma\eta L)\frac{2}{\mu}\E[f(x_f^{t}) - f(x^*)]
                \\&+
                \sum\limits_{t=\tau}^{T}p_t2p\gamma\eta\Bigg(\frac{\e d}{m\sqrt{\delta^2+1}L} + 4p\gamma\eta\frac{d^2}{m^2} + 23p\gamma^3\eta\tau^4\frac{d^3}{m^3}L^2 
                + p\gamma\eta\tau^2\frac{d^2}{m^2}\Bigg(16\frac{\tau^2d^2L^2p^2\gamma^2}{m^2\beta}+\frac{5}{8} \Bigg)\Bigg)\sigma^2.
        \end{align*}
        Taking 
        \begin{align*}
            \gamma\leq\frac{\mu^{\frac{1}{3}}m^{\frac{1}{2}}}{2\tau L^{\frac{4}{3}}d^{\frac{1}{2}}}
            ~,~~~~
            &p\leq\frac{m^2}{13d^2(\delta^2+1)\tau^2},
            \\
            \e\leq\min\Big\{\frac{m^{\frac{7}{4}}}{6d^{\frac{7}{4}}\tau^{\frac{5}{4}}L(\delta^2 + 1)}; &\frac{m^{\frac{5}{4}}}{\sqrt{2}\tau^{\frac{3}{4}}\mu^{\frac{1}{3}}L^{\frac{2}{3}}d^{\frac{5}{4}}}; \frac{m^{\frac{15}{4}}}{6d^{\frac{15}{4}}\tau^{\frac{13}{4}}(\delta^2 + 1)^2}\Big\}, 
        \end{align*}
        we get 
        \begin{align*}
                &10\frac{d^2}{m^2}(\delta^2+1) - \frac{1}{p} + \tau^2(\delta^2 + 1)\frac{d^2}{m^2}\Bigg(32\frac{\tau^2d^2L^2p^2\gamma^2}{m^2\beta}+\frac{5}{4} \Bigg)\leq0,\\&
                6\e\gamma\eta L\Big(3p\frac{d}{m}\sqrt{\delta^2+1}+2\gamma\eta L\Big)\leq\frac{\beta}{4},\\&
                12\frac{\e\gamma\eta L}{\mu}(3p\frac{d}{m}\sqrt{\delta^2+1}+2\gamma\eta L)\leq 2p\gamma\eta^2\frac{p}{2\eta}, 
        \end{align*}
        \vspace{-0.1cm}and therefore with $\beta = \frac{p}{\eta}$\newline
        \begin{align*}
                &\sum\limits_{t=\tau}^{T}p_t\E[\|x^{t+1} - x^*\|^2 + 2\gamma\eta^2( f(x_f^{t+1}) - f(x^*))]
                \leq
                \sum\limits_{t=\tau}^{T}p_t\Big(1 - \frac{\beta}{2}\Big)\E[\|x^t - x^*\|^2 + 2\gamma\eta^2( f(x_f^t) - f(x^*))] 
                \\&\quad+
                \sum\limits_{t=0}^{\tau}p_{t+\tau}8p^2\gamma^4\eta^2(\delta^2+1)\frac{d^3}{m^3}\tau^3L^2\Bigg(\frac{2p^2d}{m\beta} + 5\Bigg)\sum\limits_{r = t - \tau}^{t-1}\E\norm{\nabla f(x_g^r)}
                \\&\quad+
                \sum\limits_{t=0}^{\tau}p_{t+\tau}\e\gamma\eta L(3p\frac{d}{m}\sqrt{\delta^2+1}+2\gamma\eta L)\expect{\norm{x^{t} - x^*}^2}
                \\&\quad+
                \sum\limits_{t=0}^{\tau}p_{t+\tau}\e\gamma\eta L(3p\frac{d}{m}\sqrt{\delta^2+1}+2\gamma\eta L)\frac{2}{\mu}\E[f(x_f^{t}) - f(x^*)]
                \\&\quad+
                \sum\limits_{t=\tau}^{T}p_t2p\gamma\eta\Bigg(\frac{\e d}{m\sqrt{\delta^2+1}L} + 4p\gamma\eta\frac{d^2}{m^2} 
                + 23p\gamma^3\eta\tau^4\frac{d^3}{m^3}L^2 + p\gamma\eta\tau\frac{d^2}{m^2}\Bigg(16\frac{\tau^2d^2L^2p^2\gamma^2}{m^2\beta}+\frac{5}{8} \Bigg)\Bigg)\sigma^2.
        \end{align*}
        Assume the following notation
        \begin{align*}
                \Delta_{\tau} &:= \sum\limits_{t=0}^{\tau}p_{t+\tau}8p^2\gamma^4\eta^2(\delta^2+1)\frac{d^3}{m^3}\tau^3L^2\Bigg(\frac{2p^2d}{m\beta} + 5\Bigg)\sum\limits_{r = t - \tau}^{t-1}\E\norm{\nabla f(x_g^r)}
                \\&\quad+
                \sum\limits_{t=0}^{\tau}p_{t+\tau}\e\gamma\eta L(3p\frac{d}{m}\sqrt{\delta^2+1}+2\gamma\eta L)\expect{\norm{x^{t} - x^*}^2}
                \\&\quad+
                \sum\limits_{t=0}^{\tau}p_{t+\tau}\e\gamma\eta L(3p\frac{d}{m}\sqrt{\delta^2+1}+2\gamma\eta L)\frac{2}{\mu}\E[f(x_f^{t}) - f(x^*)]
                \\&\leq \frac{\sqrt{\gamma}}{\tau^{\frac{4}{3}}\mu^{\frac{1}{3}}}\sum\limits_{t=0}^{\tau}\Bigg(\E\norm{\nabla f(x_g^t)} + \E\norm{x^{t} - x^*}^2 + \E[f(x_f^{t}) - f(x^*)] \Bigg)
        \end{align*}
        Now we substitute $p_t$, this lead us to 
        \begin{equation*}
            \begin{split}
                \sum\limits_{t=\tau}^{T}&\Big(1 - \frac{\beta}{2}\Big)^{-t}\E[\|x^{t+1} - x^*\|^2 + 2\gamma\eta^2( f(x_f^{t+1}) - f(x^*))] 
                \leq \sum\limits_{t=\tau}^{T}\Big(1 - \frac{\beta}{2}\Big)^{-t+1}\E[\|x^t - x^*\|^2 + 2\gamma\eta^2( f(x_f^t) - f(x^*))]
                \\&+ \Delta_{\tau} +
                \sum\limits_{t=\tau}^{T}\Big(1 - \frac{\beta}{2}\Big)^{-t}2p\gamma\eta\Bigg(\frac{\e d}{m\sqrt{\delta^2+1}L} + 4p\gamma\eta\frac{d^2}{m^2} 
                + 23p\gamma^3\eta\tau^4\frac{d^3}{m^3}L^2 + p\gamma\eta\tau\frac{d^2}{m^2}\Bigg(16\frac{\tau^2d^2L^2p^2\gamma^2}{m^2\beta}+\frac{5}{8} \Bigg)\Bigg)\sigma^2.
            \end{split}
        \end{equation*}
        This implies
        \begin{equation*}
            \begin{split}
                \Big(1 - \frac{\beta}{2}&\Big)^{-T}\E[\|x^{T+1} - x^*\|^2 + 2\gamma\eta^2( f(x_f^{T+1}) - f(x^*))] 
                \leq \Big(1 - \frac{\beta}{2}\Big)^{\tau}\E[\|x^{\tau} - x^*\|^2 
                + 2\gamma\eta^2( f(x_f^{\tau}) - f(x^*))]
                +
                \Delta_{\tau}
                \\&+
                \sum\limits_{t=\tau}^{T}\Big(1 - \frac{\beta}{2}\Big)^{-t}2p\gamma\eta\Bigg(\frac{\e d}{m\sqrt{\delta^2+1}L} + 4p\gamma\eta\frac{d^2}{m^2} 
                + 23p\gamma^3\eta\tau^4\frac{d^3}{m^3}L^2 + p\gamma\eta\tau\frac{d^2}{m^2}\Bigg(16\frac{\tau^2d^2L^2p^2\gamma^2}{m^2\beta}+\frac{5}{8} \Bigg)\Bigg)\sigma^2.
            \end{split}
        \end{equation*}
        Rearranging this inequality and taking $\e \leq 
 \frac{\sqrt{\gamma}m}{\sqrt{\mu}d}$ we obtain
        \begin{equation*}
            \begin{split}
                \E[\|x^{T+1} - x^*\|^2 + 2\gamma\eta^2( f(x_f^{T+1}) - f(x^*))] 
                &\leq \Big(1 - \frac{\beta}{2}\Big)^{T -\tau}\E[\|x^{\tau} - x^*\|^2 + 2\gamma\eta^2( f(x_f^{\tau}) - f(x^*))]
                \\&\quad+
                \Big(1 - \frac{\beta}{2}\Big)^{T}\Delta_{\tau}
                +
                6\sqrt{\frac{\gamma}{\mu}} \sigma^2.
            \end{split}
        \end{equation*}
        This finishes the proof.
    \end{proof}

\newpage
\section{Experiments}
\label{appx:experiments}
This section provides description of the experiment setup, presents and analyses results of logistic regression experiments on LIBSVM datasets, studies dependence of history size over convergence. Moreover, experiments with neural networks optimization for data-parallelism and model-parallelism are presented and discussed.

\subsection{Technical details}
\label{appx:experiments_setup}

Our implementation of compression operators and algorithms is written in Python 3.10, with the use of PyTorch optimization library. We implement a simulation of distributed optimization system on a single machine, which is equivalent in terms of convergence analysis. Our server is AMD Ryzen Threadripper 2950X 16-Core Processor @ 2.2 GHz CPU and x2 NVIDIA GeForce GTX 1080 Ti GPU. We use Weights\&Biases~\cite{wandb} for experiments tracking and hyperparameters tuning.

\subsection{Logistic Regression experiments}
\label{appendix_subsec:exp_libsvm}
We conduct experiments on classification with logistic regression on four datasets: \texttt{Mushrooms}, \texttt{A9A}, \texttt{W8A}, \texttt{MNIST}. We apply the following optimization algorithms: proposed \texttt{MQSGD} and its accelerated version \texttt{AMQSGD}, and also use Markovian compressors with popular \texttt{DIANA}~\cite{mishchenko2019distributed} algorithm. 
In all of our experiments, we do not utilize the steps of the optimizer, but rather the information that is transmitted by each worker at the current timestamp $t$. This implies that there are $n$ workers, with each worker sending $m$ coordinates at each iteration of the optimization step. Consequently, the $x$-axis displays numbers of the form $m n \cdot 1, mn \cdot 2, \dots, m n \cdot t, \dots , m n \cdot T$. This allows us to understand the performance of compressors with varying values of $m$ and $n$.

We use convex logistic regression loss with a regularization term $\lambda=0.05$. Each dataset is split horizontally (by rows) equally between $N=10$ clients. The feature dimension is denoted as $d$ in the figures, varying from hundreds to almost a thousand between datasets. The underlying sparsification compressors in Rand-10\% for all logistic regression experiments. Learning rate initial value and decay rate are fine-tuned for each problem and compressor. Markovian-specific parameters such as history size $K$, forgetting rate $b$ are taken from theory or a reasonable default.

\begin{figure}[H]
    \centering
    \begin{tabular}{cccc}
        \includegraphics[width=0.2\linewidth]{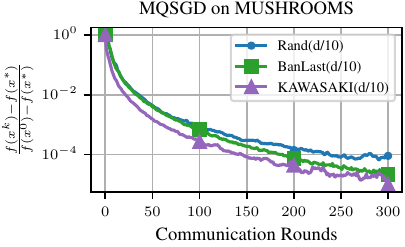} & 
        \includegraphics[width=0.21\linewidth]{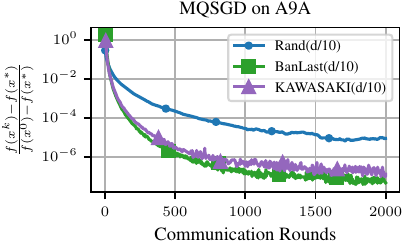} & \includegraphics[width=0.21\linewidth]{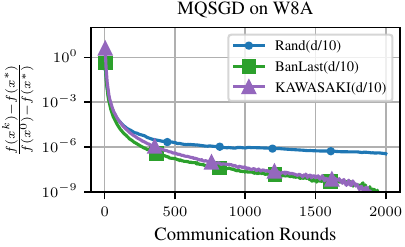} & \includegraphics[width=0.21\linewidth]{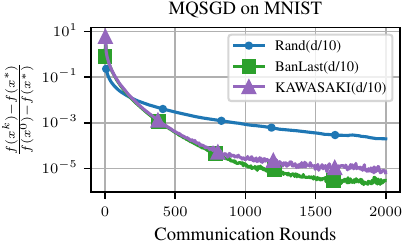} \\
        \includegraphics[width=0.21\linewidth]{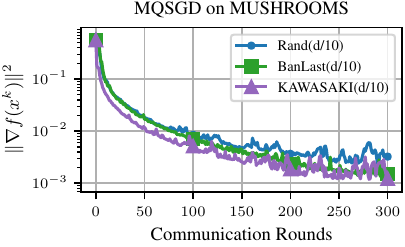} & 
        \includegraphics[width=0.21\linewidth]{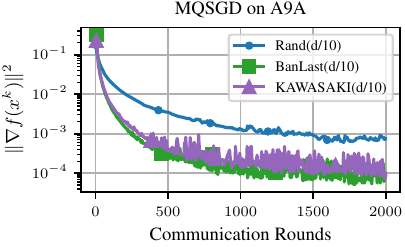} & \includegraphics[width=0.21\linewidth]{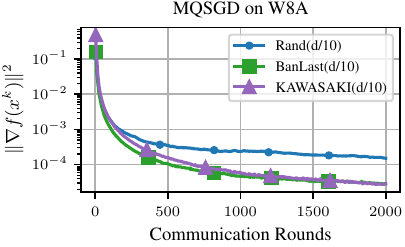} & \includegraphics[width=0.21\linewidth]{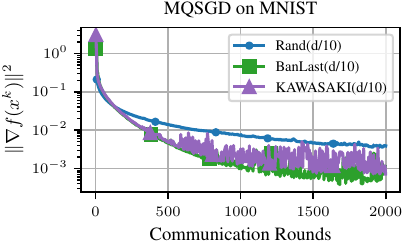} \\
    \end{tabular}
    \caption{\texttt{MQSGD} LIBSVM logistic regression experiments.}
    \label{figure:exp_gd_images}
\end{figure}

\begin{figure}[H]
    \centering
    \begin{tabular}{cccc}
        \includegraphics[width=0.2\linewidth]{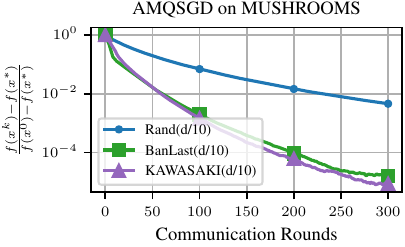} & 
        \includegraphics[width=0.21\linewidth]{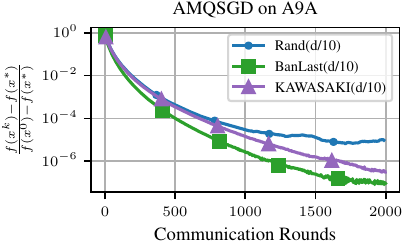} & \includegraphics[width=0.21\linewidth]{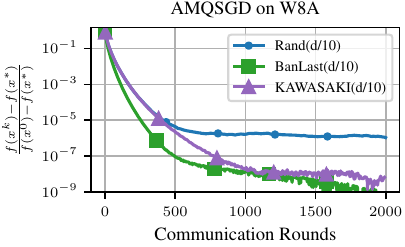} & \includegraphics[width=0.21\linewidth]{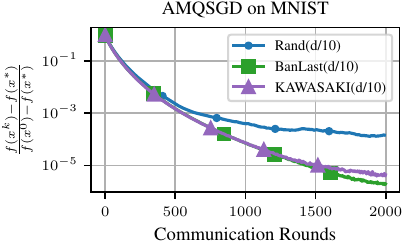} \\
        \includegraphics[width=0.21\linewidth]{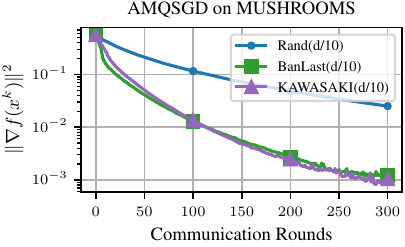} & 
        \includegraphics[width=0.21\linewidth]{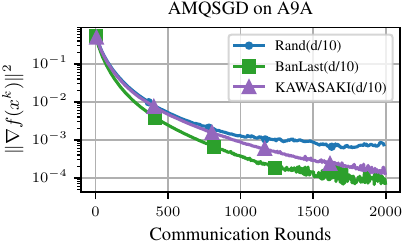} & \includegraphics[width=0.21\linewidth]{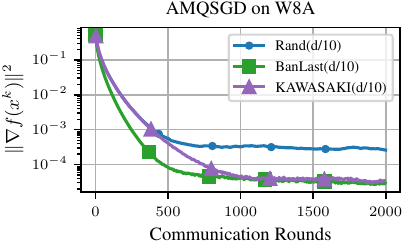} & \includegraphics[width=0.21\linewidth]{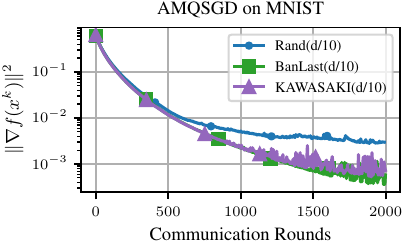} \\
    \end{tabular}
    \caption{\texttt{AMQSGD} LIBSVM logistic regression experiments.}
    \label{figure:exp_accgd_images}
\end{figure}

\begin{figure}[H]
    \centering
    \begin{tabular}{cccc}
         \includegraphics[width=0.21\linewidth]{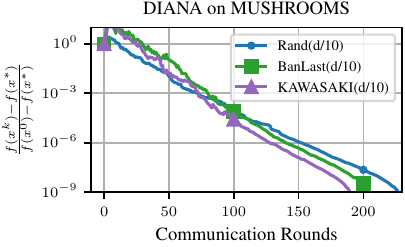} & 
        \includegraphics[width=0.21\linewidth]{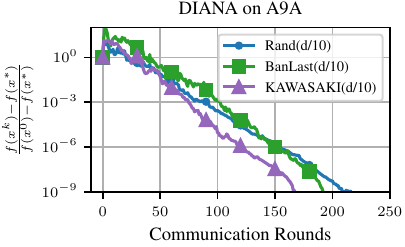} & \includegraphics[width=0.21\linewidth]{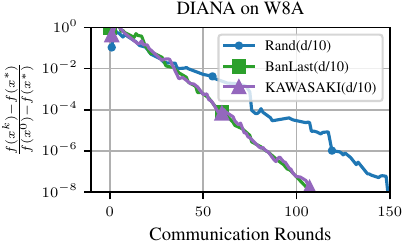} & \includegraphics[width=0.21\linewidth]{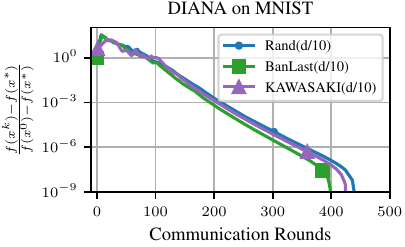} \\
        \includegraphics[width=0.21\linewidth]{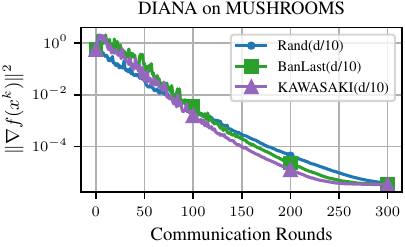} & 
        \includegraphics[width=0.21\linewidth]{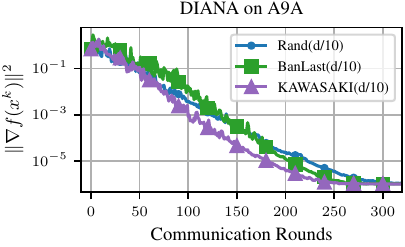} & \includegraphics[width=0.21\linewidth]{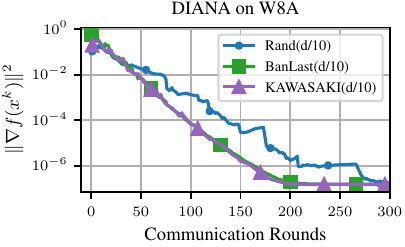} & \includegraphics[width=0.21\linewidth]{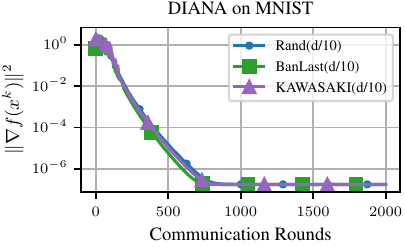} \\
    \end{tabular}
    \caption{\texttt{DIANA} LIBSVM logistic regression experiments.}
    \label{figure:exp_diana_images}
\end{figure}
Figures~\ref{figure:exp_gd_images}, ~\ref{figure:exp_accgd_images} and ~\ref{figure:exp_diana_images} present relative distance to the optimum and gradient norm for the best runs on \texttt{MQSGD}, \texttt{AMQSGD} and \texttt{DIANA}, respectively. We observe that Markovian compressors consistently outperform the Rand-10\% baseline in all scenarios, as the diverging trend can be seen. Only in some experiments with DIANA (MNIST) the advantage is negligible although present. We also observe that simpler and computational-effective \texttt{BanLast} compressor is often enough to achieve substantial convergence improvement. Notably, fine-tuned hyperparameters are similar across datasets and algorithms: for example, \texttt{BanLast} tends to perform best with largest possible values of history size $K$, and \texttt{KAWASAKI} forgetting rate $b$ is large. Notice that \texttt{BanLast} compressor with largest $K$ turns into round-robin compressor with (almost) no stochasticity in coordinates choice.

\subsection{Dependence on size history}
As a part of hyperparameter tuning, we additionally analyze how history size $K$ affects the convergence of Markovian compression-based methods. Figure~\ref{figure:exp_khistory_all} presents dependence of distance to optimum metric on history size for logistic regression experiments. We observe that \texttt{BanLast} performs better around larger values of $K=8$ or $K=9$. In such case for Rand10\% used along with \texttt{BanLast}(9), the compression procedure resembles a permutation: for each 10 iterations, no indices are repeated, and the transmission cycle repeats after that. \texttt{KAWASAKI} history size seems to have periodical spikes and drops, achieving minimum at around $K=25$. However, statistics for \texttt{DIANA} differ drastically, indicating that history size should be adjusted for each problem independently.

\begin{figure}[h]
    \centering
    \begin{tabular}{cccc}
          \includegraphics[width=0.21\linewidth]{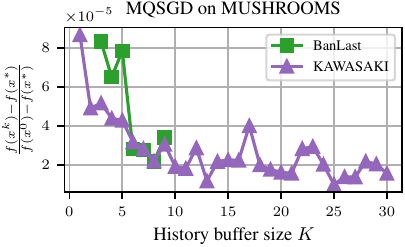} & 
        \includegraphics[width=0.21\linewidth]{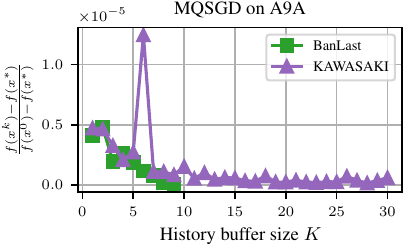} & \includegraphics[width=0.21\linewidth]{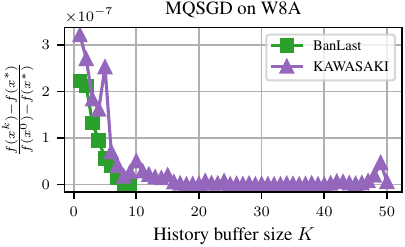} & \includegraphics[width=0.21\linewidth]{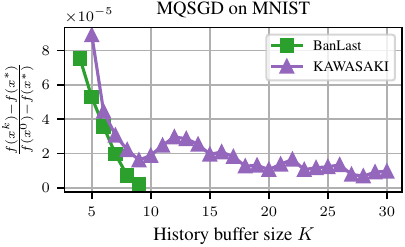} \\
       \includegraphics[width=0.21\linewidth]{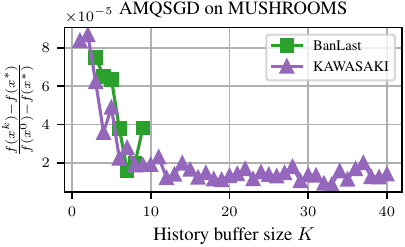} & 
        \includegraphics[width=0.21\linewidth]{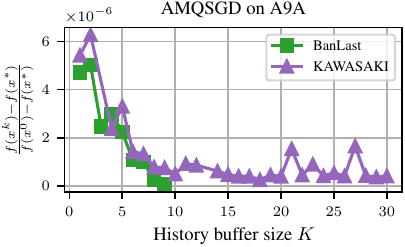} & \includegraphics[width=0.21\linewidth]{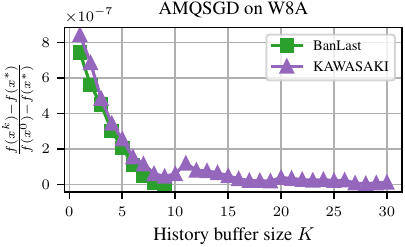} & \includegraphics[width=0.21\linewidth]{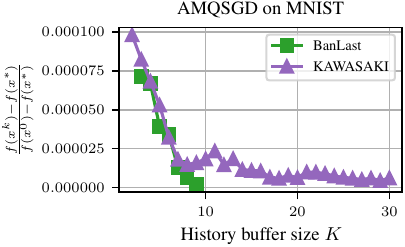} \\
        \includegraphics[width=0.21\linewidth]{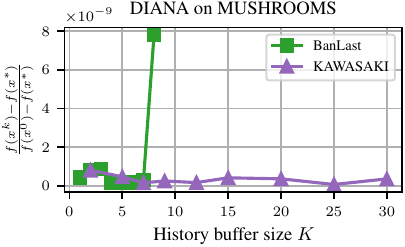} & 
        \includegraphics[width=0.21\linewidth]{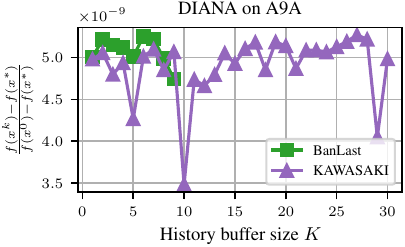} & \includegraphics[width=0.21\linewidth]{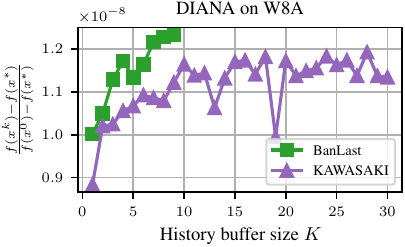} & \includegraphics[width=0.21\linewidth]{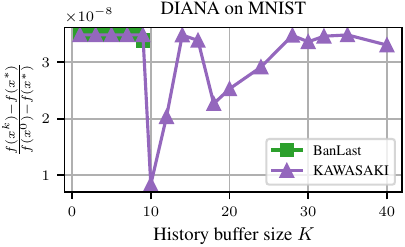} \\
    \end{tabular}
    \caption{Convergence of Markovian-based algorithms on history size $K$}
    \label{figure:exp_khistory_all}
\end{figure}

\subsection{Comparison with Permutation \& Natural Compression}
\label{appendix-comparison-permk-natural}
In this section, we provide empirical comparison of the proposed compressors with other complex compression schemes.

Markovian compressors proposed in the paper compress vector coordinates dependently over optimization epochs. A similar idea of distributed compression is proposed in PermK~\cite{szlendak2021permutation}, where coordinates are arranged between workers at each iteration. Another compressor in the consideration is Natural compression~\cite{horvath2022natural}, an unbiased randomized compressor.

Results of comparison of these compressors on MNIST dataset are presented in Figure~\ref{figure:exp_permk_and_natural}.  The results justify that Markovian compressors tend to converge faster  than the competitors, allowing larger learning rates.

\begin{figure}[h]
    \centering
    \begin{tabular}{ccc}
         \includegraphics[width=0.3\linewidth]{Figures/GD-MNIST-VS-OTHERS-fdist_ratio.pdf} & 
        \includegraphics[width=0.3\linewidth]{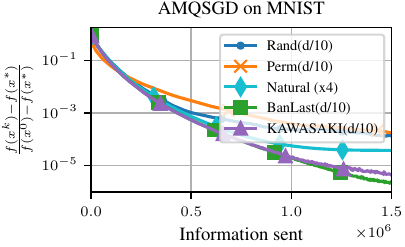} & \includegraphics[width=0.3\linewidth]{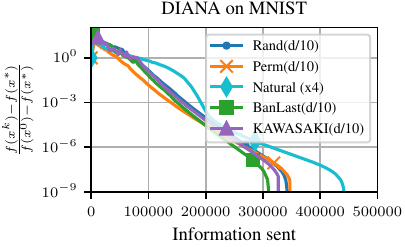}  \\
        \includegraphics[width=0.3\linewidth]{Figures/GD-MNIST-VS-OTHERS-grad_norm.pdf} & 
        \includegraphics[width=0.3\linewidth]{Figures/ACCGD-MNIST-VS-OTHERS-grad_norm.pdf} & \includegraphics[width=0.3\linewidth]{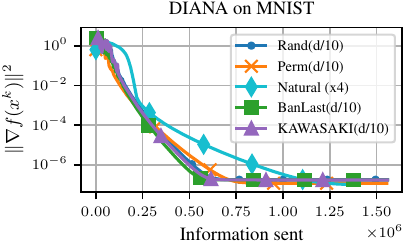} \\
    \end{tabular}
    \caption{Comparison with PermK compressor and Natural compression. PermK compression factor is 10, Natural compression factor is 4. Logistic regression with L2 regularization on MNIST dataset for MQSGD, AMQSGD and DIANA algorithms on $N=5$ clients. Best run is shown after fine-tuning learning rate and its decay. X axis represent amount of information communicated.}
    \label{figure:exp_permk_and_natural}
\end{figure}


\newpage
\subsection{Combination with other compressors}
\label{appendix-comparison-combination}
Although markovian compressors are initially targeted to work with sparsification-based compressors, refining coordinates selection probabilities, they are fully compatible with other compressors afterwards. To illustrate this, and to conduct additional comparison with PermK compressor, we setup experiments combined with Natural Compression~\cite{}. Precisely, we compare RandK+Natural, PermK+Natural, BanLast+Natural and \texttt{KAWASAKI}+Natural compressors on logistic regression on MNIST dataset.

\begin{figure}[h]
    \centering
    \begin{tabular}{ccc}
         \includegraphics[width=0.3\linewidth]{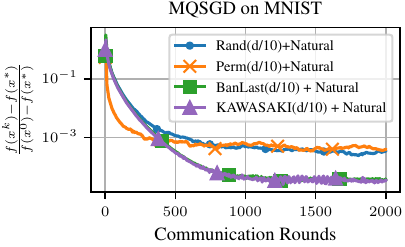} & 
        \includegraphics[width=0.3\linewidth]{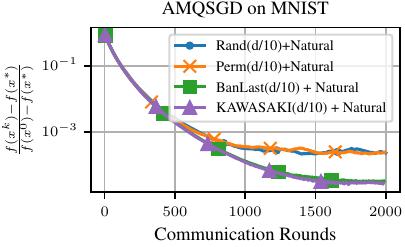} & \includegraphics[width=0.3\linewidth]{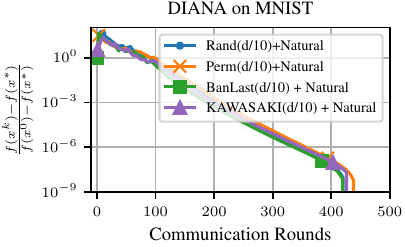}  \\
        \includegraphics[width=0.3\linewidth]{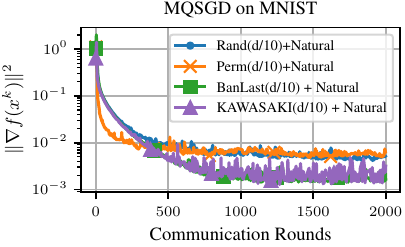} & 
        \includegraphics[width=0.3\linewidth]{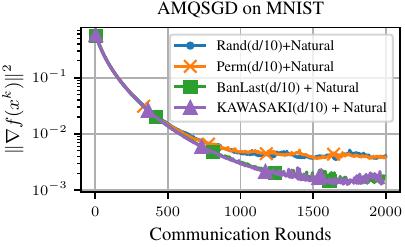} & \includegraphics[width=0.3\linewidth]{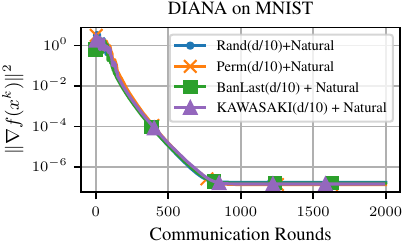} \\
    \end{tabular}
    \caption{Experiments with Natural compression, MNIST logistic regression experiments.}
    \label{figure:exp_natural}
\end{figure}

Figure~\ref{figure:exp_natural} shows the results of the combination of the mentioned sparsification compressors with natural compression.

\subsection{Neural Networks Experiments: Data Parallelism Case}
\label{appx:experiments_nn_dp}
To adopt Markovian compression to a more complex task, we perform image classification on CIFAR-10~\cite{krizhevsky2009cifar} with Resnet-18~\cite{he2016deepresnset} convolutional neural network. We split the training set of size $50,000$ equally between $N=5$ clients. We use SGD optimizer with momentum $0.9$ and weight decay $5 \cdot 10^{-4}$. Hyperparameters such as batch size and learning rate are fine-tuned. Markovian compresors hyperparameters, such as history size $K$, forgetting rate $b$, normalization function are initialized with theoretical optimum or a default. Experiments are conducted with several sparsification compressors, such as Rand-5\%, Rand-7\%, and Rand-10\%, with number of epochs adjusted for each case.

Figures~\ref{figure:exp_nn_rand5_all},~\ref{figure:exp_nn_rand7_all} and ~\ref{figure:exp_nn_rand10_all} present train loss, gradient norm and test accuracy for each baseline method and Markovian compressors for Rand-5\%, Rand-7\% and Rand-10\% scenarios, respectively. Summary on best test accuracy is presented in Table~\ref{table:exp_nn_test_acc}, and extended numerical results for Rand-5\% compressor were presented in main experiments Table~\ref{table:exp_nn_res}. We observe that in such complex, batched optimization problem only \texttt{KAWASAKI} obtains a substantial convergence improvement, as opposed to simpler logistic regression. . In terms of achieved test set accuracy, methods differ significantly only on higher compression rates like Rand-5\%. This may imply that Markovian compression tolerates stronger compression, which is useful in practice. To summarize, Markovian compressors can be successfully applied in neural networks training, with \texttt{KAWASAKI} compressor significantly improving convergence.

Finally, we also conduct the comparison with Permutation and Natural compression, both independently and in combination. Figure~\ref{figure:exp_natural_nn} shows learning curves for training with $N=20$ clients. \texttt{KAWASAKI} compressor appears to have best convergence in both independently and in combination with Natural compression againt Permutation compressor.

\begin{figure}[H]
    \centering
        \includegraphics[width=0.3\linewidth]{Figures/CIFAR-RESNET-RAND5-train_loss.pdf}
        \includegraphics[width=0.3\linewidth]{Figures/CIFAR-RESNET-RAND5-grad_norm.pdf}
        \includegraphics[width=0.3\linewidth]{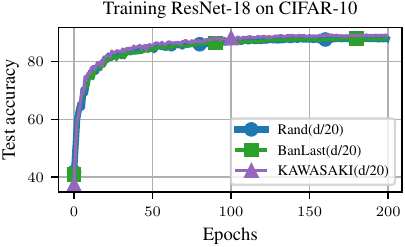}
    \caption{Resnet-18 on CIFAR-10 training results for Rand-5\% sparsification.}
    \label{figure:exp_nn_rand5_all}
\end{figure}

\begin{figure}[H]
    \centering
        \includegraphics[width=0.3\linewidth]{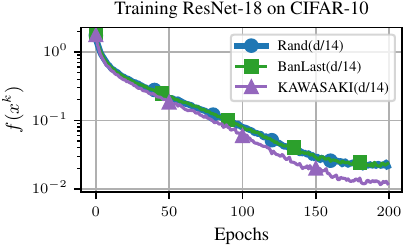}
        \includegraphics[width=0.3\linewidth]{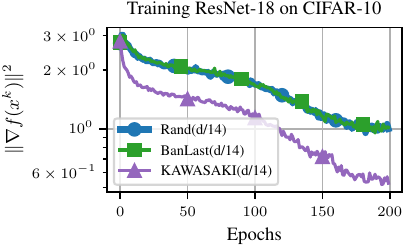}
        \includegraphics[width=0.3\linewidth]{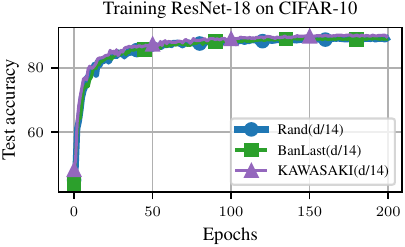}
    \caption{Resnet-18 on CIFAR-10 training results for Rand-7\% sparsification.}
    \label{figure:exp_nn_rand7_all}
\end{figure}

\begin{figure}[H]
    \centering
        \includegraphics[width=0.3\linewidth]{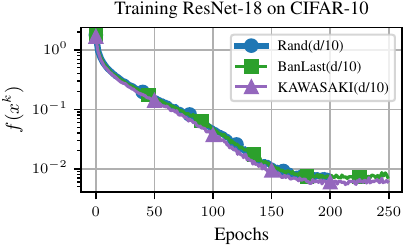}
        \includegraphics[width=0.3\linewidth]{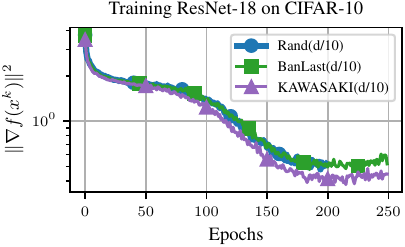}
        \includegraphics[width=0.3\linewidth]{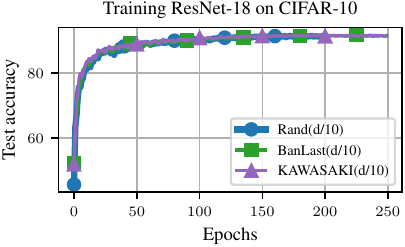}
    \caption{Resnet-18 on CIFAR-10 training results for Rand-10\% sparsification.}
    \label{figure:exp_nn_rand10_all}
\end{figure}

\begin{figure}[H]
    \centering
        \includegraphics[width=0.45\linewidth]{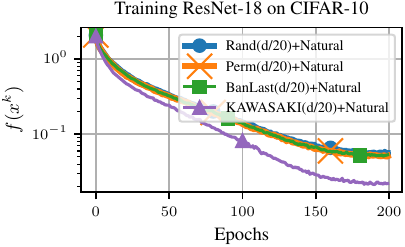}
        \includegraphics[width=0.45\linewidth]{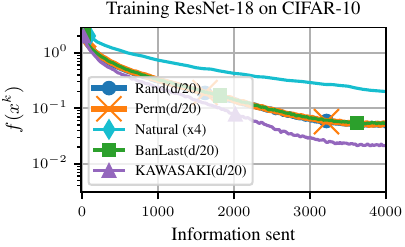}
    \caption{Comparison with other compressors on Resnet-18 training on CIFAR-10 dataset for Rand-5\% sparsification on $N=20$ clients. Natural compression factor is 4. Left figure is sequential combination with Natural compression. Right figure is comparison against PermK and Natural compressors independently, with information sent on x-axis.}
    \label{figure:exp_natural_nn}
\end{figure}


\begin{table}[htbp]
    \centering
    \begin{tabular}{lccc}
    \mbox{} & Rand-K\% & \texttt{Banlast} & \texttt{KAWASAKI} \\ \hline
    Rand-5\% & 88.03 & 88.1 & \textbf{89.27} \\ \hline
    Rand-7\% & 89.31 & 89.38 & \textbf{90.28} \\ \hline
    Rand-10\% & \textbf{91.46} & \textbf{91.72} & \textbf{91.78} \\ \hline
    \end{tabular}
    \caption{\centering Best test accuracy \% of training ResNet-18 on CIFAR-10 with different compressors}
    \label{table:exp_nn_test_acc}
\end{table}

\subsection{Neural Networks Experiments: Model Parallelism Case}
\label{appx:experiments_nn_mp}
As opposed to data-parallel setting, model parallelism is paradigm which splits the model (typically a deep neural network) to a pipeline of layers between workers. Such distributed scenario is especially relevant for large language models (LLM), which consist of billions of trainable parameters. As communication is a typical bottleneck in such systems~\cite{diskin2021distributed}, various compression techniques are applied to layer activations and their respective gradients that are transferred between adjacent pipeline workers. Such techniques include quantization and  sparsification~\cite{dettmers2022llm,bian2023doesactivationcompression}, as well as low-rank compression~\cite{song2023optimus} techniques.

We perform training of Resnet-18~\cite{he2016deepresnset} convolutional neural network on CIFAR-10 dataset~\cite{krizhevsky2009cifar}. We split the ResNet onto 4 workers by resnet blocks, simulated on a single device with compression of activations and their respective gradients in the places of communication.
We apply Markovian compressors only to gradients in model-parallel setup, using same RandK compression for both activations and gradients independently for each compression block.

\begin{table}[htbp]
    \begin{center}
    \begin{tabular}{l l l}
        \toprule
        \textbf{Compressor} & \textbf{Compression ON} & \textbf{Compression OFF}\\
        \midrule
        No compression & 92.8 & 92.8 \\
        Rand10\% & 84.6 & 86.1 \\
        BanLastK+Rand10\% & 85.2 & 86.4 \\
        KAWASAKI(simplex projection)+Rand10\% & 84.5 & 85.0 \\
        KAWASAKI(normalize)+Rand10\% & 85.2 & \textbf{86.8} \\
        KAWASAKI(softmax)+Rand10\% & 85.3 & \textbf{87.3} \\
        \bottomrule
        \end{tabular}
    \end{center}
    \caption{Best test accuracy \% for model parallelism experiments with Resnet-18 classification of CIFAR-10}
    \label{table:model-parallel}
\end{table}

Table ~\ref{table:model-parallel} presents best test set accuracy achieved for training with different compressors. While compression indeed decreases accuracy for Rand-10\%, application of Markov compressors favors the final test accuracy on a whole one percent. Note that compression is not applied during inference, only on training phase. This case illustrates potential of Markov compressors beyond data-parallelism setup considered in theory. In practical training of large neural networks, where both data-parallelism and model-parallelism are often applied simultaneously, Markov compressors could also be useful, as per shown efficiency on both these setups in separate.

\subsection{Fine-tuning DeBERTaV3-base for GLUE benchmark}
    \label{appx:exp_nn_deberta}

    In this series of experiments, we examine a distributed approach to fine-tuning language models using LoRA \citep{hu2021loralowrankadaptationlarge}. This method is based on freezing the model weights that are pre-trained on a large dataset, and add a low rank adapter with matrices $A \in \R^{n \times r}$ and $ B \in \R^{r \times m}$ to some selected layers $W_{\text{old}} \in \R^{n \times m}$ of this model, such that $W_{\text{new}} = W_{\text{old}} + A \cdot B$. Since in practice the parameter $r$ is chosen to be much smaller than $n$ and $m$, the new model has much fewer trainable parameters and can be efficiently trained on downsteram tasks. 

    In our experiments, we apply LoRA adapters with fixed rank $r = 8$ to the attention layers of the DeBERTaV3-base model \citep{he2021debertadecodingenhancedbertdisentangled}. The downsteram task is the classical GLUE benchmark for natural language understanding \citep{wang2019gluemultitaskbenchmarkanalysis}. We consider only random sparsification compressors (Definition \ref{def:rand_spars}) with $25\%$ compression rate, due to the large computational cost of this experiment. Figure~\ref{figure:peft} shows learning curves for training with $N=10$ clients. Our Markovian compressors 
    appears to have best convergence against independent Rand$m$ compressor.

\begin{figure}[H]
    \centering
        \includegraphics[width=0.32\linewidth]{Figures/DEBERTA-SST2-train-loss.pdf}
        \includegraphics[width=0.32\linewidth]{Figures/DEBERTA-QNLI-train-loss.pdf}
        \includegraphics[width=0.32\linewidth]{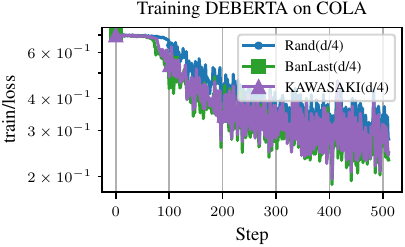}
    \caption{Comparison with other compressors on fine-tuning task on GLUE benchmark on $N=10$ clients. We performed experiments on SST2, QNLI and COLA tasks.}
    \label{figure:peft}
\end{figure}



\end{appendixpart}
\end{document}